\documentclass[final,hidelinks,onefignum,onetabnum]{siamart220329}

\ifpdf
\hypersetup{
	pdftitle={Space Mapping for PDE Constrained Shape Optimization},
	pdfauthor={S. Blauth}
}
\fi

\usepackage{lipsum}
\usepackage{amsfonts}
\usepackage{graphicx}
\usepackage{epstopdf}
\usepackage{algorithmic}
\ifpdf
\DeclareGraphicsExtensions{.eps,.pdf,.png,.jpg}
\else
\DeclareGraphicsExtensions{.eps}
\fi

\newsiamremark{remark}{Remark}
\newsiamremark{example}{Example}
\newsiamremark{hypothesis}{Hypothesis}
\crefname{hypothesis}{Hypothesis}{Hypotheses}
\crefname{algocf}{Algorithm}{Algorithms}
\newsiamthm{claim}{Claim}

\headers{Space Mapping for Shape Optimization}{S. Blauth}

\title{Space Mapping for PDE Constrained Shape Optimization\thanks{Submitted to the editors on \today.
		\funding{This work was funded by the ERDF - European Regional Development Fund for Rhineland-Palatinate as part of the reaction of the European Union on the COVID-19 pandemic (REACT-EU) in the project \textit{ENERDIG} under grant no.~84008454.}}}

\author{Sebastian Blauth\thanks{Fraunhofer ITWM, Kaiserslautern, Germany 
		(\email{sebastian.blauth@itwm.fraunhofer.de}.}}

\usepackage{amsopn}

\usepackage[algo2e, linesnumbered, noline, noend, ruled]{algorithm2e}
\usepackage{amsmath}
\usepackage{amssymb}
\usepackage[english]{babel}
\usepackage{caption}
\usepackage{braket}
\usepackage{float}
\usepackage[T1]{fontenc}
\usepackage{graphicx}
\usepackage[utf8]{inputenc}
\usepackage{mathtools}
\usepackage{nicefrac}
\usepackage[defaultlines=2,all]{nowidow}
\usepackage[per-mode=symbol, group-digits=integer]{siunitx}
\usepackage{enumitem}
\usepackage{subcaption}
\usepackage{stmaryrd}
\usepackage{textcomp}

\usepackage{cite}
\usepackage{tikz}

\definecolor{myblue}{RGB}{0, 107, 164}
\definecolor{myorange}{RGB}{255, 128, 14}
\definecolor{grayone}{RGB}{171, 171, 171}

	\newcommand{\qe}[1]{``#1''}

	\DeclareMathOperator*{\argmin}{argmin}
	
	\newcommand{\R}{\mathbb{R}}
	
	\newcommand{\trace}[1]{\mathrm{tr}\left( #1 \right)}
	\newcommand{\transposed}{^\top}
	
	\newcommand{\norm}[2]{\left\lvert \left\lvert #1 \right\rvert \right\rvert_{#2}}
	\newcommand{\abs}[1]{\left\lvert #1 \right\rvert}
	\newcommand{\dmeas}[1]{\ \mathrm{d}#1}
	
	\newcommand{\grad}{\nabla}
	\renewcommand{\div}[1]{\mathrm{div} \left( #1 \right)}
	\newcommand{\integral}[1]{\int_{#1}}
	\newcommand{\dual}[3]{\left\langle #1 , #2 \right\rangle_{#3}}
	\newcommand{\inner}[3]{\left( #1 , #2 \right)_{#3}}

	\newcommand{\tgrad}{\grad_\Gamma}

	\newcommand{\admissiblegeom}{\mathcal{A}}
	\newcommand{\costfunctional}{\mathcal{J}}
	\newcommand{\reducedcostfunctional}{J}
	\newcommand{\holdall}{D}
	\newcommand{\vectorfield}{\mathcal{V}}
	\newcommand{\flow}{\Phi}
	
	\newcommand{\smallspace}{\thinspace}
	
	\newcommand{\interfacecoefficient}{\alpha}
	
	\newcommand{\subdes}{_\mathrm{des}}
	\newcommand{\subin}{_\mathrm{in}}
	\newcommand{\subout}{_\mathrm{out}}
	\newcommand{\subfine}{_\mathrm{f}}
	\newcommand{\subcoarse}{_\mathrm{c}}
	\newcommand{\fineresponse}{f}
	\newcommand{\coarseresponse}{c}
	\newcommand{\spacemappingfunction}{s}
	\newcommand{\misalignmentfunction}{r}
	
	\newcommand{\jump}[2]{\left\llbracket #1 \right\rrbracket_{#2}}

	\newcommand{\shapedistro}{\gamma}
	\newcommand{\gradientdefo}{\mathcal{G}}
	
	\newcommand{\scalarproduct}{g^S}
	\newcommand{\stekpoinc}{S^p}
	\newcommand{\rieshapegrad}{\mathrm{grad}\reducedcostfunctional}
	\newcommand{\geodesic}{\gamma}
	\newcommand{\retraction}{\mathrm{R}}
	\newcommand{\numretraction}{\tilde{\mathrm{R}}}
	\newcommand{\vectortransport}{\mathcal{T}}
	\newcommand{\numtransport}{\tilde{\mathcal{T}}}

	\newcommand{\shapemanifold}{B_e}
	\newcommand{\tangentspace}{T}

	\numberwithin{equation}{section}
	
	\newcommand{\normal}{n}
	
	\newcommand*\closure[1]{%
		\hbox{%
			\vbox{%
				\hrule height 0.5pt % The actual bar
				\kern0.5ex%         % Distance between bar and symbol
				\hbox{%
					\ensuremath{#1}%
				}%
			}%
		}%
	} 
	
\SetCommentSty{mycommfont}
\SetKwComment{Comment}{\# }{}

\begin{document}

{\noindent\footnotesize This is a post-peer-review, pre-copyedit version of an article published in SIAM Journal on Optimization. The final version is available online at \url{https://doi.org/10.1137/22M1515665}.
}

\maketitle

\begin{abstract}
	The space mapping technique is used to efficiently solve complex optimization problems. It combines the accuracy of fine model simulations with the speed of coarse model optimizations to approximate the solution of the fine model optimization problem. In this paper, we propose novel space mapping methods for solving shape optimization problems constrained by partial differential equations (PDEs). We present the methods in a Riemannian setting based on Steklov-Poincar\'e-type metrics and discuss their numerical discretization and implementation. We investigate the numerical performance of the space mapping methods on several model problems. Our numerical results highlight the methods' great efficiency for solving complex shape optimization problems.
\end{abstract}

% REQUIRED

% REQUIRED

\section{Introduction}
\label{sec:introduction}

The aim of the space mapping technique is to reduce the cost of solving complex optimization problems by iteratively correcting a sequence of rougher approximations. It utilizes a model hierarchy consisting of two models: First, a fine (complex, detailed, costly) model, which captures all relevant effects and enables detailed, but costly, numerical simulations. Second, a coarse (cheap, approximate, simple) model, which, e.g., contains simplifications, neglects some effects or uses a coarser discretization, but is evaluated cheaply. The fine model is well-suited for detailed simulations, but is expensive to optimize. The coarse model is substantially easier to optimize, but only gives a rough approximation. Consequently, an optimization with the fine model is infeasible, whereas an optimization with the coarse model may yield inaccurate results. The space mapping technique aims at overcoming these limitations by combining the advantages of both models: It utilizes the speed of the coarse model optimization to generate a sequence of iterates, which are corrected by the high accuracy of the fine model simulations. Hence, the expensive direct optimization of the fine model is approximated by repeatedly optimizing and adapting the coarse model.

Conceived in 1994 in the context of microwave filter design \cite{Bandler1994Space}, the space mapping technique has received lots of attention, particularly in the engineering community, where it has been used, e.g., in connection with electromagnetic applications \cite{Bandler1994Space,Bandler1995Electromagnetic,Bandler2004space_mapping}, car design \cite{Redhe2006multipoint}, and structural optimization \cite{Leary2001constraint}. Moreover, space mapping techniques have been applied in the literature for solving optimization problems constrained by partial differential equations (PDEs) in the context of optimal control \cite{Marheineke2012Model,Hintermueller2005Space,Weisen2021Space}.

The aim of PDE constrained shape optimization is to optimize the shape of a system to improve its performance. To do so, a cost functional which measures the performance or quality of the system is defined and subsequently optimized w.r.t.\ the shape of the system in order to optimize the latter. Shape optimization has many interesting and practically relevant applications: It is, e.g., used for the optimization of electromagnetic devices \cite{Gangl2015Shape}, aircraft \cite{Schmidt2013Three}, polymer spin packs \cite{Hohmann2019Shape}, microchannel heat exchangers \cite{Blauth2021Model,Blauth2020Shape}, and in the context of electric impedance tomography \cite{Hintermueller2015Shape,Laurain2016Distributed}. For the efficient solution of shape optimization problems, shape calculus is a valuable tool. It is used to calculate the sensitivity of a shape functional w.r.t. infinitesimal variations of the shape and enables the use of gradient-based optimization algorithms \cite{Blauth2021Nonlinear, Schulz2016Efficient}. Recently, shape optimization based on shape calculus attracted lots of attention, particularly regarding the development of efficient algorithms. This can be, e.g., seen in \cite{Schulz2014Riemannian,Schulz2016Efficient}, where (quasi-)Newton methods for shape optimization are introduced, in \cite{Mueller2021novel,Deckelnick2022novel}, where the $W^{1,\infty}$-topology is used to compute mesh deformations, and in \cite{Blauth2021Nonlinear,Blauth2022Shape}, where nonlinear conjugate gradient methods for shape optimization are proposed. 

In this paper, we extend the space mapping technique to the field of shape optimization. We detail an aggressive space mapping method with the help of the Steklov-Poincar\'e-type metrics from \cite{Schulz2016Efficient}. We investigate this method numerically for a shape identification problem constrained by a semi-linear transmission problem and a shape optimization problem of uniform flow distribution constrained by the Navier-Stokes equations.  Moreover, we show the practical applicability of our approach by coupling our method with commercial solvers for fine model simulations. Our numerical results show the methods' great efficiency and fast convergence. 

In the literature, space mapping techniques have already been applied to the field of shape optimization, e.g., in \cite{Toman2019Blade,Jonsson2015Shape,Leifsson2013Aerodynamic}. However, in the previous literature, the shapes under consideration were always parametrized and, thus, the problems were reduced to finite-dimensional optimal control problems. In this paper, we consider shape optimization based on shape calculus and, therefore, investigate infinite-dimensional shape optimization problems without parametrizations of the shape. To the best of our knowledge, such a space mapping technique for PDE constrained shape optimization based on shape calculus has not yet been investigated in the literature. 

This paper is structured as follows. In \cref{sec:preliminaries}, we briefly recall the space mapping technique, basic results from shape calculus, and the Steklov-Poincar\'e-type metrics from \cite{Schulz2016Efficient}. We propose our space mapping technique for shape optimization in \cref{sec:space_mapping} and present our numerical implementation, which is publicly available in our open-source software cashocs \cite{Blauth2021cashocs,Blauth2022Software}. In \cref{sec:numerics}, we investigate the proposed method numerically and demonstrate its efficiency for solving complex shape optimization problems.

\section{Preliminaries}
\label{sec:preliminaries}

In this section, we recall the main ideas of the space mapping technique, present fundamental concepts from shape calculus and Riemannian geometry as well as the Steklov-Poincar\'e-type metric from \cite{Schulz2016Efficient}.

\subsection{The Space Mapping Technique}
\label{ssec:space_mapping}

In this section, we recall the basic ideas of the space mapping technique, following the presentation of \cite{Echeverria2005Space, Bakr2001introduction, Marheineke2012Model}.

Assume that we want to optimize some system for which we have a model hierarchy consisting of a fine and a coarse model. The fine model is accurate, but expensive to evaluate, whereas the coarse model is less accurate, but cheap to evaluate. Additionally, we assume that the optimization of the coarse model is possible and comparatively cheap, whereas the optimization of the fine model is very expensive or even impossible due to the high cost of the fine model simulations.

Let $\fineresponse$ be the fine model response, i.e., $\fineresponse\colon X \to Y$, where $X$ is the control space and $Y$ is some Banach space. Both $X$ and $Y$ are Banach spaces which depend on the application under consideration (some examples for this can be found in \cref{sec:numerics}. We can write the fine model optimization problem as
\begin{equation}
	\label{eq:opt_fine}
	\min_{x\in X} \reducedcostfunctional(\fineresponse(x)),
\end{equation}
i.e., we want to minimize a cost functional $\reducedcostfunctional\colon Y\to \R$ which depends on the fine model response $\fineresponse(x)$. 

Problem \cref{eq:opt_fine} is, in general, too complex and too expensive to be solved directly. For this reason, we consider the following: Let $\coarseresponse$ be the coarse model response, i.e., $\coarseresponse\colon X \to Y$. For the simplicity of presentation, we assume that the control spaces as well as the co-domains of the model responses coincide and we refer the reader to, e.g, \cite{Echeverria2005Space,Hintermueller2005Space} for a more general setting. We approximate the fine model optimization problem \cref{eq:opt_fine} by the following coarse model optimization problem
\begin{equation}
	\label{eq:opt_coarse}
	\min_{x\in X} \reducedcostfunctional(\coarseresponse(x)).
\end{equation}
To obtain well-posed problems, we assume that both \cref{eq:opt_fine} and \cref{eq:opt_coarse} have a unique minimizer. Consequently, we can introduce
\begin{equation}
	\label{eq:definition_minimizers}
	x\subfine^* := \argmin_{x \in X}\smallspace \reducedcostfunctional(\fineresponse(x)) \qquad \text{ and } \qquad x\subcoarse^* := \argmin_{x \in X}\smallspace \reducedcostfunctional(\coarseresponse(x)).
%	\begin{aligned}
%		u\subfine^* &:= \argmin_{u \in U}\smallspace \reducedcostfunctional(\fineresponse(u)), \\
%		u\subcoarse^* &:= \argmin_{u \in U}\smallspace \reducedcostfunctional(\coarseresponse(u)),
%	\end{aligned}
\end{equation}

As described earlier, we assume that solving the coarse model optimization problem \cref{eq:opt_coarse} is much cheaper than solving the fine model optimization problem \cref{eq:opt_fine}. Therefore, we do not attempt to solve \cref{eq:opt_fine} directly, but try to approximately solve it utilizing the coarse model. A naive approach for approximating $x\subfine^*$ is to solve \cref{eq:opt_coarse} and use its minimizer $x\subcoarse^*$ as approximation of $x\subfine^*$. The space mapping technique extends and generalizes this approach. To do so, we first introduce the so-called space mapping function $\spacemappingfunction\colon X \to X;\ x\subfine \mapsto \spacemappingfunction(x\subfine)$, where
\begin{equation}
	\label{eq:def_mapping_function}
	\spacemappingfunction(x\subfine) := \argmin_{x\subcoarse \in X}\smallspace \misalignmentfunction\left(\coarseresponse(x\subcoarse), \fineresponse(x\subfine)\right),
\end{equation}
with some misalignment function $\misalignmentfunction\colon Y \times Y \to \R$, which is used to measure the discrepancy between the fine and coarse model responses. To get a well-defined space mapping function, we assume that the minimization problem in \cref{eq:def_mapping_function} is well-posed, i.e., that it has a unique minimizer for all $x\subfine \in X$. Further, we assume that the misalignment function is exact in the sense that 
\begin{equation*}
	\argmin_{x\subcoarse\in X}\smallspace \misalignmentfunction\left( \coarseresponse(x\subcoarse), \coarseresponse(z) \right) = z \quad \text{ for all } z\in X.
\end{equation*}
For a given fine model control $x\subfine$, the space mapping function computes the best coarse model control $x\subcoarse$ such that the discrepancy between the fine and coarse model responses $\fineresponse(x\subfine)$ and $\coarseresponse(x\subcoarse)$ is minimized. 

We require that the coarse model is a suitable approximation of the fine model, at least in the vicinity of their respective minimizers, i.e., $\coarseresponse(x\subcoarse^*) \approx \fineresponse(x\subfine^*)$. Under this assumption, we expect that
\begin{equation}
	\label{eq:mapping}
	\spacemappingfunction(x\subfine^*) = \argmin_{x\subcoarse \in X} \misalignmentfunction(\coarseresponse(x\subcoarse), \fineresponse(x\subfine^*)) \approx \argmin_{x\subcoarse \in X} \misalignmentfunction\left( \coarseresponse(x\subcoarse), \coarseresponse(x\subcoarse^*) \right) = x\subcoarse^*,
\end{equation}
due to the exactness of the of the misalignment function (cf.~\cite{Echeverria2005Space}). Motivated by this, the fundamental idea of the space mapping technique is to solve the equation
%the space mapping technique aims at solving the equation
\begin{equation}
	\label{eq:sm_root}
	\spacemappingfunction(x\subfine^*) - x\subcoarse^* = 0.
\end{equation}

\begin{remark}
	An appropriate choice of the coarse model is crucial for the success of the space mapping technique. Note that the coarse model does not have to represent the same PDE as the fine model. Instead, the coarse model can be given by a similar PDE, which is easier to solve numerically. Here, \qe{similar} means that the PDEs should model similar phenomena. In case the coarse model is not a sufficiently good approximation of the fine model, the space mapping technique cannot be expected to work properly. For example, neither the Navier-Stokes equations for fluid flow nor the equations of linear elasticity are, of course, appropriate coarse models for a radiative heat transfer problem, as the models describe completely different physical phenomena.
	
	For complex optimization problems arising from industrial applications, it is often straightforward to derive a coarse model by using, e.g., a coarser discretization, a linearization, neglecting some physical effects, weakening nonlinear effects, etc.. For example, when solving an optimization problem constrained by the Navier-Stokes equations with a high Reynolds number, one could, e.g., use the Navier-Stokes equations with a sufficiently low Reynolds number (weakening the nonlinearity), a linearized Navier-Stokes system, or even the linear Stokes system (which completely neglects inertial effects) as coarse models. Some examples for using a coarse model are shown in the numerical experiments in \cref{sec:numerics}. Additionally, for an overview of such model hierarchies for space mapping in the context of transport processes we refer the reader to \cite{Marheineke2012Model}.	
\end{remark}

\begin{remark}
	In theory, the space mapping approach requires that both the fine and the coarse model have a unique global minimizer (see \cref{eq:definition_minimizers}). Additionally, one has to be able to compute the global minimizer of the coarse model efficiently to apply the space mapping technique. However, for our application of PDE constrained shape optimization, one can typically only compute local minimizers due to the fact that such problems are highly nonlinear and non-convex. Nevertheless, the space mapping approach can also be applied formally in this setting and still yields convincing results, as shown in \cref{sec:numerics}.
\end{remark}

In the following we briefly present the so-called aggressive space mapping (ASM) method \cite{Bandler1995Electromagnetic}, which is one of the most popular space mapping methods. Let us assume that the control space $X$ is a Hilbert space. The ASM method uses Broyden's method to solve equation \cref{eq:sm_root} for $x\subfine^*$. This can be interpreted as successively approximating the space mapping function by $\spacemappingfunction^k(x) = \spacemappingfunction(x\subfine^k) + B^k(x - x\subfine^k)$, where $B^k$ is an approximate Jacobian of $\spacemappingfunction$. Therefore, the ASM method solves \cref{eq:sm_root} with the iteration $\spacemappingfunction^k(x\subfine^{k+1}) = x\subcoarse^*$, such that we get
\begin{equation*}
	B^k h^k =  - \left(\spacemappingfunction(x\subfine^k) - x\subcoarse^*\right), \qquad x\subfine^{k+1} = x\subfine^k + h^k.
\end{equation*}
This is typically supplemented with the initial guess $x\subfine^0 = x\subcoarse^*$ and $B^0 = I$. Further, the update formula for the approximate Jacobian $B^k$ is given by
\begin{equation*}
	B^{k+1} v = B^k v + \frac{\spacemappingfunction(x\subfine^{k+1}) - \spacemappingfunction(x\subfine^k) - B^k h^k}{\norm{h^k}{X}^2} \inner{h^k}{v}{X} \quad \text{ for all } v\in X,
\end{equation*}
%or, equivalently,
%\begin{equation*}
%	B^{k+1} = B^k + \frac{\spacemappingfunction(u\subfine^{k+1}) - \spacemappingfunction(u\subfine^{k}) - B^k h^k}{\norm{h^k}{U}^2} \otimes h^k
%\end{equation*}
where $h^k = x\subfine^{k+1} - x\subfine^k$, 
%$\otimes$ denotes the dyadic product, 
and $\inner{a}{b}{H}$ denotes the inner product of $a$ and $b$ in the Hilbert space $H$. The iteration can be terminated once the relative stopping criterion
\begin{equation*}
	\frac{\norm{\spacemappingfunction(x\subfine^k) - x\subcoarse^*}{X}}{\norm{x\subcoarse^*}{X}} \leq \tau
\end{equation*}
is satisfied. 

\begin{remark}
	Note that there exist many other space mapping approaches in the literature, such as primal, dual, and trust region aggressive space mapping. Additionally, there exist alternative methods for solving \cref{eq:sm_root}, e.g., steepest descent, BFGS, or nonlinear conjugate gradient methods, all of which are implemented in our software cashocs \cite{Blauth2021cashocs}. For the sake of brevity, we restrict our attention to the ASM method in this paper and refer the reader to, e.g., \cite{Bandler2004Space,Bakr2001introduction,Echeverria2005Space} and the references therein for an overview of other space mapping approaches.
\end{remark}

\begin{remark}
	One of the main features of the space mapping method is that it is able to solve the fine model optimization problem efficiently while only using simulations of the fine model. No sensitivity or adjoint information of the fine model is required. This is very relevant for practical applications, where very complex models are required to properly model the relevant physical effects. Typically, commercial solvers are used for the simulation of such models as they offer lots of relevant models and sophisticated solvers for such problems. In such a setting, an optimization of the fine model can become prohibitively costly - both due to the cost of simulating the fine model as well as the amount of effort required to implement a solver for the optimization problem, which would involve the calculation of sensitivities or adjoint information of the fine model. However, the space mapping technique enables the optimization of such complex models even in such a setting by only requiring fine model simulations, which can be performed with commercial solvers, whereas all optimizations are performed only w.r.t.\ the coarse model. 
\end{remark}

\subsection{Shape Calculus}
\label{ssec:shape_calculus}

A PDE constrained shape optimization problem can be written in the form
\begin{equation*}
	\min_{\Omega \in \admissiblegeom} \costfunctional(\Omega, u) \quad \text{ s.t. } \quad e(\Omega, u) = 0,
\end{equation*}
where $\costfunctional$ is a cost functional depending on the geometry $\Omega$ and the state variable $u$. For the optimization, the cost functional $\costfunctional$ shall be minimized w.r.t.\ the geometry $\Omega$ over some set of admissible geometries $\admissiblegeom$. Additionally, $e(\Omega, \cdot)\colon U(\Omega) \to V(\Omega)^*$ is an operator which models the state equation, i.e., a PDE constraint, which we consider in the weak form
\begin{equation*}
	\text{Find } u \in U(\Omega) \text{ such that } \quad \dual{e(\Omega, u)}{v}{V(\Omega)^*, V(\Omega)} = 0 \quad \text{ for all } v\in V(\Omega).
\end{equation*}
We assume that the state equation has a unique weak solution $u=u(\Omega) \in U(\Omega)$ for all $\Omega \in \admissiblegeom$ such that $e(\Omega, u(\Omega)) = 0$. Hence, we can define the reduced cost functional $\reducedcostfunctional(\Omega) = \costfunctional(\Omega, u(\Omega))$ and consider the following equivalent unconstrained problem
\begin{equation}
	\label{eq:reduced_optimization_problem}
	\min_{\Omega \in \admissiblegeom} \reducedcostfunctional(\Omega).
\end{equation}
To illustrate this, let us consider an example similar to the one discussed in \cite{Schulz2016Efficient}.

\begin{example}
	\label{ex:transmission}
	We consider a fixed, open, and bounded hold-all domain $\holdall \subset \R^d$ for $d\in \mathbb{N}, d \geq 2$ with Lipschitz boundary $\partial\holdall$. The hold-all domain is partitioned into two disjoint open subsets: an inclusion $\Omega \subset \holdall$ and an exterior part $\Omega\subout := \holdall \setminus \closure{\Omega}$, where $\closure{\Omega}$ denotes the closure of $\Omega$. We denote by $\Gamma = \partial \Omega$ the boundary of $\Omega$, which is assumed to be smooth, and $\normal$ denotes the outward unit normal on $\Gamma$. This setting is illustrated in \cref{fig:model_setup}. 
	
	\begin{figure}[!t]
		\centering
		\begin{tikzpicture}[line width=2, scale=2.5]
		\draw [fill=myblue] (0,0) -- (2,0) -- (2,1) -- (0,1) -- cycle;
		\draw [myorange, fill=myorange] plot [smooth cycle] coordinates {(1,0.8) (1.5,0.45) (1.2, 0.2) (0.8,0.2) (0.6, 0.3) (0.7,0.6) (0.5, 0.8)};
		
		\node at (0.95,0.475) {\LARGE $\Omega$};
		\node at (0.2, 0.8) {\textcolor{white}{\LARGE $\holdall$}};
		\node at (0.35, 0.25) {\textcolor{white}{\Large $\Omega\subout$}};
		\node at (1.5, 0.7) {\textcolor{white}{\LARGE $\Gamma$}};
		
		\draw[-latex, white, line width=1] (1.4,0.3) -- (1.55,0.15);
		\node at (1.65, 0.25) {\textcolor{white}{\LARGE $\normal$}};
		\end{tikzpicture}
		\caption{Sketch of the problem setup for $d=2$.}
		\label{fig:model_setup}
	\end{figure}
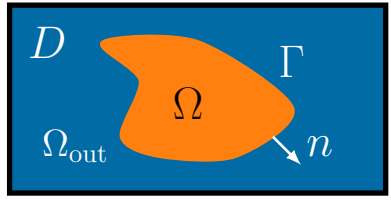
	
	Our elliptic shape interface problem is given by
	\begin{equation}
		\label{eq:interface_strong}
		\begin{aligned}
			&\min_{\Omega \in \admissiblegeom}\smallspace \costfunctional(\Omega, u) = \integral{D} \left(u - u\subdes\right)^2 \dmeas{x} \\
			&\text{s.t.} \quad \begin{alignedat}[t]{2}
				- \nabla \cdot (\interfacecoefficient\smallspace \grad u) &= f \quad &&\text{ in } D, \\
				u &= 0 \quad &&\text{ on } \partial D, \\
				\jump{u}{\Gamma} &= 0 \quad &&\text{ on } \Gamma,\\
				\jump{\interfacecoefficient\smallspace \partial_\normal u}{\Gamma} &= 0 \quad &&\text{ on } \Gamma.
			\end{alignedat}
		\end{aligned}
	\end{equation}
	Here, the parameter $\interfacecoefficient$ is given by a positive constant in each subdomain, i.e.,
	\begin{equation*}
		\interfacecoefficient = \begin{cases}
			\interfacecoefficient\subin = \text{const.} \quad &\text{ in } \Omega, \\
			\interfacecoefficient\subout = \text{const.} \quad &\text{ in } \Omega\subout,
		\end{cases}
	\end{equation*}
	such that the strong form of the PDE constraint \cref{eq:interface_strong} has to be understood formally. Moreover, the jump operator $\jump{\cdot}{\Gamma}$ is given by $\Gamma$ as $\jump{v}{\Gamma} = \left. v\right\rvert_{\Omega\subout} - \left. v\right\rvert_{\Omega}$.
	Hence, the final equations of system \cref{eq:interface_strong} enforce the continuity of the state variable and the flux. We assume that $f \in H^1(\holdall)$ and $u\subdes \in H^1(\holdall)$, which we require later to derive the shape derivative. The variational formulation of the PDE constraint in \cref{eq:interface_strong} is given by
	\begin{equation}
		\label{eq:interface_weak}
		\text{Find } u \in H^1_0(\holdall) \text{ such that } \quad \integral{\holdall} \interfacecoefficient \grad u \cdot \grad v \dmeas{x} = \integral{\holdall} f v \dmeas{x} \quad \text{ for all } v\in H^1_0(\holdall).
	\end{equation}
	As $\interfacecoefficient\subin, \interfacecoefficient\subout > 0$, \cref{eq:interface_weak} has a unique solution due to the Lax-Milgram lemma. The goal of the optimization problem \cref{eq:interface_strong} is to find the optimal shape of the inclusion in the sense that the $L^2$-distance between $u(\Omega)$ and some desired state $u\subdes$ is minimized. Typically, the desired state is obtained by measurements, such that \cref{eq:interface_strong} can be interpreted as a shape identification problem for the measurement $u\subdes$.
\end{example}

Now, we can use techniques from shape calculus to calculate sensitivities of the reduced cost functional $\reducedcostfunctional$ w.r.t. infinitesimal changes in the geometry, which is the foundation of efficient gradient-based optimization algorithms for problem \cref{eq:reduced_optimization_problem}. For a detailed introduction to shape calculus we refer, e.g., to the textbooks \cite{Delfour2011Shapes,Sokolowski1992Introduction}. We consider the so-called speed method, which is used to transform a domain $\Omega \subset \holdall$, where $\holdall$ is some hold-all domain, by the flow of a vector field $\vectorfield$. Here, we consider vector fields $\vectorfield \in C^k_0(\holdall;\R^d)$ for $k\geq 1$, i.e., the space of all $k$-times continuously differentiable functions from $\holdall$ to $\R^d$ with compact support in $\holdall$. A point $x_0 \in \Omega$ is moved along the flow of $\vectorfield$, which is described by the following initial value problem
\begin{equation}
	\label{eq:ivp}
	\dot{x}(t) = \vectorfield(x(t)), \quad x(0) = x_0.
\end{equation}
It is well-known (cf.\ \cite{Delfour2011Shapes}) that \cref{eq:ivp} has a unique solution $x(t)$ for all $t\in [0,\tau]$ if $\tau > 0$ is sufficiently small due to the regularity of $\vectorfield$. Therefore, we define the flow of $\vectorfield$ as
\begin{equation*}
	\flow_t^\vectorfield \colon  \R^d \to \R^d; \quad x_0 \mapsto \flow_t^\vectorfield x_0 = x(t).
\end{equation*}
%In fact, $\flow_t^\vectorfield$ is a diffeomorphism (see, e.g., \cite{Delfour2011Shapes}), which we utilize in the following to define the shape derivative.

\begin{definition}
	Let $\mathcal{S} \subset \Set{\Omega | \Omega \subset D}$, $J\colon \mathcal{S} \to \R$, $\Omega \in \mathcal{S}$, and $\vectorfield \in C^k_0(D;\R^d)$ for some $k\geq 1$. We denote by $\flow_t^\vectorfield$ the flow of $\vectorfield$ and we assume that $\flow_t^\vectorfield (\Omega) \in \mathcal{S}$ for all $t \in [0,\tau]$ for a sufficiently small $\tau > 0$. The Eulerian semiderivative of the shape functional $\reducedcostfunctional$ at $\Omega$ in direction $\vectorfield$ is given by the limit
	\begin{equation*}
		d\reducedcostfunctional(\Omega)[\vectorfield] := \lim\limits_{t\searrow 0} \frac{\reducedcostfunctional(\flow_t^\vectorfield(\Omega)) - \reducedcostfunctional(\Omega)}{t} = \left. \frac{d}{dt} \reducedcostfunctional(\flow_t^\vectorfield(\Omega)) \right\rvert_{t=0^+}
	\end{equation*}
	if it exists. 
	Further, let $\Xi$ be a topological vector subspace of $C^\infty_0(D;\R^d)$. The functional $\reducedcostfunctional$ is called shape differentiable at $\Omega$ (w.r.t. $\Xi$) if it has an Eulerian semiderivative at $\Omega$ in all directions $\vectorfield \in \Xi$ and, additionally, the mapping $d\reducedcostfunctional(\Omega): \Xi \to \R;\ \vectorfield \mapsto d\reducedcostfunctional(\Omega)[\vectorfield]$ is linear and continuous. Then, we call $d\reducedcostfunctional(\Omega)[\vectorfield]$ the shape derivative of $\reducedcostfunctional$ at $\Omega$ (w.r.t. $\Xi$) in direction $\vectorfield$.
\end{definition}

\begin{remark}
	The subspace $\Xi$ in the previous definition is used to incorporate additional geometrical constraints of the optimization problem \cref{eq:reduced_optimization_problem} into the definition of the shape derivative, i.e., the optimization is considered over all $\Omega \in \admissiblegeom$, which can require that some boundaries have to be fixed during the optimization (cf., e.g., \cite{Blauth2020Shape}).
\end{remark}

For our model problem, we have the following shape derivative.
\begin{proposition}
	\label{prop:volume}
	The reduced cost functional $\reducedcostfunctional$ corresponding to problem \cref{eq:interface_strong} is shape differentiable with shape derivative
	\begin{equation}
		\label{eq:shape_derivative_interface}
		\begin{aligned}
			d\reducedcostfunctional(\Omega)[\mathcal{\vectorfield}] = &\integral{\holdall} \left( u - u\subdes \right)^2 \div{\vectorfield} - 2 \left( u - u\subdes \right) \left(\grad u\subdes \cdot \vectorfield \right) \dmeas{x} \\
			&+ \integral{\holdall} \interfacecoefficient \left( \left( \div{\vectorfield} I - (D\vectorfield + D\vectorfield\transposed) \right) \grad u  \right) \cdot \grad p - \div{f \vectorfield} p \dmeas{x},
		\end{aligned}
	\end{equation}
	where $u \in H^1_0(\holdall)$ is the solution of \cref{eq:interface_weak}, $D\mathcal{V}$ denotes the Jacobian of $\mathcal{V}$, and the adjoint state $p \in H^1_0(\holdall)$ solves
	\begin{equation}
		\label{eq:interface_adjoint}
		\text{Find } p \in H^1_0(\holdall) \text{ s.t. } \quad \integral{\holdall} \interfacecoefficient\smallspace \grad p \cdot \grad v \dmeas{x} = -2\integral{\holdall} \left( u - u\subdes \right) v \dmeas{x} \quad \text{ f.a. } v\in H^1_0(\holdall).
	\end{equation}
\end{proposition}

\Cref{eq:shape_derivative_interface} is also known as volume formulation of the shape derivative as it only involves integrals over the domain $\holdall$. Under certain regularity assumptions, there exists an equivalent boundary formulation of the shape derivative due to the structure theorem, which will be summarized in the following and whose proof can be found in \cite[Theorem~3.6 and Corollary~1]{Delfour2011Shapes}.
%Due to Hadamard's structure theorem, which we recall below, there exist two equivalent formulations of the shape derivative, namely the volume and boundary formulations. For a proof of the structure theorem, we refer the reader to \cite[Theorem~3.6 and Corollary~1]{Delfour2011Shapes}.
\begin{theorem}[structure theorem]
	\label{thm:structure}
	Let $\reducedcostfunctional$ be a shape functional which is shape differentiable at some $\Omega \subset \R^d$ and let $\Gamma = \partial \Omega$ be compact. Further, let $k \geq 0$ be an integer for which the mapping $d\reducedcostfunctional(\Omega)\colon C^\infty_0(\holdall;\R^d) \to \R; \quad \vectorfield \mapsto d\reducedcostfunctional(\Omega)[\vectorfield]$
	is continuous w.r.t.\ the $C^k_0(\holdall;\R^d)$ topology and assume that $\Gamma$ is of class $C^{k+1}$. Then, there exists a continuous functional $\shapedistro\colon C^k(\Gamma) \to \R$ such that
	\begin{equation*}
		d\reducedcostfunctional(\Omega)[\vectorfield] = \shapedistro[\vectorfield \cdot n],
	\end{equation*}
	where $n$ is the outer unit normal on $\Gamma$. In particular, if $\shapedistro \in L^1(\Gamma)$, we have 
	\begin{equation*}
		d\reducedcostfunctional(\Omega)[\vectorfield] = \integral{\Gamma} \shapedistro \vectorfield\cdot \normal \dmeas{s}.
	\end{equation*}
\end{theorem}

For our model problem, we obtain the following.
\begin{proposition}
	\label{prop:boundary}
	Let $\Omega\subset \holdall$ be a compactly contained subset of $\holdall$ and assume that $\Gamma = \partial \Omega$ is of class $C^2$. Then, the shape derivative \cref{eq:shape_derivative_interface} has the following equivalent boundary representation
	\begin{equation*}
	d\reducedcostfunctional(\Omega)[\vectorfield] = \integral{\Gamma} \left( \vectorfield\cdot \normal \right) \jump{\interfacecoefficient\smallspace \tgrad u \cdot \tgrad p - \interfacecoefficient\smallspace \partial_\normal u\smallspace \partial_\normal p}{\Gamma} \dmeas{s},
	\end{equation*}
	where $u$ solves \cref{eq:interface_weak}, $p$ solves \cref{eq:interface_adjoint}, and $\tgrad v = \grad v - (\partial_\normal v) \smallspace \normal$ is the tangential gradient of $v$.
\end{proposition}
\begin{remark}
	The proofs of \cref{prop:volume,prop:boundary} can be found in \cite[Section~4]{Sturm2015Minimax} for $\beta(\abs{\grad u}^2, x) = \alpha$ such that $\beta'(\abs{\grad u}^2, x) = 0$ in \cite{Sturm2015Minimax}.
\end{remark}

\subsection{Steklov-Poincar\'e-Type Metrics for Shape Optimization}
\label{ssec:stekpoinc}

Let us now briefly recall the Steklov-Poincar\'e-type metrics from \cite{Schulz2016Efficient}.

We start by considering connected and compact subsets $\Omega \subset \holdall \subset \R^2$ with $C^\infty$ boundary, where $\holdall$ is a bounded hold-all domain. We define the space of all smooth two-dimensional shapes as in \cite{Michor2006Riemannian} by
\begin{equation*}
	\shapemanifold(S^1;\R^2) := \text{Emb}(S^1;\R^2) / \text{Diff}(S^1),
\end{equation*}
which is the set of all equivalence classes of $C^\infty$ embeddings of the unit circle $S^1\subset \R^2$ into $R^2$, i.e., $\text{Emb}(S^1, \R^2)$, where the equivalence relation is defined via the set of all $C^\infty$ diffeomorphisms of $S^1$ into itself, i.e., $\text{Diff}(S^1)$. Due to the equivalence relation, reparametrizations of geometries are not considered since they do not change the underlying shape. It is shown in \cite{Kriegl1997convenient} that $\shapemanifold$ is a smooth manifold. An element of $\shapemanifold(S^1, \R^2)$ is represented by a smooth curve $\Gamma\colon S^1 \to \R^2; \theta \mapsto \Gamma(\theta)$. The tangent space of $\shapemanifold$ at $\Gamma \in \shapemanifold$ is isomorphic to the set of all $C^\infty$ normal vector fields along $\Gamma$, i.e.,
\begin{equation*}
	\tangentspace_\Gamma \shapemanifold \cong \Set{h | h=\alpha \normal, \alpha\in C^\infty(\Gamma;\R)} \cong \Set{\alpha | \alpha \in C^\infty(\Gamma;\R)},
\end{equation*}
where $\normal$ is the outer unit normal vector on $\Gamma$.

As in \cite{Schulz2016Efficient}, we consider the following Steklov-Poincar\'e-type metric $\scalarproduct_\Gamma$ at some $\Gamma \in \shapemanifold$:
\begin{equation}
	\label{eq:def_scalarproduct}
	\scalarproduct_\Gamma \colon H^{-\frac{1}{2}}(\Gamma) \times H^{\frac{1}{2}}(\Gamma) \to \R; \quad (\alpha, \beta) \mapsto \integral{\Gamma} \alpha \left(\left( \stekpoinc_\Gamma \right)^{-1} \beta \right) \dmeas{s}.
\end{equation}
Here, $\stekpoinc_\Gamma$ is a symmetric and coercive operator defined by
\begin{equation*}
	\stekpoinc_\Gamma \colon H^{-\frac{1}{2}}(\Gamma) \to H^{\frac{1}{2}}(\Gamma); \quad \alpha \mapsto \trace{U}\cdot \normal,
\end{equation*}
with the trace operator $\mathrm{tr}\colon H^1_0(D; \R^d) \to H^{\frac{1}{2}}(\Gamma)$ and where $U\in H^1_0(\holdall;\R^d)$ solves
\begin{equation}
	\label{eq:stekpoinc_subproblem}
	\text{Find } U \in H^1_0(\holdall;\R^d) \text{ s.t. } \quad a_\Omega(U, V) = \integral{\Gamma} \alpha \left( \trace{V}\cdot n \right) \dmeas{s} \quad \text{ for all } V\in H^1_0(\holdall;\R^d)
\end{equation}
for a symmetric, continuous, and coercive bilinear form $a_\Omega \colon H^1_0(\holdall;\R^d) \times H^1_0(\holdall;\R^d) \to \R$. As discussed in \cite{Schulz2016Efficient}, $\stekpoinc_\Gamma$ corresponds to a projected Steklov-Poincar\'e operator, and we refer the reader to \cite{Schulz2016Efficient} for a discussion regarding the invertibility of the operator $\stekpoinc_\Gamma$. Note that \cref{eq:stekpoinc_subproblem} is well-posed due to the Lax-Milgram lemma. The metric $\scalarproduct_\Gamma$ induces the following norm
\begin{equation*}
	\norm{\alpha}{\Gamma} = \sqrt{\scalarproduct_\Gamma(\alpha, \alpha)} \quad \text{ for } \alpha \in H^{\frac{1}{2}}(\Gamma).
\end{equation*}
The Steklov-Poincar\'e metric $\scalarproduct_\Gamma$ yields a Riemannian metric on $\shapemanifold$ if we restrict it to the tangent space $\tangentspace_\Gamma\shapemanifold$. As discussed in \cite{Blauth2021Nonlinear}, this setting can be readily adapted to the case where parts of $\Gamma$ are fixed during the optimization.

Let us now consider the relation between the metric $\scalarproduct$ and shape calculus. We assume that $\reducedcostfunctional$ has the shape derivative
\begin{equation*}
	d\reducedcostfunctional(\Omega)[\vectorfield] = \integral{\Gamma} \shapedistro\smallspace \trace{\vectorfield}\cdot \normal \dmeas{s},
\end{equation*}
where $\shapedistro\in L^2(\Gamma)$. By definition, the Riemannian shape gradient w.r.t.\ $\scalarproduct_\Gamma$ solves
\begin{equation*}
	\text{Find } \rieshapegrad \in \tangentspace_\Gamma\shapemanifold \text{ such that } \quad \scalarproduct_\Gamma\left( \rieshapegrad, \phi \right) = \integral{\Gamma} \shapedistro \phi \dmeas{s} \quad \text{ for all } \phi \in \tangentspace_\Gamma\shapemanifold.
\end{equation*}
Due to the definition of $\scalarproduct_\Gamma$, it is easy to see that we have $\rieshapegrad = \stekpoinc_\Gamma \shapedistro$ and, therefore, it holds that $\rieshapegrad = \trace{\gradientdefo} \cdot \normal$, where $\gradientdefo$ solves the problem
\begin{equation}
	\label{eq:def_gradientdefo}
	\text{Find } \gradientdefo \in H^1_0(\holdall;\R^d) \text{ such that } \quad a_\Omega(\gradientdefo, \vectorfield) = d\reducedcostfunctional(\Omega)[\vectorfield] \quad \text{ for all } \vectorfield \in H^1_0(\holdall;\R^d).
\end{equation}
We call the unique solution $\gradientdefo$ of \cref{eq:def_gradientdefo} the gradient deformation of $\reducedcostfunctional$ at $\Omega$. It can be interpreted as extension of the Riemannian shape gradient $\rieshapegrad$ to the entire domain $\Omega$. Note that the gradient deformation is often used for gradient-based shape optimization algorithms (cf. \cite{Schulz2016Efficient,Blauth2021Nonlinear,Gangl2015Shape,Hohmann2019Shape,Blauth2021Model}).

%Note that due to the coercivity of $a_\Omega$, there exists a constant $C>0$ such that
%\begin{equation*}
%	d\reducedcostfunctional(\Omega)[-\gradientdefo] = a_\Omega(-\gradientdefo, \gradientdefo) \leq - C \norm{\gradientdefo}{H^1(\Omega;\R^d)} \leq 0.
%\end{equation*}
%Hence, an infinitesimal transformation of $\Omega$ along the flow associated to the negative gradient deformation yields a descent in the shape functional $\reducedcostfunctional$. This is often used in the sense of a steepest descent method in the context of numerical shape optimization (see, e.g., \cite{Gangl2015Shape,Blauth2021Model,Hohmann2019Shape}). Additionally, the gradient deformation is an important part of other gradient-based shape optimization algorithms, such as the quasi-Newton methods introduced in \cite{Schulz2016Efficient} and the nonlinear conjugate gradient methods proposed in \cite{Blauth2021Nonlinear}.

\begin{remark}
	As discussed in \cite{Schulz2016Efficient}, the Riemannian shape gradient $\rieshapegrad = \trace{\gradientdefo} \cdot \normal$ may not be a part of the tangent space $\tangentspace_\Gamma \shapemanifold$ as $\gradientdefo$ might not have a $C^\infty$ regularity. However, if the bilinear form $a_\Omega$ is associated to a second order elliptic differential operator with sufficiently smooth coefficients and if $\shapedistro$ and $\Omega$ are sufficiently smooth, the gradient deformation $\gradientdefo$ is an element of $C^\infty(\closure{\Omega};\R^d)$ by the theorem of infinite differentiability (cf.~\cite[Theorem~6, Section~6.3]{Evans2010Partial}). In this case, we have $\rieshapegrad \in \tangentspace_\Gamma\shapemanifold$.
\end{remark}

\begin{remark}
	For our numerical examples presented in \cref{ssec:uniform_flow_distribution,ssec:fluent}, the requirements of the structure theorem are not satisfied, as the considered geometries only have a Lipschitz boundary. However, we can still formally apply the Riemannian framework presented in this section and obtain very convincing numerical results.
\end{remark}

%\begin{Remark}
%	The Steklov-Poincar\'e metric $\scalarproduct_\Gamma$ is an appropriate inner product for shape optimization. It enables us to define scalar product between deformation vector fields, which is required to derive limited-memory quasi-Newton methods and, therefore, allows us to derive the space mapping methods for shape optimization in \cref{sec:space_mapping}. 
%\end{Remark}

\subsection{Basics of Riemannian Geometry}
\label{ssec:riemannian_geometry}

Finally, we briefly recall some basic ideas and notations from Riemannian geometry, following \cite{Absil2008Optimization}. For some $\Gamma \in \shapemanifold$, a retraction $\retraction_\Gamma$ is a smooth mapping $\retraction_\Gamma \colon \tangentspace_\Gamma \shapemanifold \to \shapemanifold$ which satisfies the relations $\retraction_\Gamma(0_\Gamma) = \Gamma$ and $D\retraction_\Gamma(0_\Gamma) = \text{id}_{\tangentspace_\Gamma \shapemanifold}$. Here, $D\retraction_\Gamma$ is the derivative of $\retraction_\Gamma$ and $0_\Gamma$ and $\text{id}_{\tangentspace_\Gamma \shapemanifold}$ are the zero element and identity mapping in $\tangentspace_\Gamma \shapemanifold$, respectively. A special retraction is given by the exponential map $\exp_\Gamma \colon \tangentspace_\Gamma \shapemanifold \to \shapemanifold; \alpha \mapsto \exp_\Gamma \alpha = \geodesic_\Gamma^\alpha(1)$, where $\geodesic_\Gamma^\alpha \colon [0,1] \to \shapemanifold$ is the unique geodesic starting at $\geodesic_\Gamma^\alpha(0) = \Gamma$ with $\dot{\geodesic}_\Gamma^\alpha(0) = \alpha$. Moreover, we denote by $\vectortransport \colon \tangentspace \shapemanifold \oplus \tangentspace \shapemanifold \to \tangentspace \shapemanifold; (\alpha, \beta) \mapsto \vectortransport_{\alpha} (\beta)$ a vector transport which satisfies the following properties. For $\alpha, \beta \in \tangentspace_\Gamma \shapemanifold$ it holds that $\vectortransport_\alpha (\beta)$ is an element of $\tangentspace_{\retraction_\Gamma (\alpha)} \shapemanifold$. Further, we have the relations $\vectortransport_{0_\Gamma} (\alpha) = \alpha$ and $\vectortransport_{\gamma} (a \alpha + b \beta) = a \vectortransport_\gamma (\alpha) + b \vectortransport_{\gamma} (\beta)$ for $a, b \in \R$. Note that the parallel transport, as defined in, e.g., \cite{Absil2008Optimization, Ring2012Optimization}, is a vector transport associated to the exponential map.

\section{Space Mapping Methods for Shape Optimization}
\label{sec:space_mapping}

In this section, we propose novel space mapping methods for PDE constrained shape optimization based on shape calculus. We first introduce the methods in a Riemannian setting and then investigate them in a volume-based setting suitable for numerical discretization.

\subsection{A Riemannian Framework for Space Mapping Shape Optimization}
\label{ssec:general_sm}

We consider the optimization of a system where the system's shape is the optimization variable. As in \cref{ssec:space_mapping}, we assume that we have a fine model, which is very detailed but expensive to evaluate, and a coarse model, which is a rough approximation but cheap to evaluate, for the system. Associated to both models we have the fine and coarse model response functions $\fineresponse \colon \shapemanifold \to Y$ and $\coarseresponse \colon \shapemanifold \to Y$. Here, the control spaces for the models are given by the shape space $\shapemanifold$, i.e., the models depend on the geometry, as usual in shape optimization.

Like in \cref{ssec:space_mapping}, we consider the fine model optimization problem
\begin{equation}
	\label{eq:fine_opt}
	\min_{\Gamma\in \shapemanifold}\smallspace \reducedcostfunctional(\fineresponse(\Gamma)),
\end{equation}
as well as the corresponding coarse model optimization problem
\begin{equation}
	\label{eq:coarse_opt}
	\min_{\Gamma \in \shapemanifold}\smallspace \reducedcostfunctional(\coarseresponse(\Gamma)),
\end{equation}
with a cost functional $\reducedcostfunctional\colon Y \to \R$. To obtain well-posed problems, we assume that \cref{eq:fine_opt} and \cref{eq:coarse_opt} each have a unique minimizer, which we denote by
\begin{equation*}
	\Gamma\subfine^* := \argmin_{\Gamma\in \shapemanifold}\smallspace \reducedcostfunctional(\fineresponse(\Gamma)) \qquad \text{ and } \qquad \Gamma\subcoarse^* := \argmin_{\Gamma\in \shapemanifold}\smallspace \reducedcostfunctional(\coarseresponse(\Gamma)).
\end{equation*}
To be compatible with the Riemannian framework for shape optimization from \cref{ssec:stekpoinc}, the geometry is denoted by its boundary $\Gamma$ instead of the domain $\Omega$, in contrast to our notation in \cref{ssec:shape_calculus}.

\begin{remark}
	\label{rem:reduced_shape}
	For PDE constrained shape optimization, a typical setting looks as follows: The fine model shape optimization problem is given by
	\begin{equation*}
		\min_{\Omega \in \admissiblegeom}\smallspace \costfunctional(\Omega, u) \quad \text{ subject to } \quad e\subfine(\Omega, u) = 0.
	\end{equation*}
	Correspondingly, the coarse model optimization problem is given by
	\begin{equation*}
	\min_{\Omega\in \admissiblegeom}\smallspace \costfunctional(\Omega, u) \quad \text{ subject to } \quad e\subcoarse(\Omega, u) = 0.
	\end{equation*}
	Here, the difference of the models comes from the different PDE constraints $e\subfine$ and $e\subcoarse$. The fine model may involve complex nonlinearities or large systems of equations, whereas the coarse model may consist of a less complex and smaller system. 
%	We assume that the PDEs admit a unique solution for all $\Omega\in\admissiblegeom$, which we denote by $u\subfine(\Omega)$ and $u\subcoarse(\Omega)$, respectively. Finally, we arrive at the reduced optimization problems
%	\begin{equation*}
%		\min_{\Omega \in \admissiblegeom} \reducedcostfunctional\subfine(\Omega) = \costfunctional(\Omega, u\subfine(\Omega)) \qquad \text{ and } \qquad \min_{\Omega\in \admissiblegeom} \reducedcostfunctional\subcoarse(\Omega) = \costfunctional(\Omega, u\subcoarse(\Omega)),
%	\end{equation*}
%	which play the role of the fine and coarse model optimization problems \cref{eq:fine_opt} and \cref{eq:coarse_opt}, respectively. 
	However, to compare the fine and coarse model, their respective system responses have to be restricted to a common subdomain, which is invariant during the optimization, such as a fixed boundary. We encounter prototypical examples for such settings in \cref{sec:numerics}.
\end{remark}

Now, we introduce a space mapping function $\spacemappingfunction$ which is defined as the vector field
\begin{equation*}
	\spacemappingfunction\colon \shapemanifold \to \tangentspace\shapemanifold; \quad \shapemanifold \ni \Gamma\subfine \mapsto \spacemappingfunction(\Gamma\subfine) \in \tangentspace_{\Gamma\subfine} \shapemanifold,
\end{equation*}
where
\begin{equation*}
	\spacemappingfunction(\Gamma\subfine) := \vectortransport_\eta \left( \argmin_{\alpha \in \tangentspace_{\hat{\Gamma}}\shapemanifold}\smallspace \misalignmentfunction\left( \coarseresponse(\retraction_{\hat{\Gamma}}(\alpha)), \fineresponse(\Gamma\subfine) \right) \right)
\end{equation*}
with $\eta = \retraction_{\hat{\Gamma}}^{-1} (\Gamma\subfine)$. Here, $\hat{\Gamma} \in \shapemanifold$ is a fixed reference domain and we have a misalignment function $\misalignmentfunction\colon Y\times Y \to \R$, which measures the discrepancy between its arguments. Typically, the misalignment function $\misalignmentfunction$ has the form $\misalignmentfunction\left( \coarseresponse(\Gamma\subcoarse), \fineresponse(\Gamma\subfine) \right) = \norm{\coarseresponse(\Gamma\subcoarse) - \fineresponse(\Gamma\subfine)}{}$ for some suitable norm $\norm{\cdot}{}$. We assume that the misalignment function is exact such that
\begin{equation*}
	\argmin_{\alpha \in \tangentspace_{\hat{\Gamma}}\shapemanifold}\smallspace \misalignmentfunction\left( \coarseresponse(\retraction_{\hat{\Gamma}}(\alpha)), \coarseresponse(\retraction_{\hat{\Gamma}} (\beta)) \right) = \beta \quad \text{ for all } \beta\in \tangentspace_{\hat{\Gamma}} \shapemanifold.
\end{equation*}
Hence, the space mapping function computes the geometry update $\alpha$, which minimizes the discrepancy between $\coarseresponse(\retraction_{\hat{\Gamma}}(\alpha))$ and $\fineresponse(\Gamma\subfine)$.

We assume that any two shapes $\Gamma, \Gamma^* \in \shapemanifold$ can be connected via retractions, i.e., there exists an element $\alpha \in \tangentspace_\Gamma\shapemanifold$ such that $\Gamma^* = \retraction_\Gamma (\alpha)$ (cf.~\cite{Ring2012Optimization}). Therefore, we can rewrite the coarse model optimum in terms of the fixed reference domain $\hat{\Gamma}$ as $\Gamma\subcoarse^* = \retraction_{\hat{\Gamma}}(\alpha^*)$ for some $\alpha^* \in \tangentspace_{\hat{\Gamma}} \shapemanifold$. Similarly to \cref{ssec:space_mapping}, we have the following relation: Under the assumption that the fine and coarse model responses are similar near their respective optima, i.e., $\coarseresponse(\Gamma\subcoarse^*) \approx \fineresponse(\Gamma\subfine^*)$, we can expect that
\begin{equation*}
	\begin{aligned}
		\spacemappingfunction(\Gamma\subfine^*) &= \vectortransport_{\eta^*} \left( \argmin_{\alpha \in \tangentspace_{\hat{\Gamma}}\shapemanifold}\smallspace \misalignmentfunction\left( \coarseresponse(\retraction_{\hat{\Gamma}}(\alpha)), \fineresponse(\Gamma\subfine^*) \right) \right) \\
		&\approx \vectortransport_{\eta^*} \left( \argmin_{\alpha \in \tangentspace_{\hat{\Gamma}}\shapemanifold}\smallspace \misalignmentfunction\left( \coarseresponse(\retraction_{\hat{\Gamma}}(\alpha)), \coarseresponse(\Gamma\subcoarse^*) \right) \right) \\
		&= \vectortransport_{\eta^*} \left( \argmin_{\alpha \in \tangentspace_{\hat{\Gamma}}\shapemanifold}\smallspace \misalignmentfunction\left( \coarseresponse(\retraction_{\hat{\Gamma}}(\alpha)), \coarseresponse(\retraction_{\hat{\Gamma}} (\alpha^*)) \right) \right) = \vectortransport_{\eta^*} (\alpha^*),
	\end{aligned}
\end{equation*}
where $\eta^* = \retraction_{\hat{\Gamma}}^{-1} (\Gamma\subfine^*)$, due to the fact that $\Gamma\subcoarse^* = \retraction_{\hat{\Gamma}}(\alpha^*)$ and the exactness of the misalignment function. Therefore, our space mapping method first computes the optimal coarse model geometry $\Gamma\subcoarse^* = \retraction_{\hat{\Gamma}}(\alpha^*)$, and then approximately solves 
\begin{equation}
	\label{eq:sm_equation}
	\spacemappingfunction(\Gamma\subfine^*) = \vectortransport_{\eta^*} (\alpha^*)
\end{equation}
to compute an approximation of the optimal fine model geometry $\Gamma\subfine^*$.

\begin{remark}
	The assumption that any two shapes $\Gamma, \Gamma^*$ can be connected via retractions is common when considering optimization on Riemannian manifolds (see, e.g., \cite{Ring2012Optimization}). Additionally, for the implementation of the space mapping methods, we only require that there exists a retraction connecting the reference shape with the shapes obtained from optimizing the coarse model. 
\end{remark}

\begin{algorithm2e}[!b]
	\KwIn{Reference geometry $\hat{\Gamma}$, tolerance $\tau \in (0,1)$, maximum number of iterations $k_{\text{max}} \in \mathbb{N}$ \label{line:input}, memory size $m \geq 0$}
	Compute $\Gamma\subcoarse^* = \argmin_{\Gamma\in \shapemanifold}\smallspace \reducedcostfunctional(\coarseresponse(\Gamma))$ and define $\alpha^* = \retraction_{\hat{\Gamma}}^{-1} (\Gamma\subcoarse^*) \in \tangentspace_{\hat{\Gamma}} \shapemanifold$ \\
	Set $\Gamma^0 = \Gamma\subcoarse^*$\\
	\For{$k=0,1,2,\dots, k_\textrm{max}$}{
		Compute $\spacemappingfunction(\Gamma^k) \in \tangentspace_{\Gamma^k} \shapemanifold$ \\
		Set $\eta^k = \retraction_{\hat{\Gamma}}^{-1} (\Gamma^k)$ \\
		\If{$\norm{\spacemappingfunction(\Gamma^k) - \vectortransport_{\eta^k} (\alpha^*)}{\Gamma^k} \leq \tau \norm{\vectortransport_{\eta^k} (\alpha^*)}{\Gamma^k}$}{
			\KwRet{\upshape $\Gamma\subfine^* = \Gamma^k$ as (approximate) solution}}
		Compute $\delta^k = -(B^k)^{-1} \left( \spacemappingfunction(\Gamma^k) - \vectortransport_{\eta^k} (\alpha^*) \right) \in \tangentspace_{\Gamma^k} \shapemanifold$ with \cref{algo:riemannian_loop} using the stored $s^i$ and $y^i$ for $i=k-m, \dots, k-1$ \\
		Update the geometry by $\Gamma^{k+1} = \retraction_{\Gamma^k} (\delta^k)$ \\
		%
		%		\Comment{Update the approximation of $(B^k)^{-1}$}
		%
		Compute $\nu = (B^k)^{-1}\left( \vectortransport_{\retraction_{\Gamma^{k+1}}^{-1} (\Gamma^{k})} (\spacemappingfunction(\Gamma^{k+1})) - \spacemappingfunction(\Gamma^{k}) \right)$ with \cref{algo:riemannian_loop} using the stored $s^i$ and $y^i$ for $i=k-m, \dots, k-1$ \\
		Store $s^k = \frac{1}{\scalarproduct_{\Gamma^k}\left( \delta^k, \nu \right)}(\delta^k - \nu) \in \tangentspace_{\Gamma^k}\shapemanifold$ and $y^k = \delta^k \in \tangentspace_{\Gamma^k} \shapemanifold$\\
		\If{$k > m$}{Discard $s^{k-m}$ and $y^{k-m}$ from storage}
	}
	\caption{Aggressive Space Mapping Method for Shape Optimization.}
	\label{algo:riemannian_asm}
\end{algorithm2e}

\begin{algorithm2e}[!b]
	\KwIn{Right-hand side $\beta \in \tangentspace_{\Gamma^k} \shapemanifold$ and stored elements $s^i, y^i \in \tangentspace_{\Gamma^i} \shapemanifold$ for $i=k-m, \dots, k-1$}
	Set $\gamma = \beta$ \\
	\For{$i=k-m,\dots,k-1$}{
		Set $\mu^i = \retraction_{\Gamma^i}^{-1} (\Gamma^k)$ \\
		Compute $\theta = \scalarproduct_{\Gamma^{k}}\left( \gamma, \vectortransport_{\mu^i} (y^i) \right)$ \\
		$\gamma = \gamma + \theta\smallspace \vectortransport_{\mu^i} (s^i)$ \\
	}
	\KwRet{$\gamma = (B^k)^{-1} \beta$}
	\caption{Limited memory loop to compute $\gamma = (B^k)^{-1} \beta$}
	\label{algo:riemannian_loop}
\end{algorithm2e}

We now propose an aggressive space mapping (ASM) method in the context of shape optimization, which uses a quasi-Newton iteration to solve equation \cref{eq:sm_equation}. This amounts to the iteration (cf.~\cite{Absil2008Optimization,Ring2012Optimization})
\begin{equation*}
	B^k \delta^k = -\left(\spacemappingfunction(\Gamma^k) - \vectortransport_{\eta^k} (\alpha^*) \right) , \qquad \Gamma^{k+1} = \retraction_{\Gamma^k} (\delta^k),
\end{equation*}
where $\eta^k = \retraction_{\hat{\Gamma}}^{-1} (\Gamma^k)$ and $B^k\colon \tangentspace_{\Gamma^k} \shapemanifold \to \tangentspace_{\Gamma^k} \shapemanifold$ is a linear operator. For the iteration, we typically use the initial guesses $\Gamma^0 = \Gamma\subcoarse^*$ and $B^0 = \text{id}_{\tangentspace_{\Gamma^0} \shapemanifold}$. As update formula for the approximate Jacobian we have
\begin{equation*}
	B^{k+1} \gamma = \tilde{B}^{k} \gamma + \frac{\spacemappingfunction(\Gamma^{k+1}) - \vectortransport_{\delta^k} \left(\spacemappingfunction(\Gamma^{k}) + B^k \delta^k\right) }{\norm{\vectortransport_{\delta^k} (\delta^k)}{\Gamma^{k+1}}^2} \scalarproduct_{\Gamma^{k+1}} \left( \vectortransport_{\delta^k} (\delta^k), \gamma \right) \quad \text{ f.a. } \gamma \in \tangentspace_{\Gamma^{k+1}} \shapemanifold,
%	B^{k+1} \gamma = \tilde{B}^{k} \gamma + \frac{\spacemappingfunction(\Gamma^{k+1}) - \vectortransport_{\eta^{k+1}} (\alpha^*)}{\norm{\vectortransport_{\delta^k} (\delta^k)}{\Gamma^{k+1}}^2} \scalarproduct_{\Gamma^{k+1}} \left( \vectortransport_{\delta^k} (\delta^k), \gamma \right) \quad \text{ for all } \gamma \in \tangentspace_{\Gamma^{k+1}} \shapemanifold,
\end{equation*}
where $\tilde{B}^k = \vectortransport_{\delta^k} \circ B^{k} \circ (\vectortransport_{\delta^k})^{-1}$.
% and $\eta^{k+1} = \retraction_{\hat{\Gamma}}^{-1} (\Gamma^{k+1})$.

To summarize the above discussion, our resulting ASM method for shape optimization is detailed in \cref{algo:riemannian_asm}. For our setting in the context of PDE constrained shape optimization, we employ a limited-memory variant of Broyden's method for the space mapping iteration, which is based on the inverse of the Jacobian approximation. This formulation of Broyden's method makes use of the fact that the application of $(B^k)^{-1}$ to some $\beta \in \tangentspace_{\Gamma^k}\shapemanifold$ can be computed efficiently using a recursive algorithm based on the previous iterates $s^i$ and $y^i$, which are defined in \cref{algo:riemannian_asm}. The corresponding procedure for computing $(B^k)^{-1} \beta$ can be found in \cref{algo:riemannian_loop} and we refer the reader, e.g., to \cite[Chapter~2.1]{Deuflhard2011Newton} for a description of this recursive implementation in a classical setting.

%\begin{Remark}
%	Note that the definition of the space mapping function as vector field over $\shapemanifold$ is essential for deriving the space mapping method in the context of the Steklov-Poincar\'e metric $\scalarproduct$ as it shifts the view from the elements of the shape space $\shapemanifold$ (geometries) to elements of the tangent space $\tangentspace\shapemanifold$ (deformations). The Steklov-Poincar\'e metric then induces the tanget space $\tangentspace\shapemanifold$ with a local Hilbert space topology, which enables us to extend the ASM method to the field of shape optimization.
%\end{Remark}

\begin{remark}
	\Cref{algo:riemannian_asm} has low memory requirements as only the $2m$ elements $s^{i}, y^{i}$ for $i=k-m,\dots,k-1$ have to be stored. Hence, the algorithm is also feasible for large-scale problems as it avoids storing large dense matrices as discretizations of the operators $B^k$ or $(B^k)^{-1}$.
\end{remark}

\subsection{Volume-Based Formulation of the Space Mapping Methods for Shape Optimization}
\label{ssec:volume_sm}

Even though the space mapping framework presented in the previous section is complete from a theoretical perspective, it is not yet suited for numerical discretization (cf.~the discussion in \cite{Blauth2021Nonlinear}). For this reason, we now present a volume-based view of the method which can be discretized straightforwardly. 

For our volume-based realization of the space mapping technique, we proceed analogously to \cite[Section~3.3]{Blauth2021Nonlinear} and extend the retraction and vector transport such that they act on vector fields instead of elements of the tangent space $\tangentspace_\Gamma \shapemanifold$. We denote these extensions by
\begin{equation*}
	\numretraction_\Omega \colon H^1(\Omega;\R^d) \to \admissiblegeom \qquad \text{ and } \qquad \numtransport \colon H^1(\Omega;\R^d) \times H^1(\Omega;\R^d) \to H^1(\numretraction_\Omega(\vectorfield); \R^d),
\end{equation*}
where the retraction $\numretraction_\Omega$ maps some vector field $\vectorfield \in H^1(\Omega;\R^d)$ to a domain $\numretraction_\Omega(\vectorfield) \subset \admissiblegeom$ deformed by it, and $\numtransport$ maps two vector fields $\mathcal{V},\mathcal{W} \in H^1(\Omega,\R^d)$ to a vector field $\numtransport_{\vectorfield} \mathcal{W} \in H^1(\numretraction_\Omega(\vectorfield))$. We assume that for $\alpha = \trace{\vectorfield} \cdot \normal \in \tangentspace_\Gamma \shapemanifold$ and $\beta = \trace{\mathcal{W}} \cdot \normal \in \tangentspace_\Gamma \shapemanifold$ with $\vectorfield, \mathcal{W} \in H^1(\Omega;\R^d)$ it holds that $\retraction_\Gamma (\alpha)$ is the boundary of the domain $\numretraction_\Omega (\vectorfield)$ and that 
\begin{equation*}
	\vectortransport_\alpha \beta = \trace{\numtransport_\vectorfield \mathcal{W}} \cdot \normal \in \tangentspace_{\retraction_\Gamma (\alpha)} \shapemanifold,
\end{equation*}
to be compatible with the notations from \cref{ssec:riemannian_geometry}. By slight misuse of notation and for simplicity of presentation, we denote the inverse of $\numretraction$ as $\vectorfield = \numretraction_{\Omega}^{-1} (\tilde{\Omega}) =  \tilde{\Omega} - \Omega$, which is motivated by our numerical discretization as explained in \cref{ssec:implementation_details}. Here, the difference $\tilde{\Omega} - \Omega$ is only to be understood formally and is defined as $\tilde{\Omega} - \Omega := \numretraction_{\Omega}^{-1}(\tilde{\Omega})$. This extends the previous notions of retraction and vector transport to our volume-based setting.

Note that the Steklov-Poincar\'e-type metric from \cref{ssec:stekpoinc} is constructed such that the following holds true. Let the Riemannian shape gradient be given by $\rieshapegrad \in \tangentspace_\Gamma \shapemanifold$, where $\rieshapegrad = \trace{\gradientdefo} \cdot \normal$ as in \cref{ssec:stekpoinc}, and consider an element $\beta \in \tangentspace_\Gamma \shapemanifold$ given by $\beta = \trace{\vectorfield}\cdot \normal$ for some vector field $\vectorfield$. Then, we have that
\begin{equation}
\label{eq:relation_scalar_products}
\scalarproduct_\Gamma\left( \rieshapegrad, \beta \right) = \integral{\Gamma} \shapedistro \beta \dmeas{s} = \integral{\Gamma} \shapedistro\ \trace{\vectorfield} \cdot \normal \dmeas{s} = d\reducedcostfunctional(\Omega)[\vectorfield] = a_\Omega(\gradientdefo, \vectorfield).
\end{equation}
As discussed in \cite{Schulz2016Efficient,Blauth2021Nonlinear}, this relation also holds true in case we consider a linear combination of (transported) gradient-type vectors. Particularly, for (a linear combination of) gradient-type vectors $\rieshapegrad = \trace{\gradientdefo} \cdot \normal$, the Steklov-Poincar\'e metric $\scalarproduct$ induces the following norm
\begin{equation*}
	\norm{\rieshapegrad}{\Gamma} = \sqrt{\scalarproduct_\Gamma\left( \rieshapegrad, \rieshapegrad \right)} = \sqrt{a_\Omega\left( \gradientdefo, \gradientdefo \right)} = \norm{\gradientdefo}{a_\Omega}.
\end{equation*}
For the numerical solution of shape optimization problems, typically gradient-based algorithms are used. Such algorithms have the form as in \cite[Algorithms~1 and~2]{Blauth2021Nonlinear}. If we employ such gradient-based algorithms to solve the coarse model optimization problems (which we do in the numerical experiments in \cref{sec:numerics}), the search directions in each inner iteration are given as linear combinations of gradient-type vectors. Moreover, the complete deformation of such a gradient-based algorithm is just the linear sum of all iterations. This implies that all elements $\spacemappingfunction(\Gamma\subfine) \in \tangentspace_{\Gamma\subfine}\shapemanifold$ as well as $\alpha^* \in \tangentspace_{\hat{\Gamma}} \shapemanifold$ can be interpreted as linear combinations of gradient-type vectors. This allows us to compute the required scalar products for the ASM method in \cref{algo:riemannian_asm,algo:riemannian_loop} via scalar products between domain deformations with the help of the bilinear form $a$ and \cref{eq:relation_scalar_products}. Hence, we can recast the space mapping framework in a volume-based way.

Now, we formulate our volume-based variant of the ASM method, where we replace the elements of the tangent space by vector fields and the scalar product $\scalarproduct$ with the one induced by the bilinear form $a$. Analogously to before, assume that we have a fine and a coarse model for a system we want to optimize. We have response functions $\fineresponse\colon \admissiblegeom \to Y$ and $\coarseresponse\colon \admissiblegeom \to Y$ for the fine and coarse model, respectively, which map from some set of admissible geometries to a Banach space $Y$. We define the fine and coarse model shape optimization problems as 
\begin{equation*}
	\min_{\Omega \in \admissiblegeom} \reducedcostfunctional(\fineresponse(\Omega)) \qquad \text{ and } \qquad \min_{\Omega \in \admissiblegeom} \reducedcostfunctional(\coarseresponse(\Omega)),
\end{equation*}
respectively, and we assume that both problems have unique minimizers $\Omega\subfine^*$ and $\Omega\subcoarse^*$. Typically, these problems correspond to PDE constrained shape optimization problems (cf.~\cref{rem:reduced_shape}). We define a space mapping function $\spacemappingfunction$ by
\begin{equation*}
	\spacemappingfunction(\Omega\subfine) = \numtransport_{\Omega\subfine - \hat{\Omega}}\left( \argmin_{\vectorfield\in H^1(\hat{\Omega}; \R^d)} \misalignmentfunction\left( \coarseresponse\left(\numretraction_{\hat{\Omega}}(\vectorfield)\right), \fineresponse(\Omega\subfine) \right) \right),
\end{equation*}
where $\hat{\Omega}$ is a fixed reference domain and the misalignment function $\misalignmentfunction \colon Y\times Y \to \R$ measures the discrepancy between its arguments. As before, we assume that the misalignment function is exact in the sense that 
\begin{equation*}
	\argmin_{\vectorfield\in H^1(\hat{\Omega}; \R^d)} \misalignmentfunction \left( \coarseresponse\left( \numretraction_{\hat{\Omega}}(\vectorfield) \right), \coarseresponse\left( \numretraction_{\hat{\Omega}}(\mathcal{W}) \right) \right) = \mathcal{W} \quad  \text{ for all } \mathcal{W}.
\end{equation*}

To state the space mapping method, we assume that any two geometries $\Omega, \Omega^* \in \admissiblegeom$ can be connected via the retraction $\numretraction$, i.e., there exists some vector field $\vectorfield^*$ such that $\Omega^* = \numretraction_\Omega(\vectorfield^*)$. This allows us to rewrite the coarse model optimum geometry as $\Omega\subcoarse^* = \numretraction_{\hat{\Omega}}(\vectorfield^*)$. If we assume that the fine and coarse model responses are similar near their respective optima, i.e., if $\coarseresponse(\Omega\subcoarse^*) \approx \fineresponse(\Omega\subfine^*)$, then we can expect that
\begin{equation*}
	\begin{aligned}
		\spacemappingfunction(\Omega\subfine^*) &= \numtransport_{\Omega\subfine^* - \hat{\Omega}} \left( \argmin_{\vectorfield\in H^1(\hat{\Omega}; \R^d)} \misalignmentfunction\left( \coarseresponse\left( \numretraction_{\hat{\Omega}}(\vectorfield) \right), \fineresponse(\Omega\subfine^*) \right) \right) \\
		&\approx \numtransport_{\Omega\subfine^* - \hat{\Omega}} \left( \argmin_{\vectorfield\in H^1(\hat{\Omega}; \R^d)} \misalignmentfunction\left( \coarseresponse\left( \numretraction_{\hat{\Omega}}(\vectorfield) \right), \coarseresponse(\Omega\subcoarse^*) \right) \right) \\
		&= \numtransport_{\Omega\subfine^* - \hat{\Omega}} \left( \argmin_{\vectorfield\in H^1(\hat{\Omega}; \R^d)} \misalignmentfunction\left( \coarseresponse\left( \numretraction_{\hat{\Omega}}(\vectorfield) \right), \coarseresponse\left( \numretraction_{\hat{\Omega}}(\vectorfield^*)  \right) \right) \right) = \numtransport_{\Omega\subfine^* - \hat{\Omega}}\left( \vectorfield^* \right),
	\end{aligned}
\end{equation*}
analogously to our previous discussion. Hence, our space mapping method first computes $\Omega\subcoarse^*$, and, thus, $\vectorfield^* = \Omega\subcoarse^* - \hat{\Omega}$, and then solves the equation
\begin{equation}
	\label{eq:quasi_newton}
	\spacemappingfunction(\Omega\subfine^*) = \numtransport_{\Omega\subfine^* - \hat{\Omega}} (\vectorfield^*)
\end{equation}
to compute an approximation of $\Omega\subfine^*$. In the context of our ASM method, equation \cref{eq:quasi_newton} is solved by a quasi-Newton iteration with Broyden's method in complete analogy to \cref{ssec:general_sm}. The resulting method is summarized in \cref{algo:volume_asm,algo:volume_loop}, which detail our limited-memory implementation of the method.

\begin{algorithm2e}[!t]
	\KwIn{Reference geometry $\hat{\Omega}$, tolerance $\tau \in (0,1)$, maximum number of iterations $k_{\text{max}}$, memory size $m \geq 0$}
	Compute $\Omega\subcoarse^* = \argmin_{\Omega\in \admissiblegeom}\smallspace \reducedcostfunctional(\coarseresponse(\Omega))$ and define $\vectorfield^* = \numretraction_{\hat{\Omega}}^{-1} \left( \Omega\subcoarse^* \right) = \Omega\subcoarse^* - \hat{\Omega} \in H^1(\hat{\Omega};\R^d)$ \\
	Set $\Omega^0 = \Omega\subcoarse^*$\\
	\For{$k=0,1,2,\dots, k_\textrm{max}$}{
		Compute $\spacemappingfunction(\Omega^k) \in H^1(\Omega^k;\R^d)$ \\
		Set $\mathcal{W}^k = \numretraction_{\hat{\Omega}}^{-1} (\Omega^k) = \Omega^k - \hat{\Omega} \in H^1(\hat{\Omega};\R^d)$ \\
		\If{$\norm{\spacemappingfunction(\Omega^k) - \numtransport_{\mathcal{W}^k} (\vectorfield^*)}{\Omega^k} \leq \tau \norm{\vectortransport_{\mathcal{W}^k} (\vectorfield^*)}{\Omega^k}$}{
			\KwRet{\upshape $\Omega\subfine^* = \Omega^k$ as (approximate) solution}}
		Compute $\mathcal{X}^k = -(B^k)^{-1} \left( \spacemappingfunction(\Omega^k) - \numtransport_{\mathcal{W}^k} (\vectorfield^*) \right) \in H^1(\Omega^k;\R^d)$ with \cref{algo:volume_loop} using the stored $\mathcal{S}^i$ and $\mathcal{Y}^i$ for $i=k-m, \dots, k-1$ \\
		Update the geometry by $\Omega^{k+1} = \numretraction_{\Omega^k} (\mathcal{X}^k)$ \\
		%
		%		\Comment{Update the approximation of $(B^k)^{-1}$}
		%
		Compute $\nu = (B^k)^{-1}\left( \numtransport_{\numretraction_{\Omega^{k+1}}^{-1} (\Omega^{k})} (\spacemappingfunction(\Omega^{k+1})) - \spacemappingfunction(\Omega^{k}) \right)$ with \cref{algo:volume_loop} using the stored $\mathcal{S}^i$ and $\mathcal{Y}^i$ for $i=k-m, \dots, k-1$ \\
		Store $\mathcal{S}^k = \frac{1}{a_{\Omega^k}\left( \mathcal{X}^k, \nu \right)}(\mathcal{X}^k - \nu) \in H^1(\Omega^k, \R^d)$ and $\mathcal{Y}^k = \mathcal{X}^k \in H^1(\Omega^k;\R^d)$ \\
		\If{$k > m$}{Discard $\mathcal{S}^{k-m}$ and $\mathcal{Y}^{k-m}$ from storage}
	}
	\caption{Volume-Based Aggressive Space Mapping (ASM) Method for Shape Optimization.}
	\label{algo:volume_asm}
\end{algorithm2e}

\begin{algorithm2e}[!t]
	\KwIn{Right-hand side $\mathcal{B} \in H^1(\Omega^k;\R^d)$ and stored elements $\mathcal{S}^i, \mathcal{Y}^i \in H^1(\Omega^i, \R^d)$ for $i=k-m, \dots, k-1$}
	Set $\mathcal{H} = \mathcal{B}$ \\
	\For{$i=k-m,\dots,k-1$}{
		Set $\mathcal{M}^i = \numretraction_{\Omega^i}^{-1} (\Omega^k) = \Omega^k - \Omega^i$ \\
		Compute $\theta = a_{\Omega^k} \left( \mathcal{H}, \numtransport_{\mathcal{M}^i} (\mathcal{Y}^i) \right)$ \\
		$\mathcal{H} = \mathcal{H} + \theta\smallspace \numtransport_{\mathcal{M}^i} (\mathcal{S}^i)$ \\
	}
	\KwRet{$\mathcal{H} = (B^k)^{-1} \mathcal{B}$}
	\caption{Limited memory loop to compute $\mathcal{H} = (B^k)^{-1} \mathcal{B}$}
	\label{algo:volume_loop}
\end{algorithm2e}

%\begin{remark}
%	
%	
%	This is a major benefit for PDE constrained shape optimization: For the efficient direct solution of a shape optimization problem, sensitivity information is required so that gradient-based solvers for the optimization problem can be used. This sensitivity information is typically either derived by hand or with the help of automatic differentiation, usually utilizing an adjoint approach. Since deriving the necessary equations is tedious and error-prone, manually carrying out the adjoint approach is limited to models that are not overly complex. On the other hand, in order to apply automatic differentiation to derive the sensitivity information, the source code of the solver is required, which makes it not applicable when using commercial solvers. 
%	
%	For very complex problems such as the ones arising from industrial applications, commercial solvers are usually employed as lots of practically relevant models (e.g.\ turbulence models) and sophisticated solvers for these problems are available. Solving a shape optimization problem efficiently in such a setting is not practically possible using the adjoint approach by hand or with the help of automatic differentiation. However, the space mapping technique proposed in this paper (cf. \cref{sec:space_mapping}) now enables the shape optimization in such a complex setting.
%\end{remark}

\section{Numerical Solution of Shape Optimization Problems with Space Mapping Methods}
\label{sec:numerics}

Let us now numerically investigate the aggressive space mapping (ASM) method, which we proposed in \cref{sec:space_mapping}. We consider two model problems. The first is given by a shape identification problem constrained by a semi-linear transmission problem and the second is given by a problem of uniform flow distribution in a network of pipes. Moreover, we show that space mapping methods can be used to couple closed-source commercial solvers for the fine model with optimization software for the coarse model. Note that our implementation of the space mapping methods and our numerical experiments is available freely on GitHub \cite{Blauth2022Software}.

\subsection{Implementation Details}
\label{ssec:implementation_details}

We have implemented the space mapping methods for PDE constrained shape optimization in our software package cashocs \cite{Blauth2021cashocs}, which is an open-source software for solving PDE constrained optimal control and shape optimization problems. It is based on the finite element software FEniCS \cite{Logg2012Automated, Alnes2015FEniCS} and, hence, allows the user to implement their problems in the near-mathematical syntax of the Unified Form Language UFL \cite{Logg2012Automated}. Our software derives the corresponding adjoint systems and (shape) derivatives automatically with the help of automatic differentiation and it offers a number of optimization algorithms for the user to choose from. For shape optimization, cashocs only uses the volume formulation of the shape derivative, as it is shown in \cite{Hiptmair2015Comparison} that the volume formulation yields better approximation properties than the boundary formulation in the context of finite element methods. 

Let us first describe how we discretize the retraction $\numretraction$ and vector transport $\numtransport$ from \cref{ssec:volume_sm}. We proceed analogously to \cite{Schulz2016Efficient, Blauth2021Nonlinear} and use the following numerical retraction
\begin{equation}
	\label{eq:num_retraction}
	\numretraction_\Omega (\vectorfield) = (I + \vectorfield)\Omega = \Set{x + \vectorfield(x) | x\in \Omega},
\end{equation}
which is related to the perturbation of identity (cf. \cite{Blauth2021Nonlinear}). \Cref{eq:num_retraction} implies that we consider a Eulerian setting, where the underlying mesh is deformed in each iteration of the respective optimization algorithm. As we discretize the deformation vector fields with piecewise linear Lagrange elements, we can easily realize this by simply adding the values of the vector field $\vectorfield$ to the mesh nodes. Note that we obtain the inverse of $\numretraction$ by simply subtracting the nodal coordinates of the meshes corresponding to a deformed domain $\tilde{\Omega}$ and the reference domain $\Omega$, which we, again, denote as $\vectorfield = \numretraction_{\Omega}^{-1} (\tilde{\Omega}) =  \tilde{\Omega} - \Omega$.

Moreover, we use
\begin{equation}
	\label{eq:num_transport}
	\numtransport_{\vectorfield} \mathcal{W}(y) = \mathcal{W}(x) \qquad \text{ for } y = x + \vectorfield(x) \in \numretraction_\Omega (\vectorfield) \text{ with } x \in \Omega,
\end{equation}
as numerical approximation of the vector transport. Here, $\numtransport_{\vectorfield} \mathcal{W} \colon \numretraction_\Omega(\vectorfield) \to \R^d$ is a vector field on $\numretraction_\Omega(\vectorfield)$, and $\vectorfield \colon \Omega \to \R^d$ and $\mathcal{W} \colon \Omega \to \R^d$ are vector fields on $\Omega$. From the numerical point of view, a vector field $\mathcal{W}$ defined on $\Omega$ is represented by its vector of nodal values. \Cref{eq:num_transport} then states that the transported vector field $\numtransport_\vectorfield \mathcal{W}$ on the deformed domain $\numretraction_\Omega(\vectorfield)$ is represented by the same vector of nodal values as the original vector field $\mathcal{W}$, where only the position of the corresponding mesh nodes is changed according to \cref{eq:num_retraction}. Note that the retraction given by \cref{eq:num_retraction} and the vector transport given by \cref{eq:num_transport} are used in \cite{Schulz2016Efficient,Schulz2015Structured,Blauth2021Nonlinear} for the numerical realization of the Riemannian framework from \cref{ssec:stekpoinc}, and we refer the reader to these publications for further details. 

Finally, we use the equations of linear elasticity to compute the gradient deformation, i.e., for the bilinear form $a_\Omega$ we use
\begin{equation*}
	a_\Omega(V, W) = \integral{\Omega} 2 \mu\ \varepsilon(V) : \varepsilon(W) + \lambda\ \text{div}\left( V \right) \text{div}\left( W \right) + \delta\ V \cdot W \dmeas{x},
\end{equation*}
which shows good deformation properties numerically (cf.~\cite{Blauth2021Nonlinear,Schulz2016Computational,Gangl2015Shape}). Throughout this section, we use $\mu = 1$, $\lambda=0$, and $\delta=0$ for all numerical experiments.
%Note that our approach is easy to realize while still showing good approximation properties in comparison to more sophisticated discretizations of the retraction and vector transport (cf.~\cite{Schulz2015Structured}).

\subsection{A Semi-Linear Transmission Problem}
\label{ssec:semi_linear_transmission}

In this section, we consider a semi-linear version of the transmission problem, which we encountered earlier in \cref{ssec:shape_calculus}. Here, we consider the same setting as in \cref{ex:transmission} before, which is sketched in \cref{fig:model_setup}.
%Let $\holdall \subset \R^d$ be a fixed, open, and bounded hold-all domain with Lipschitz boundary $\partial\holdall$, which is partitioned into two disjoint open subsets: an inclusion $\Omega \subset \holdall$ and an exterior part $\Omega\subout := \holdall \setminus \closure{\Omega}$. The boundary of the interior part is denoted by $\Gamma = \partial \Omega$, which is assumed to be smooth, and the outer unit normal vector on $\Gamma$ is denoted by $\normal$. This setting is sketched in \cref{fig:model_setup}.

We consider the following two models for the ASM method. The fine model is given by the semi-linear transmission problem
\begin{equation}
	\label{eq:semi_linear_transmission_model}
	\begin{alignedat}[t]{2}
		- \nabla \cdot (\interfacecoefficient\smallspace \grad u) + \beta u^3 &= f \quad &&\text{ in } D, \\
		u &= 0 \quad &&\text{ on } \partial D, \\
		\jump{u}{\Gamma} &= 0 \quad &&\text{ on } \Gamma,\\
		\jump{\interfacecoefficient\smallspace \partial_\normal u}{\Gamma} &= 0 \quad &&\text{ on } \Gamma,
	\end{alignedat}
\end{equation}
and the coarse model is given by the linear transmission problem
\begin{equation}
	\label{eq:linear_transmission_model}
	\begin{alignedat}[t]{2}
		- \nabla \cdot (\interfacecoefficient\smallspace \grad u) &= f \quad &&\text{ in } D, \\
		u &= 0 \quad &&\text{ on } \partial D, \\
		\jump{u}{\Gamma} &= 0 \quad &&\text{ on } \Gamma,\\
		\jump{\interfacecoefficient\smallspace \partial_\normal u}{\Gamma} &= 0 \quad &&\text{ on } \Gamma.
	\end{alignedat}
\end{equation}
For both models, both $\interfacecoefficient$ and $f$ are constant in each subdomain, i.e.,
\begin{equation*}
	\interfacecoefficient = \begin{cases}
	\interfacecoefficient\subin = \text{const.} \quad &\text{ in } \Omega, \\
	\interfacecoefficient\subout = \text{const.} \quad &\text{ in } \Omega\subout,
	\end{cases}  \qquad f = \begin{cases}
	f\subin = \text{const.} \quad &\text{ in } \Omega, \\
	f\subout = \text{const.} \quad &\text{ in } \Omega\subout.
	\end{cases}
\end{equation*}
We consider the following cost functional
\begin{equation*}
	\costfunctional(\Omega, u) = \integral{\holdall} \left( u - u\subdes \right)^2 \dmeas{x},
\end{equation*} 
where $u\subdes \in H^1(\holdall)$ is some desired state. We have the optimization problems
\begin{equation}
	\label{eq:semi_linear_transmission}
	\min_{\Omega \in \admissiblegeom} \costfunctional(\Omega, u) \quad \text{ subject to } \cref{eq:semi_linear_transmission_model},
\end{equation}
which acts as the fine model optimization problem, as well as
\begin{equation}
	\label{eq:linear_transmission}
	\min_{\Omega\in \admissiblegeom} \costfunctional(\Omega, u) \quad \text{ subject to } \cref{eq:linear_transmission_model},
\end{equation}
which acts as the coarse model optimization problem. Since the coarse model \cref{eq:linear_transmission_model} is linear, problem \cref{eq:linear_transmission} is substantially easier to solve than \cref{eq:semi_linear_transmission}, which is constrained by the nonlinear PDE \cref{eq:semi_linear_transmission_model}.

\begin{figure}[!b]
	\centering
	\includegraphics[width=0.275\textwidth]{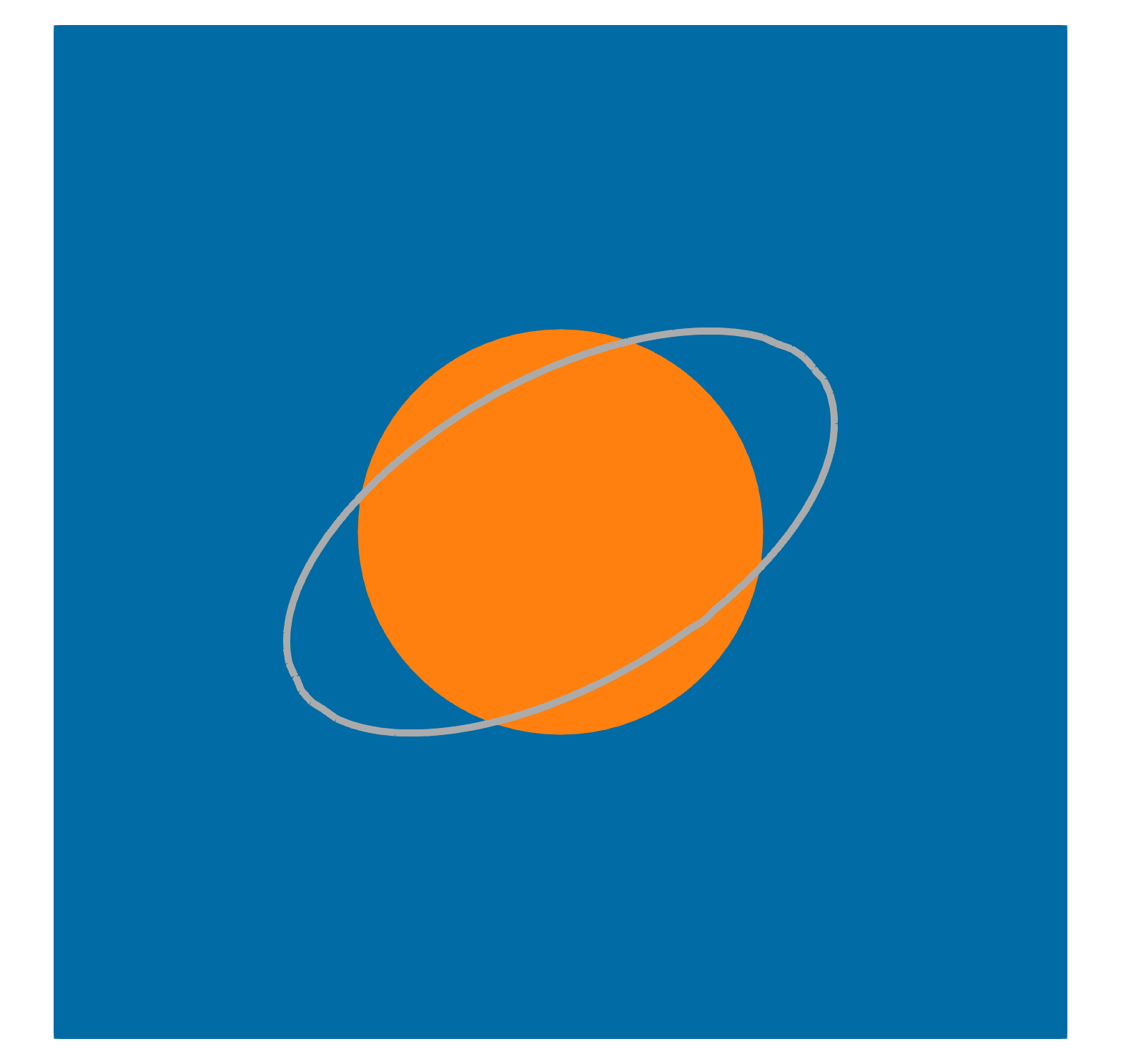}
	\caption{Initial inclusion $\Omega$ (orange) and reference ellipse (outlined in gray).}
	\label{fig:initial_setup_transmission}
\end{figure}

We solve \cref{eq:semi_linear_transmission} with the previously-proposed ASM method. Therefore, we define the fine model response by $\fineresponse \colon \admissiblegeom \to H^1_0(\holdall);\ \Omega \mapsto u\subfine(\Omega)$, where $u\subfine(\Omega)$ is the solution of \cref{eq:semi_linear_transmission_model} for the inclusion $\Omega$, and the coarse model response is given by $\coarseresponse \colon \admissiblegeom \to H^1_0(\holdall); \Omega \mapsto u\subcoarse(\Omega)$, where $u\subcoarse(\Omega)$ is the solution of \cref{eq:linear_transmission_model} for the inclusion $\Omega$. We define the misalignment function $\misalignmentfunction \colon H^1_0(\holdall) \times H^1_0(\holdall) \to \R$ as
\begin{equation*}
	\misalignmentfunction(u, v) = \integral{\holdall} \left( u - v \right)^2 \dmeas{x}.
\end{equation*}
Hence, to evaluate the space mapping function $\spacemappingfunction$, we have to solve the following subproblem
\begin{equation*}
	\tilde{\vectorfield} = \argmin_{\vectorfield \in H^1(\hat{\Omega};\R^d)} \misalignmentfunction\left( \coarseresponse\left( \numretraction_{\hat{\Omega}}(\vectorfield) \right), \fineresponse(\Omega\subfine) \right),
\end{equation*}
which is, due to the definition of the misalignment function, equivalent to
\begin{equation}
	\label{eq:parameter_extraction}
	\tilde{\Omega} = \argmin_{\Omega \in \admissiblegeom} \integral{\holdall} \left( u\subcoarse(\Omega) - u\subfine(\Omega\subfine) \right)^2 \dmeas{x},
\end{equation}
where the corresponding minimizers are related by $\tilde{\Omega} = \numretraction_{\hat{\Omega}}(\tilde{\vectorfield})$. Note that \cref{eq:parameter_extraction} is equivalent to the coarse model problem \cref{eq:linear_transmission} with $u\subdes$ replaced by $u\subfine(\Omega\subfine)$.

%\begin{Remark}
%	Note that problem \cref{eq:parameter_extraction} is very similar to the coarse model optimization problem \cref{eq:linear_transmission}, with $u\subdes$ replaced by $u\subfine(\Omega\subfine)$. This is the case due to the definition of the misalignment function and allows the user to reuse the same code for the optimization of the coarse model and the computation of the space mapping function, with $u\subdes$ replaced by $u\subfine(\Omega\subfine)$.
%\end{Remark}

\begin{figure}[!b]
	\centering
	\begin{subfigure}{0.5\textwidth}
		\centering
		\includegraphics[width=0.8\textwidth]{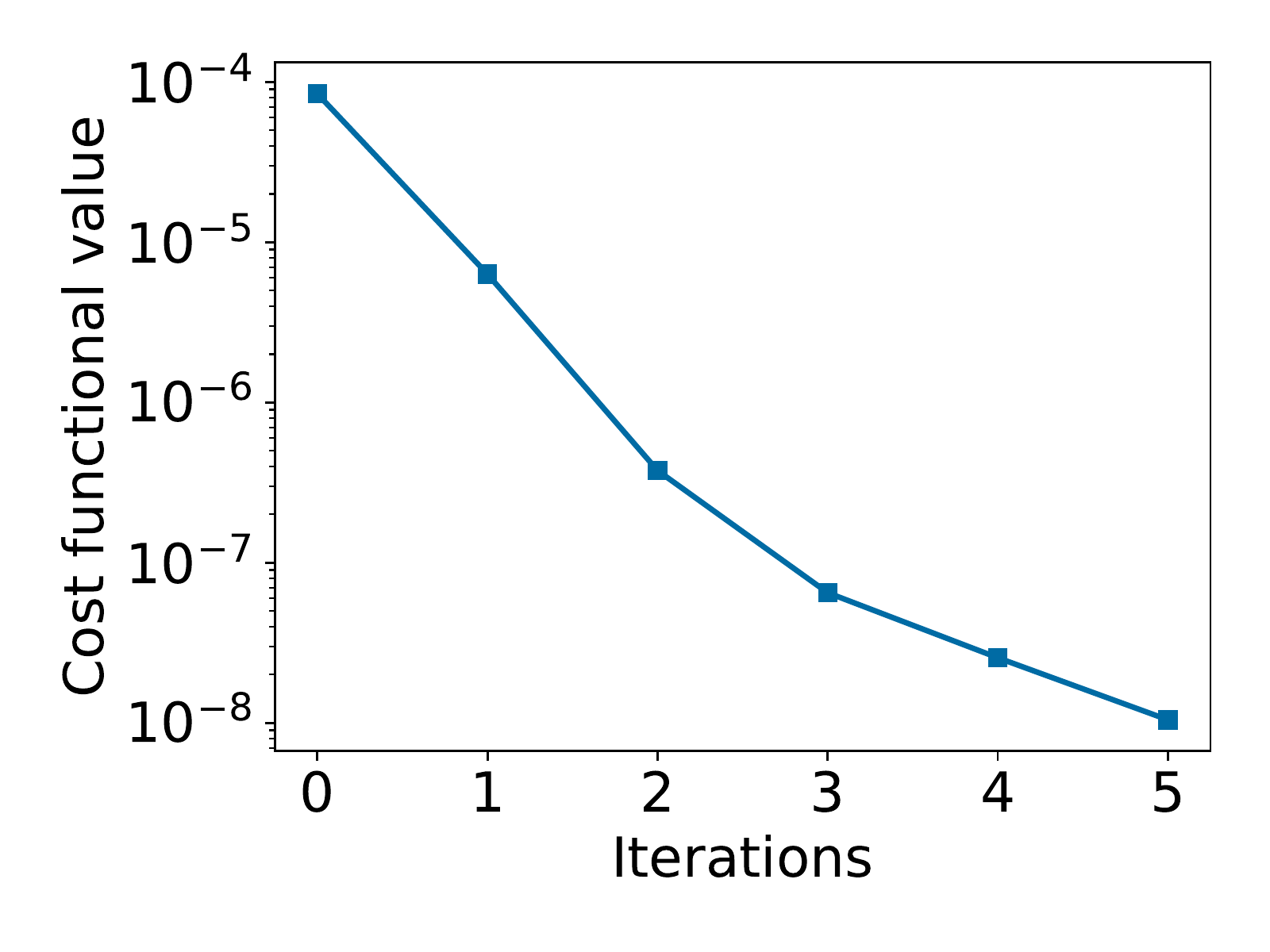}
		\caption{Evolution of the cost functional.}
	\end{subfigure}%
	\hfil%
	\begin{subfigure}{0.5\textwidth}
		\centering
		\includegraphics[width=0.8\textwidth]{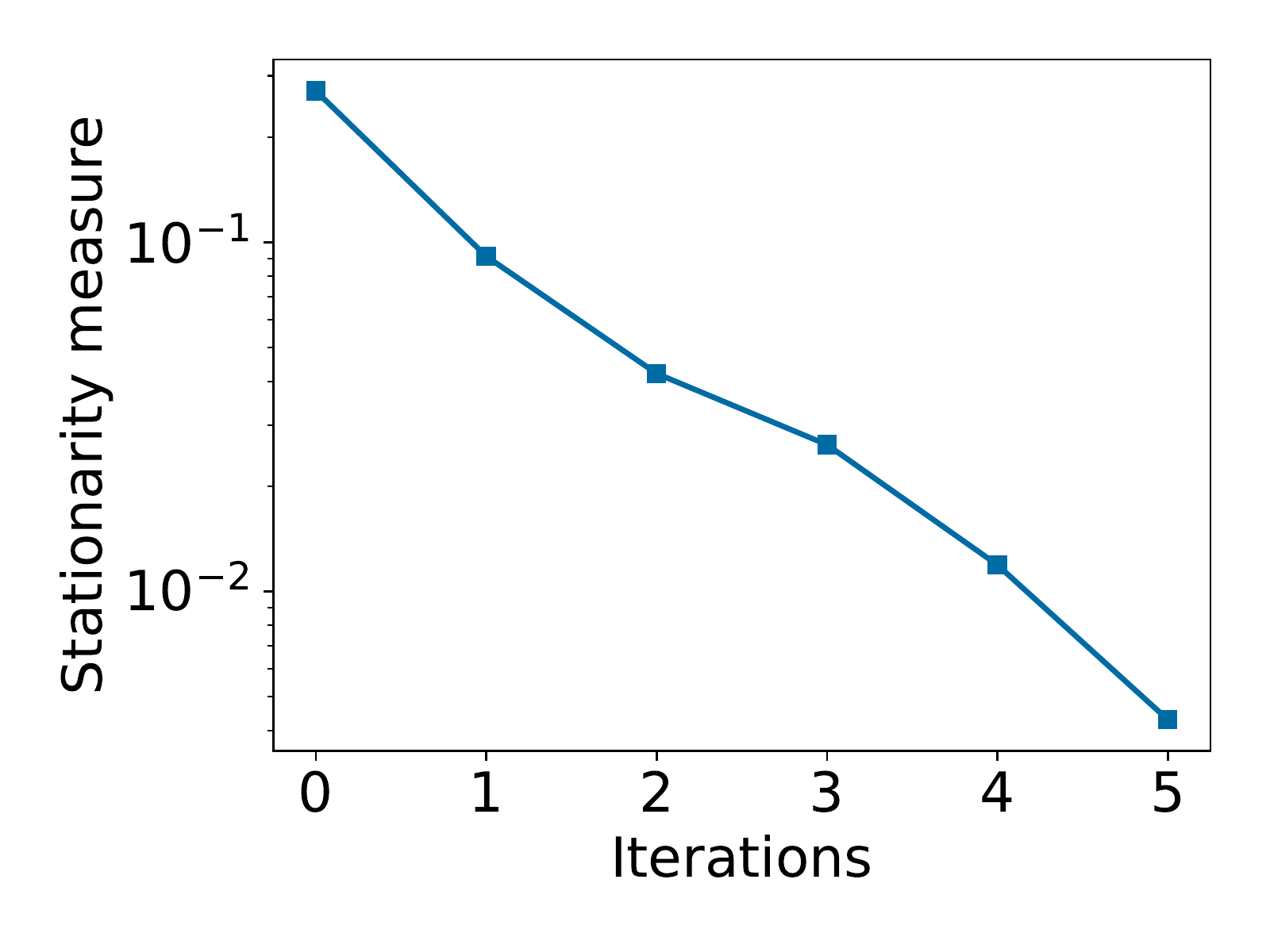}
		\caption{Evolution of the stationarity measure.}
	\end{subfigure}
	\caption{History of the ASM method for problem \cref{eq:semi_linear_transmission}.}
	\label{fig:asm_transmission}
\end{figure}

\begin{figure}[!b]
	\centering
	\begin{subfigure}{0.275\textwidth}
		\centering
		\includegraphics[width=\textwidth]{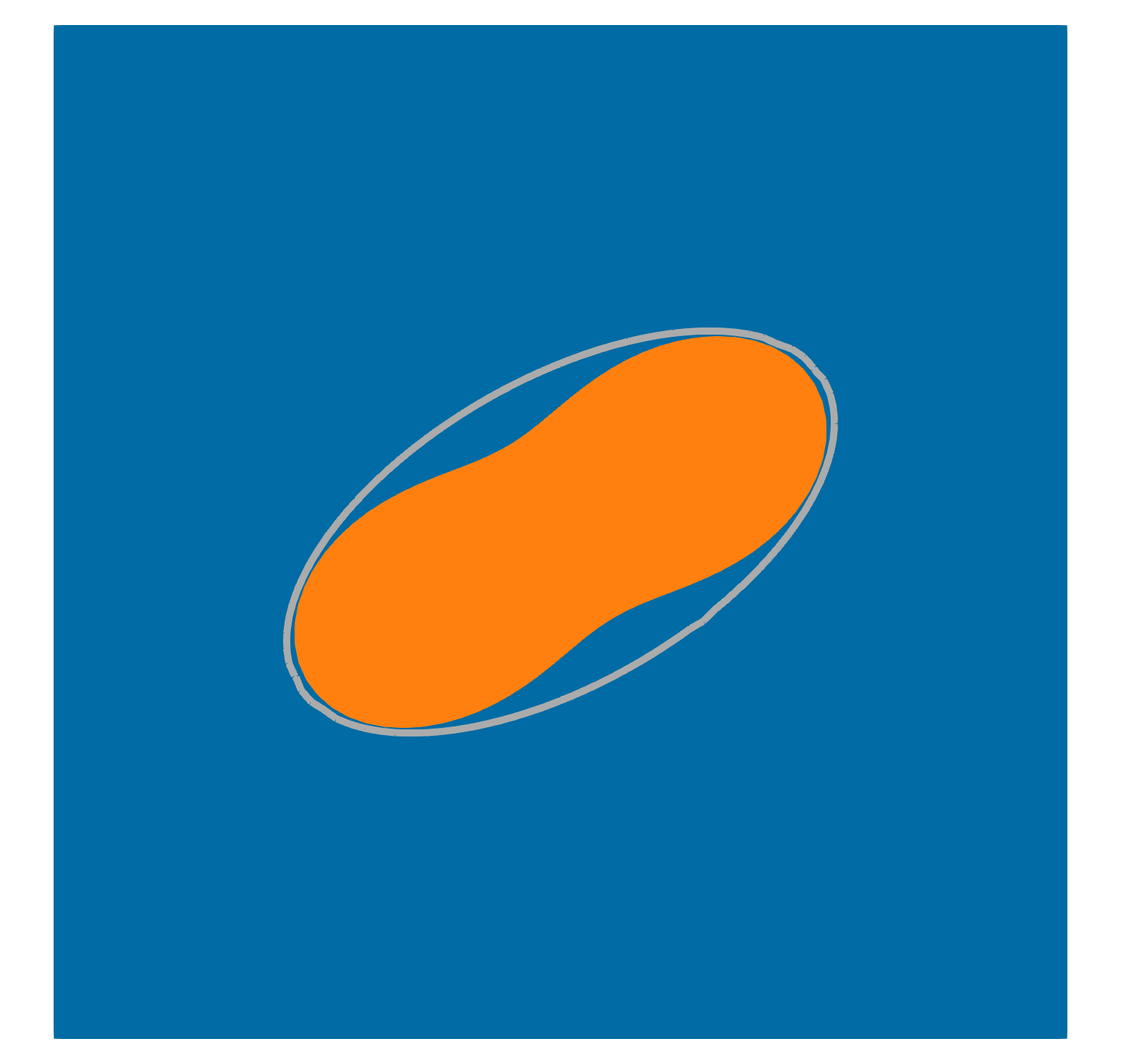}
		\caption{Iteration 0.}
		\label{fig:trans_0}
	\end{subfigure}%
	\hfil%
	\begin{subfigure}{0.275\textwidth}
		\centering
		\includegraphics[width=\textwidth]{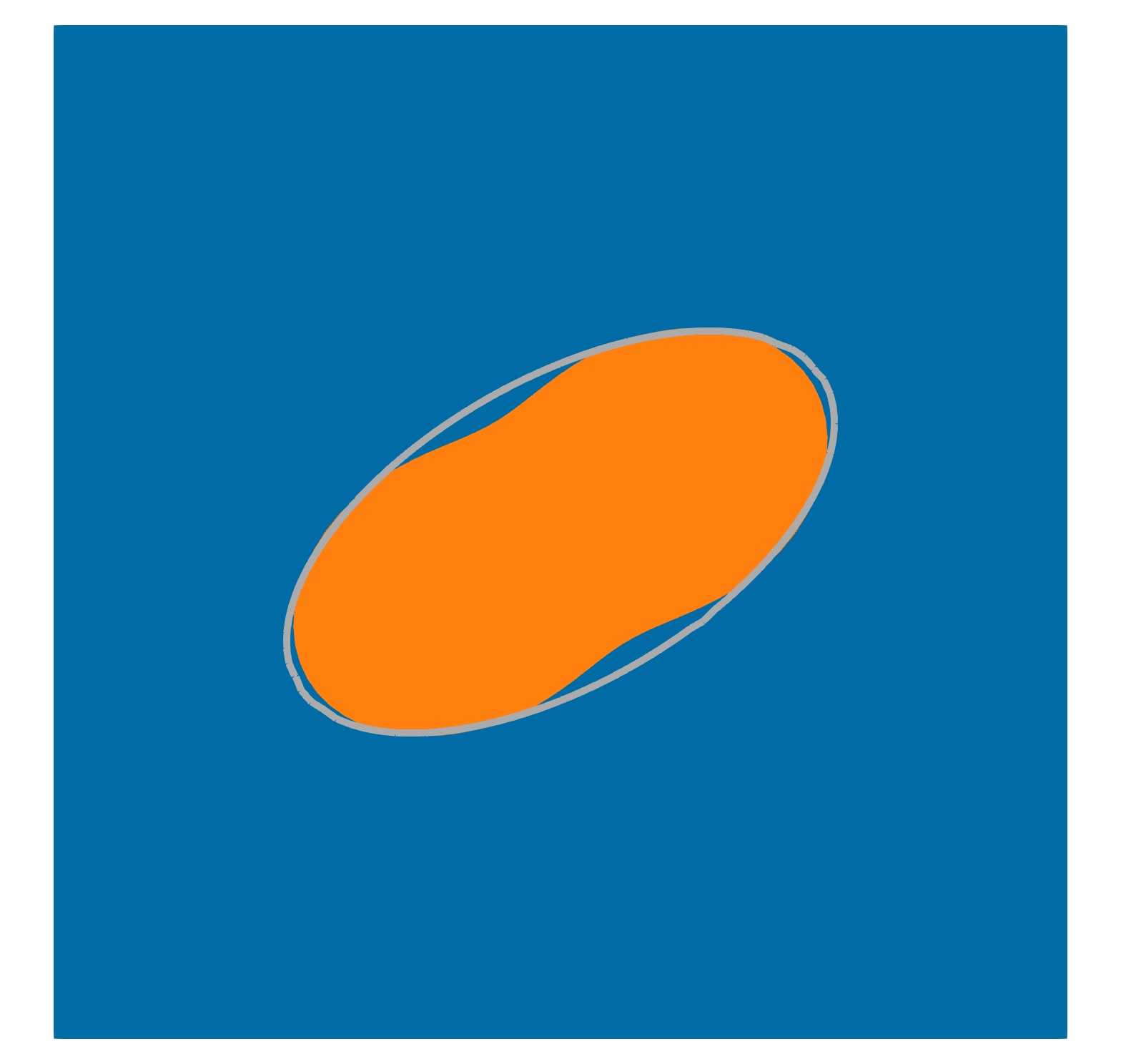}
		\caption{Iteration 1.}
	\end{subfigure}%
	\hfil%
	\begin{subfigure}{0.275\textwidth}
		\centering
		\includegraphics[width=\textwidth]{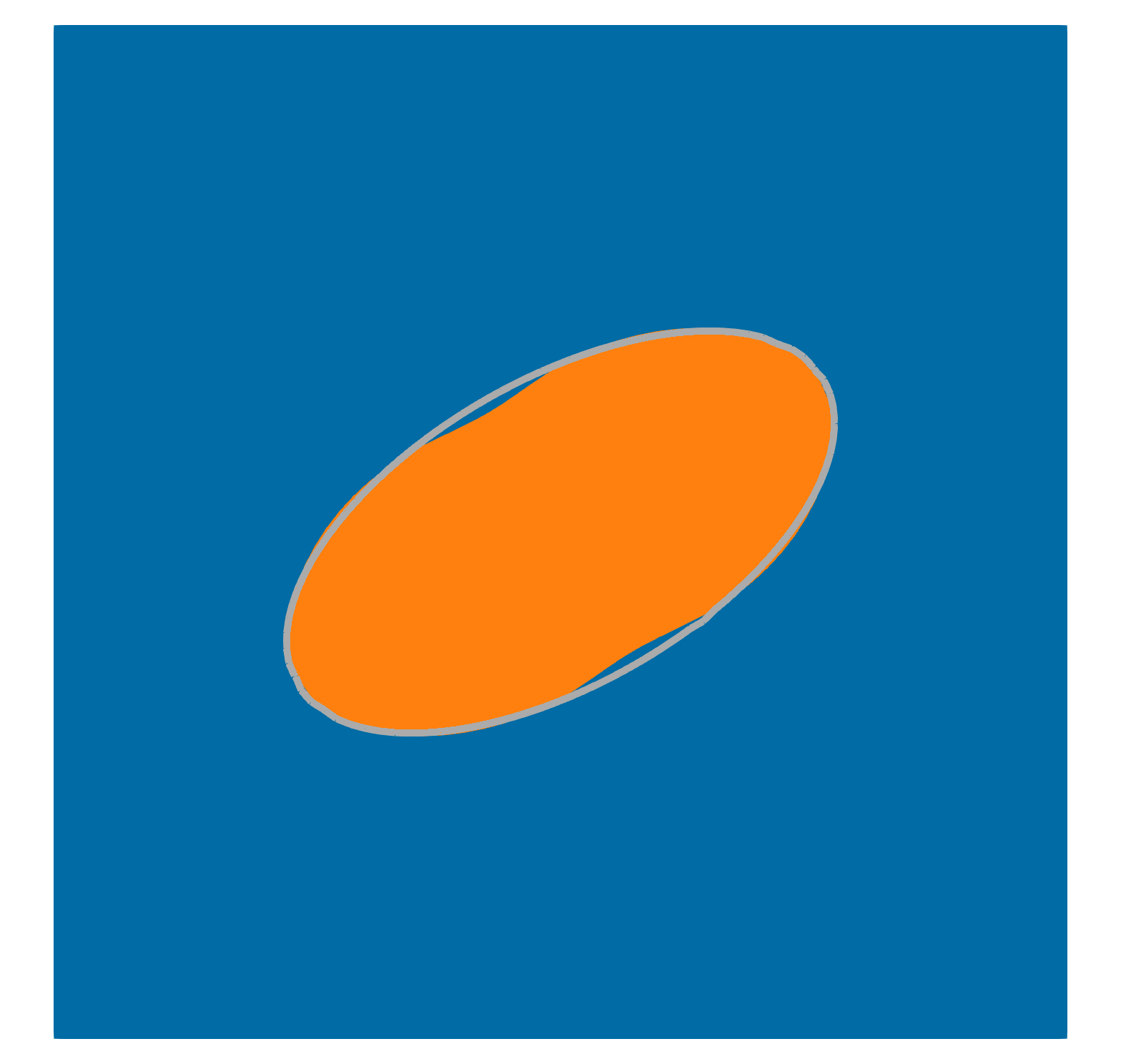}
		\caption{Iteration 2.}
	\end{subfigure}
	
	\begin{subfigure}{0.275\textwidth}
		\centering
		\includegraphics[width=\textwidth]{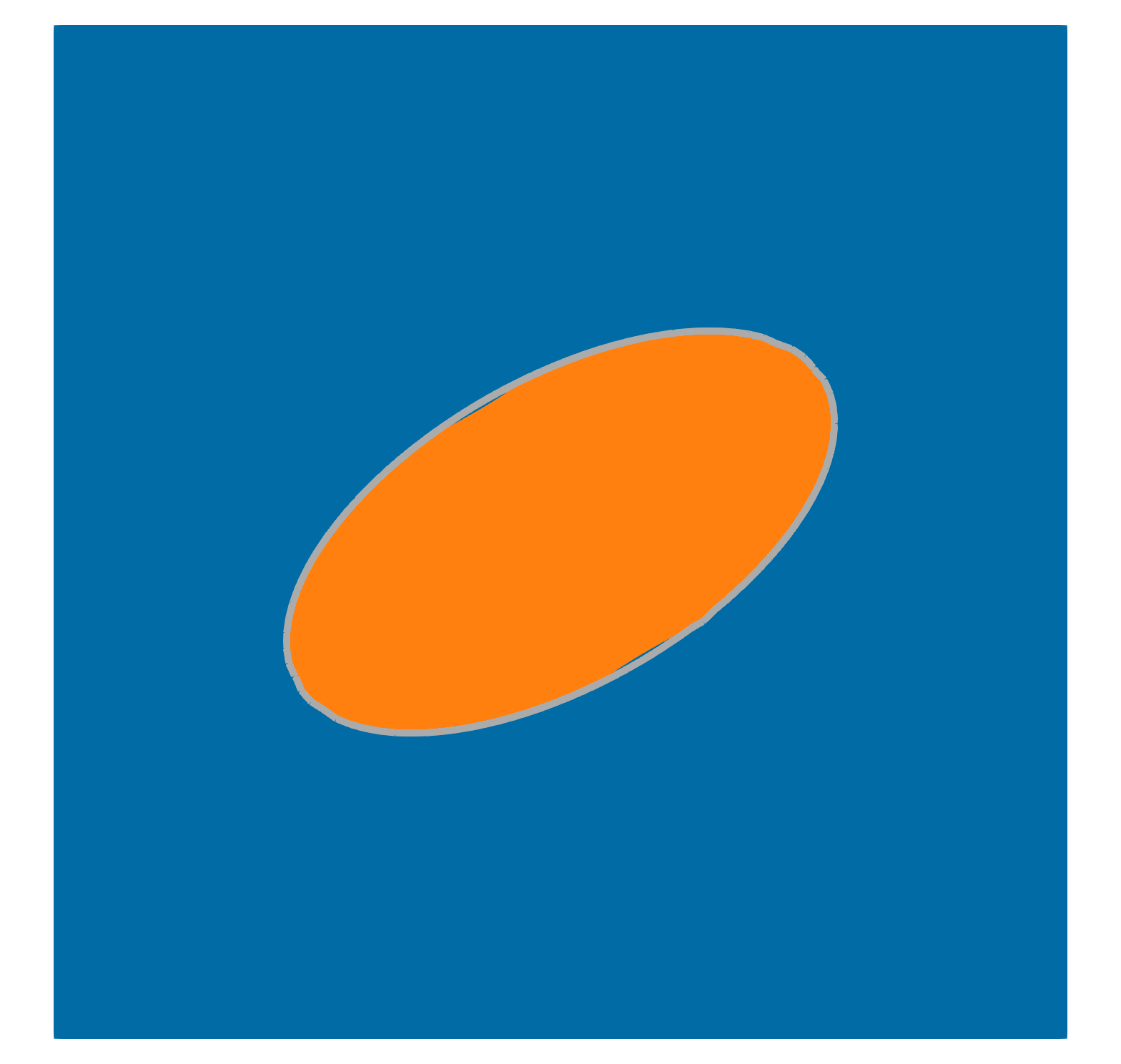}
		\caption{Iteration 3.}
	\end{subfigure}%
	\hfil%
	\begin{subfigure}{0.275\textwidth}
		\centering
		\includegraphics[width=\textwidth]{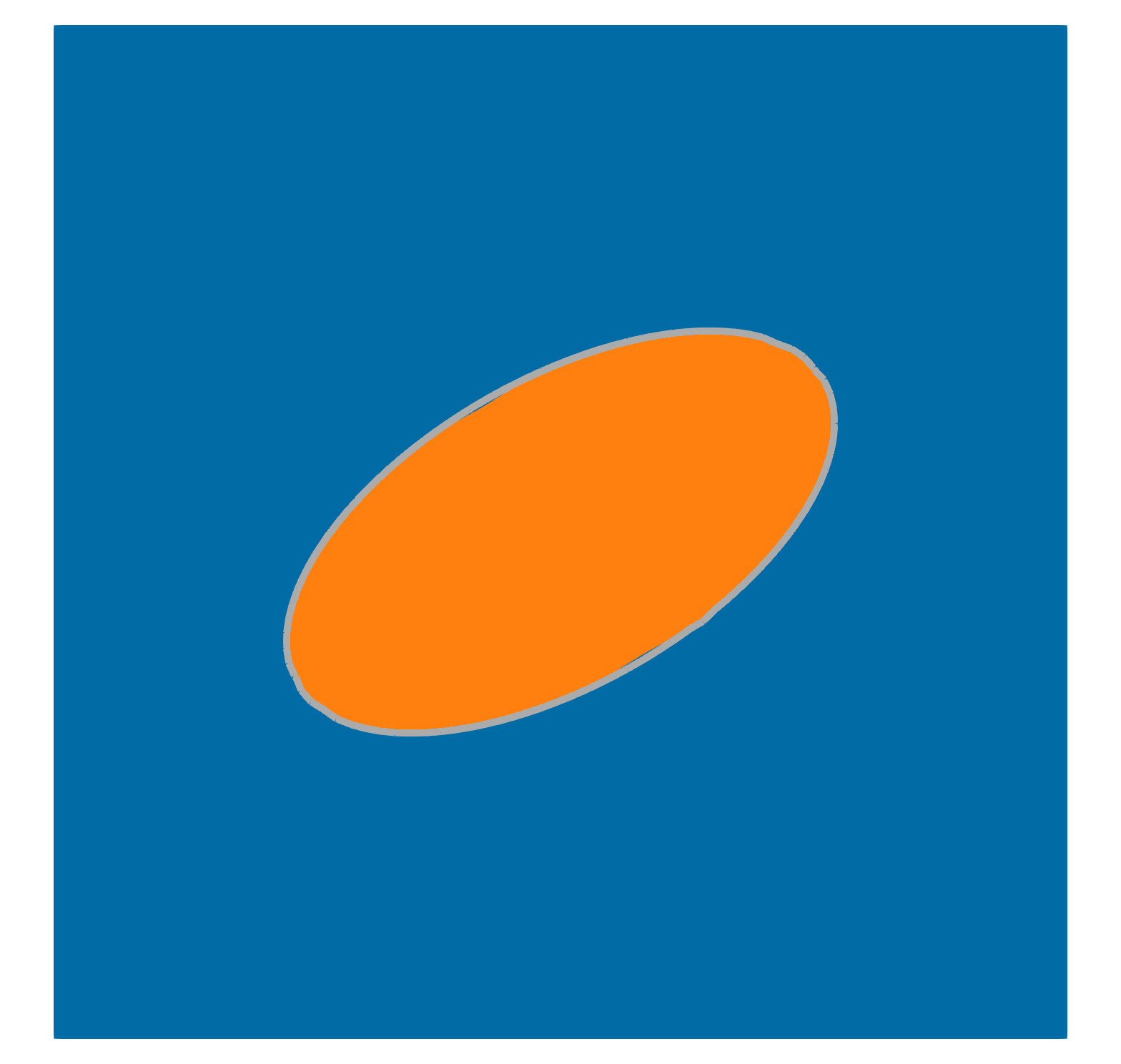}
		\caption{Iteration 4.}
	\end{subfigure}%
	\hfil%
	\begin{subfigure}{0.275\textwidth}
		\centering
		\includegraphics[width=\textwidth]{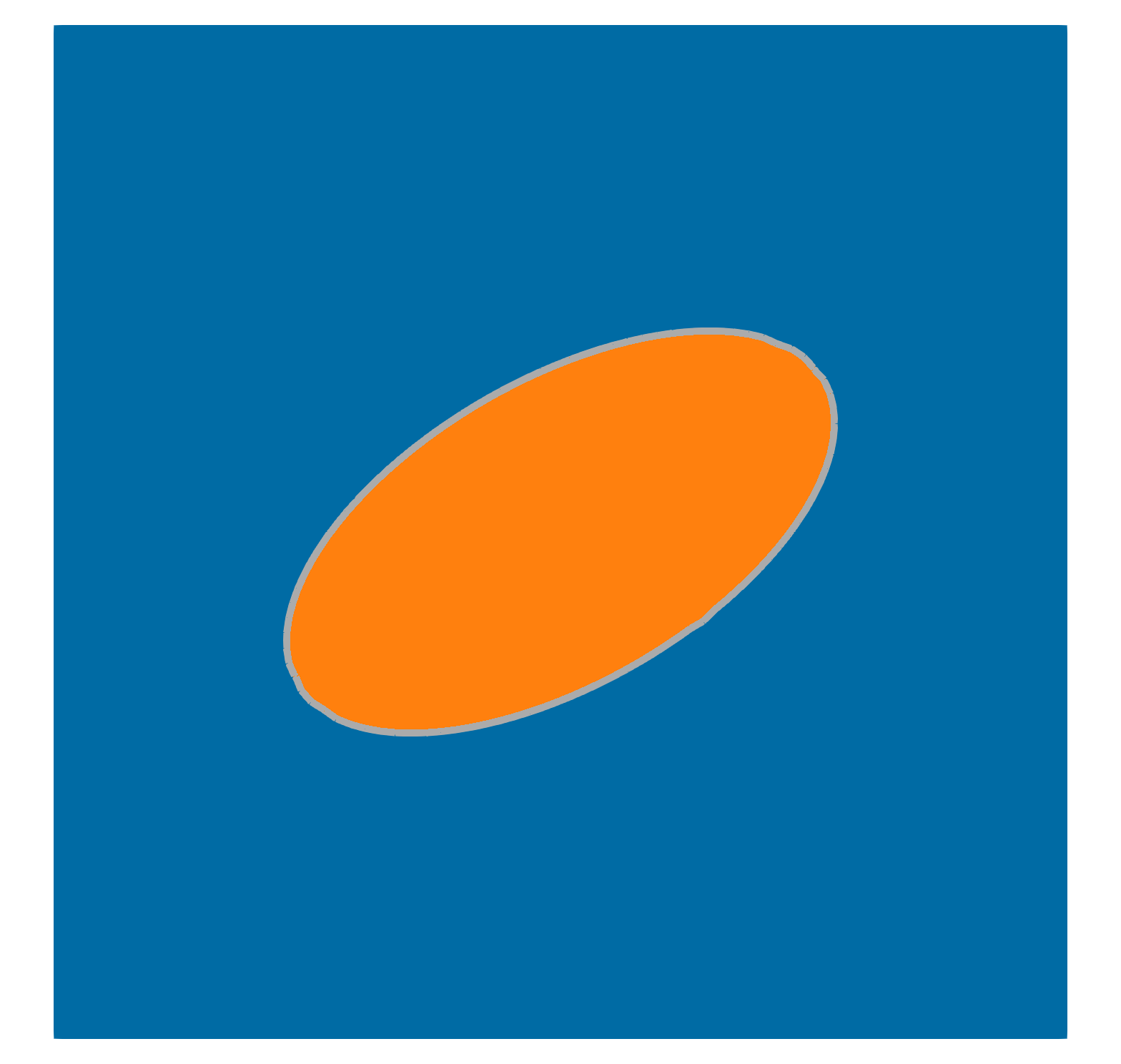}
		\caption{Iteration 5.}
	\end{subfigure}
	\caption{Evolution of the inclusion $\Omega$ over the course of the ASM method. The interior domain $\Omega$ is depicted in orange, the outer domain $\holdall\setminus \closure{\Omega}$ in blue, and the reference ellipse is outlined in gray.}
	\label{fig:transmission_geometries}
\end{figure}

Let us now describe our numerical setting. The hold-all domain is chosen as the unit square $\holdall = (0,1)^2$, which we discretize with a mesh consisting of \num{5398} nodes and \num{10526} triangles. Particularly, we use the same discretization for both the fine and the coarse model. To solve problem \cref{eq:semi_linear_transmission}, we utilize the ASM method from \cref{sec:space_mapping}, which is implemented in our software package cashocs \cite{Blauth2021cashocs}. For the solution of the state and adjoint systems, the finite element software FEniCS \cite{Logg2012Automated} is used, where both state and adjoint variables are discretized with piecewise linear Lagrange elements. Moreover, the coarse model shape optimization problems \cref{eq:linear_transmission,eq:parameter_extraction} are solved with a limited memory BFGS method (cf.~\cite{Schulz2016Efficient}) with relative tolerance \num{1e-2}, which is also implemented in cashocs. For the jumping coefficient $\alpha$, we choose $\alpha\subin = 10$ and $\alpha\subout = 1$, for the right-hand side $f$ we use $f\subin = 10$ and $f\subout = 1$, and for the coefficient $\beta$ in \cref{eq:semi_linear_transmission_model} we use $\beta = 100$. We obtain the desired state $u\subdes$ by solving the fine model \cref{eq:semi_linear_transmission_model} for a reference inclusion $\Omega$ given by an ellipse with center $(0.5, 0.5)$, a semi-major axis of \num{0.3}, a semi-minor axis of \num{0.15}, which is rotated counterclockwise by \SI{30}{\degree} (cf.~\cref{fig:transmission_geometries}). For the reference geometry $\hat{\Omega}$ we choose an inclusion $\Omega$ given by a circle with center $(0.5, 0.5)$ and radius $0.2$, which also acts as initial guess for all optimization problems. This setup can be seen in \cref{fig:initial_setup_transmission}. We stop the ASM method once the stationarity measure $\sigma$, which is defined as
\begin{equation}
	\label{eq:stationarity_measure}
	\sigma = \frac{\norm{\spacemappingfunction(\Omega^k) - \numtransport_{\mathcal{W}^k} (\vectorfield^*)}{\Omega^k}}{ \norm{\vectortransport_{\mathcal{W}^k} (\vectorfield^*)}{\Omega^k}},
\end{equation}
is smaller than the relative tolerance of $\tau = \num{1e-2}$ (cf.~\cref{algo:volume_asm}).

The history of the ASM method can be seen in \cref{fig:asm_transmission}, where we depict the evolution of the cost functional and the stationarity measure. We observe a very fast reduction of the cost functional, which decreases by about four orders of magnitude after five space mapping iterations. Additionally, we see that the stationarity measure also decreases by an order of magnitude and that the relative convergence criterion is reached after five iterations. 
%These observations show the fast convergence and efficiency the space mapping method exhibits. 

The evolution of the inclusion $\Omega$ over the course of the space mapping method, which is depicted in \cref{fig:transmission_geometries}, reinforces these observations. Here, we observe that the initial guess of the ASM method, which is obtained by solving the coarse model optimization problem \cref{eq:linear_transmission_model}, is not such a good approximation of the desired ellipse (cf.~\cref{fig:trans_0}). Instead, it resembles more a \qe{peanut} shape and is too small, but the orientation of the shape is already correct. During the space mapping iterations, the shape approximates the ellipse better in each iteration, as suggested by the steep decrease of the cost functional we discussed previously. After the third iteration, there are only very subtle differences between $\Omega$ and the desired ellipse, and in the final iteration there are no visible differences anymore. Hence, the ASM method is successful in solving the shape identification problem \cref{eq:semi_linear_transmission}. In particular, we only had to solve the nonlinear PDE \cref{eq:semi_linear_transmission} six times during the ASM method. Additionally, we had to solve the coarse model optimization problem six times, for which we only had to solve linear PDEs, which substantially simplifies the effort of solving the fine model optimization problem \cref{eq:semi_linear_transmission}.

\subsection{Uniform Flow Distribution}
\label{ssec:uniform_flow_distribution}

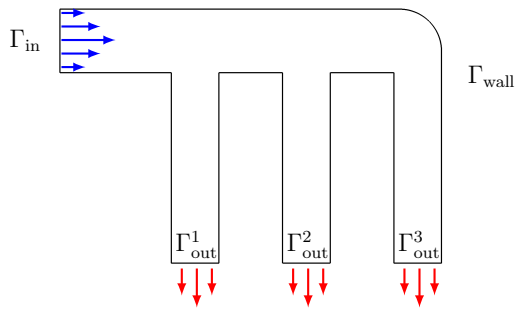
\begin{figure}[!b]
	\centering
	\scalebox{0.8}{
		\begin{tikzpicture}
		\node at (0, 0) {\includegraphics[width=0.5\textwidth, trim=0cm 0cm 3cm 0cm, clip]{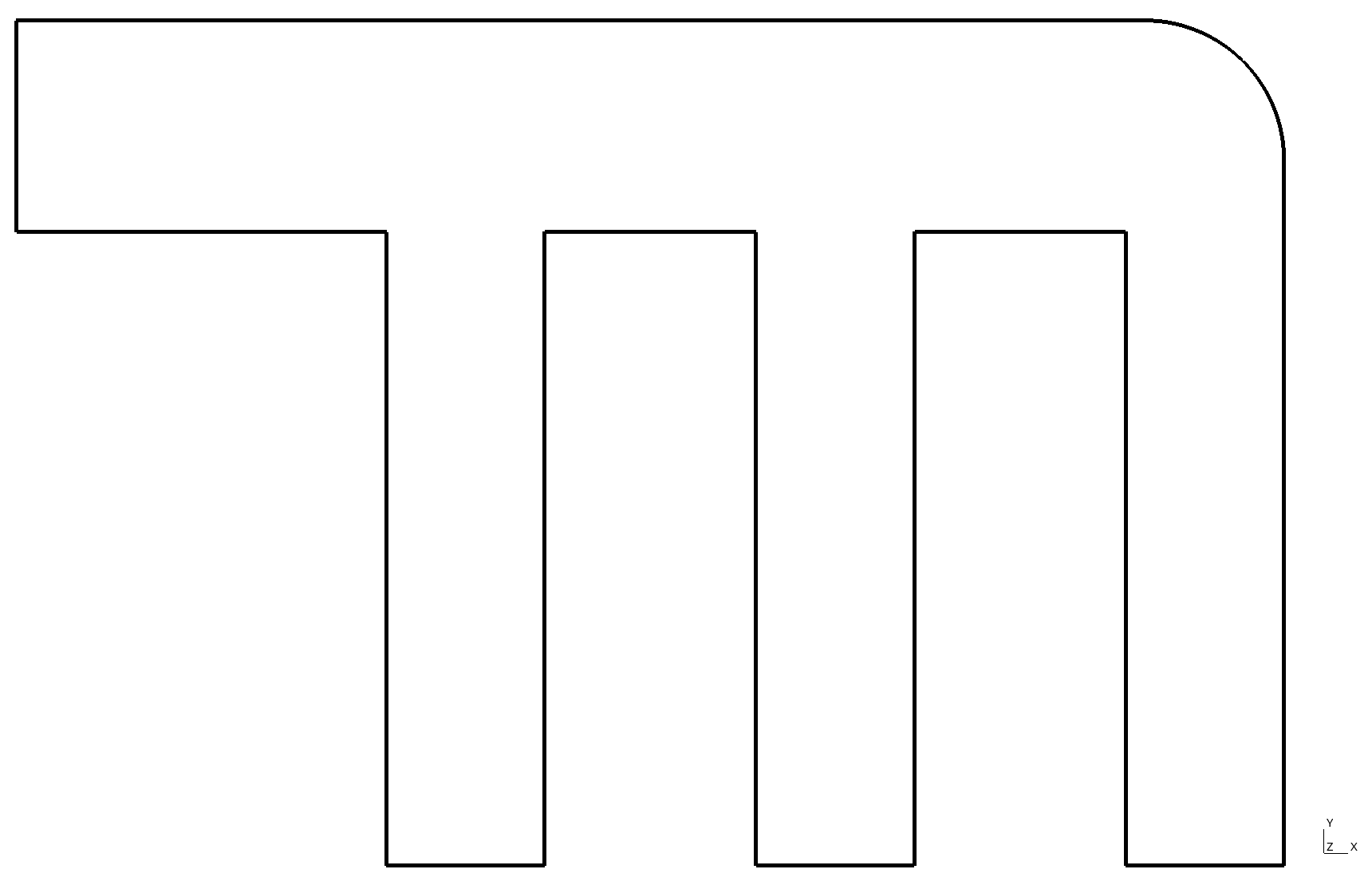}};
		
		\draw[-latex, line width=1,color=blue] (-3.15,2.05) -- (-2.75,2.05);
		\draw[-latex, line width=1,color=blue] (-3.15,1.825) -- (-2.5,1.825);
		\draw[-latex, line width=1,color=blue] (-3.15,1.6) -- (-2.25,1.6);
		\draw[-latex, line width=1,color=blue] (-3.15,1.375) -- (-2.5,1.375);
		\draw[-latex, line width=1,color=blue] (-3.15,1.15) -- (-2.75,1.15);
		
		\draw[-latex, line width=1, color=red] (-1.15,-2.2) -- (-1.15, -2.65);
		\draw[-latex, line width=1, color=red] (-0.9,-2.2) -- (-0.9, -2.85);
		\draw[-latex, line width=1, color=red] (-0.65,-2.2) -- (-0.65, -2.65);
		
		\draw[-latex, line width=1, color=red] (0.7,-2.2) -- (0.7, -2.65);
		\draw[-latex, line width=1, color=red] (0.95,-2.2) -- (0.95, -2.85);
		\draw[-latex, line width=1, color=red] (1.2,-2.2) -- (1.2, -2.65);
		
		\draw[-latex, line width=1, color=red] (2.55,-2.2) -- (2.55, -2.65);
		\draw[-latex, line width=1, color=red] (2.8,-2.2) -- (2.8, -2.85);
		\draw[-latex, line width=1, color=red] (3.05,-2.2) -- (3.05, -2.65);
		
		\node at (-3.75,1.6) {\large $\Gamma\subin$};
		\node at (-0.9, -1.8) {\large $\Gamma\subout^1$};
		\node at (0.95, -1.8) {\large $\Gamma\subout^2$};
		\node at (2.8, -1.8) {\large $\Gamma\subout^3$};
		\node at (4, 1) {\large $\Gamma_{\mathrm{wall}}$};
		\end{tikzpicture}
	}
	\caption{Reference geometry for the uniform flow distribution.}
	\label{fig:schema_flow}
\end{figure}

For our second numerical example, we consider the optimization of the flow inside a pipe network: The flow is divided from one into three pipes, and we want to achieve a uniform flow distribution over all three outlet pipes. This setting is depicted schematically in \cref{fig:schema_flow}. We denote the geometry of the pipe network by $\Omega$ and its boundary $\Gamma = \partial \Omega$ is partitioned into three parts: the inlet $\Gamma\subin$, where the flow enters the domain, the wall $\Gamma_\mathrm{wall}$, where we have a no-slip condition, and the outlet $\Gamma\subout$, where the flow leaves the domain. In particular, we denote the outlet boundary of each pipe by $\Gamma\subout^i$ for $i=1,2,3$. 

We consider the following models for the space mapping method. The fine model is given by the stationary incompressible Navier-Stokes equations
\begin{equation}
	\label{eq:flow_fine}
	\begin{alignedat}{2}
		-\Delta u + \mathrm{Re} \left( u \cdot \grad \right) u + \grad p &= 0 \quad &&\text{ in } \Omega,\\
		\div{u} &= 0 \quad &&\text{ in } \Omega,\\
		u &= u\subin \quad &&\text{ on } \Gamma\subin,\\
		u &= 0 \quad &&\text{ on } \Gamma_{\mathrm{wall}}, \\
		p &= 0 \quad &&\text{ on } \Gamma\subout, \\
		u \times \normal &= 0 \quad &&\text{ on } \Gamma\subout,
	\end{alignedat}
\end{equation}
and the coarse model is given by the linear Stokes system
\begin{equation}
	\label{eq:flow_coarse}
	\begin{alignedat}{2}
		-\Delta u + \grad p &= 0 \quad &&\text{ in } \Omega,\\
		\div{u} &= 0 \quad &&\text{ in } \Omega,\\
		u &= u\subin \quad &&\text{ on } \Gamma\subin,\\
		u &= 0 \quad &&\text{ on } \Gamma_{\mathrm{wall}}, \\
		p &= 0 \quad &&\text{ on } \Gamma\subout, \\
		u \times \normal &= 0 \quad &&\text{ on } \Gamma\subout.
	\end{alignedat}
\end{equation}
Here, $u$ denotes the fluid's velocity, $p$ its pressure, and $\mathrm{Re}$ is the Reynolds number.

To obtain a uniform flow distribution among the three pipes, we consider the cost functional
\begin{equation*}
	\costfunctional(\Omega, u) = \frac{1}{2} \sum_{i=1}^{3} \left( q\subout^i(u) - q\subdes \right)^2,
\end{equation*}
where the outlet flow rate $q\subout^i(u)$ of pipe $i$ is defined as 
\begin{equation*}
	q\subout^i(u) = \integral{\Gamma\subout^i} u\cdot \normal \dmeas{s}.
\end{equation*}
Hence, by minimizing the cost functional $\costfunctional$, we aim at achieving a uniform outlet flow rate for each pipe. Due to the incompressibility of the fluid, we have $q\subdes = \nicefrac{q\subin}{3}$, where the inlet flow rate is defined by $q\subin = -\integral{\Gamma\subin} u\subin\cdot \normal \dmeas{s}$. Analogously to before, we have the following shape optimization problems. The fine model optimization problem is given by
\begin{equation}
	\label{eq:flow_fine_model}
	\min_{\Omega \in \admissiblegeom}\smallspace\costfunctional(\Omega, u) \quad \text{ subject to } \cref{eq:flow_fine},
\end{equation}
and the coarse model optimization problem reads
\begin{equation}
	\label{eq:flow_coarse_model}
	\min_{\Omega \in \admissiblegeom}\smallspace\costfunctional(\Omega, u) \quad \text{ subject to } \cref{eq:flow_coarse}.
\end{equation}
The fine model optimization problem \cref{eq:flow_fine_model} is constrained by the nonlinear Navier-Stokes equations, whereas the coarse model optimization problem \cref{eq:flow_coarse_model} is only constrained by the linear Stokes system. Hence, the numerical solution of \cref{eq:flow_coarse_model} is substantially easier than that of \cref{eq:flow_fine_model}.

\begin{figure}[!b]
	\centering
	\noindent
	\begin{subfigure}{0.5\textwidth}
		\centering
		\includegraphics[width=0.8\textwidth]{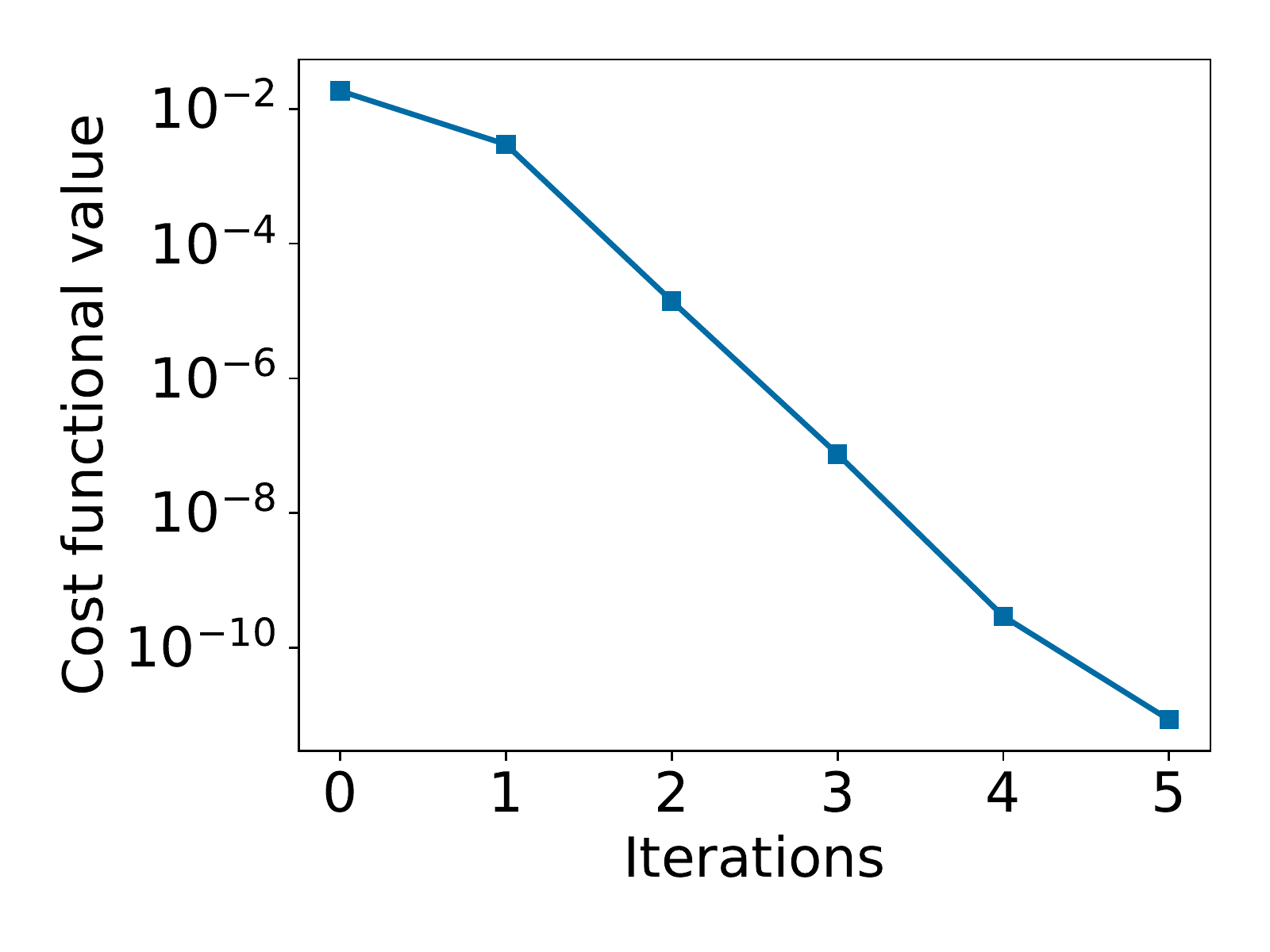}
		\caption{Evolution of the cost functional.}
	\end{subfigure}%
	\hfil%
	\begin{subfigure}{0.5\textwidth}
		\centering
		\includegraphics[width=0.8\textwidth]{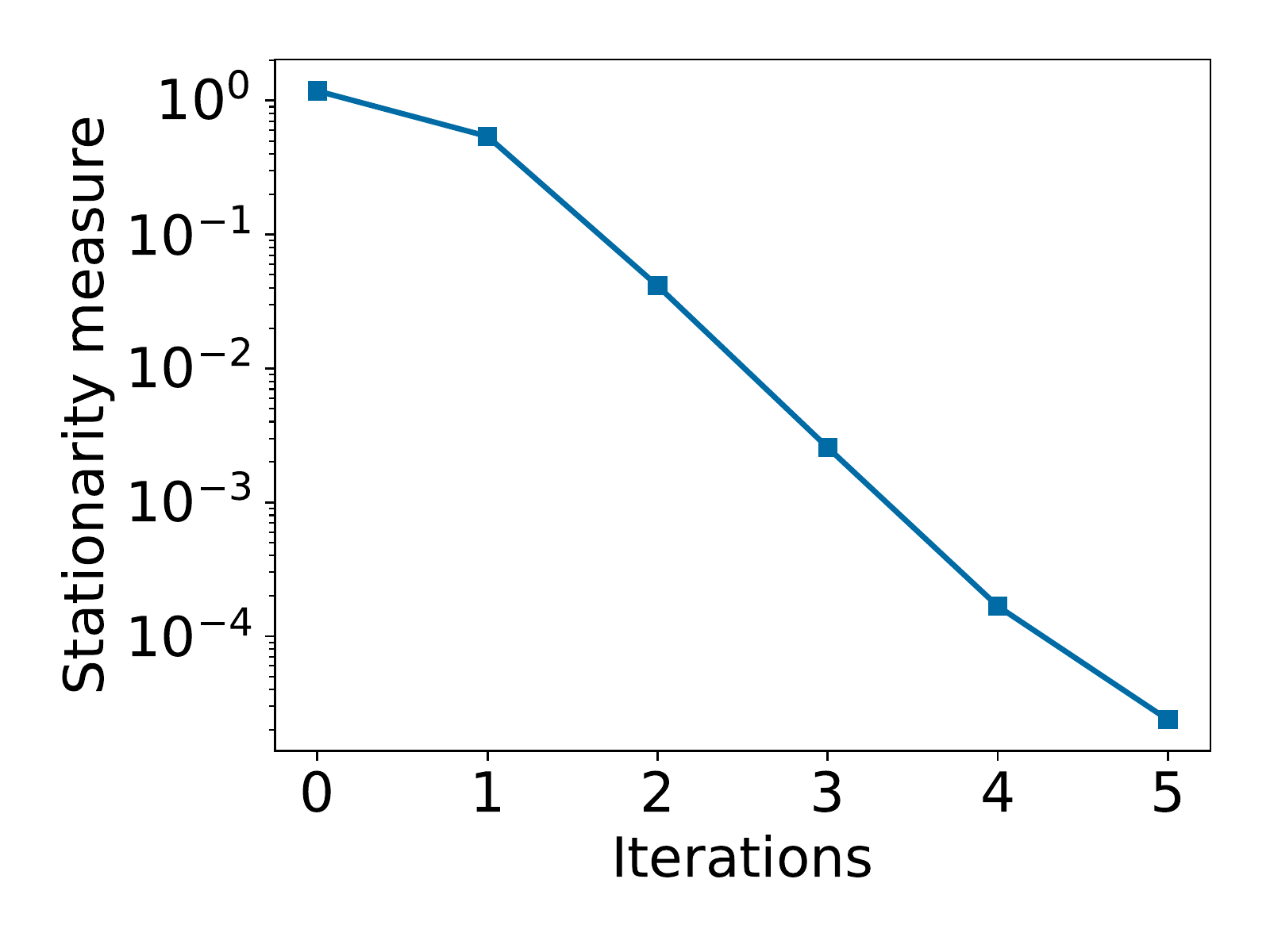}
		\caption{Evolution of the stationarity measure.}
	\end{subfigure}
	\caption{History of the ASM method for problem \cref{eq:flow_fine_model} with $\mathrm{Re} = \num{100}$.}
	\label{fig:asm_flow}
\end{figure}

As before, we solve problem \cref{eq:flow_fine_model} with the ASM method proposed in \cref{sec:space_mapping}. Let $(u\subfine(\Omega), p\subfine(\Omega))$ be the solution of the fine model \cref{eq:flow_fine} and let $(u\subcoarse(\Omega), p\subcoarse(\Omega))$ be the solution of the coarse model \cref{eq:flow_coarse}, respectively, on the domain $\Omega$. The fine model response is then given by $\fineresponse\colon \admissiblegeom \to \R^3; \Omega \mapsto q^\mathrm{f}\subout(\Omega)$, where $q^\mathrm{f}\subout(\Omega) = \begin{bsmallmatrix}
	q\subout^1(u\subfine(\Omega)), q\subout^2(u\subfine(\Omega)), q\subout^3(u\subfine(\Omega))
\end{bsmallmatrix}\transposed$, which is the vector of outlet flow rates corresponding to the fine model. The coarse model response $\coarseresponse\colon \admissiblegeom \to \R^3; \Omega \mapsto q^\mathrm{c}\subout(\Omega)$ is defined analogously, i.e., it is given by $q^\mathrm{c}\subout = \begin{bsmallmatrix}
	q\subout^1(u\subcoarse(\Omega)), q\subout^2(u\subcoarse(\Omega)), q\subout^3(u\subcoarse(\Omega))
\end{bsmallmatrix}\transposed$. We define the misalignment function $\misalignmentfunction\colon \R^3 \times \R^3 \to \R$ as
\begin{equation*}
	\misalignmentfunction(q_1, q_2) = \frac{1}{2} \norm{q_1 - q_2}{2}^2 = \frac{1}{2} \sum_{i=1}^{3} \left( q_1^{(i)} - q_2^{(i)} \right)^2.
\end{equation*}
Note that, in order to evaluate the space mapping function $\spacemappingfunction$, we need to solve the subproblem
\begin{equation*}
	\tilde{\vectorfield} = \argmin_{\vectorfield \in H^1(\hat{\Omega};\R^d)}\smallspace \misalignmentfunction\left( \coarseresponse\left( \numretraction_{\hat{\Omega}}(\vectorfield) \right), \fineresponse\left( \Omega\subfine \right) \right).
\end{equation*}
Due to the definition of $\misalignmentfunction$, this is equivalent to the problem
\begin{equation}
	\label{eq:parameter_extraction_flow}
	\tilde{\Omega} = \argmin_{\Omega\in\admissiblegeom}\smallspace \frac{1}{2} \sum_{i=1}^{3} \left( q\subout^i(u\subcoarse(\Omega)) - q\subout^i(u\subfine(\Omega\subfine)) \right)^2,
\end{equation}
in the sense that we have $\tilde{\Omega} = \numretraction_{\hat{\Omega}}(\tilde{\vectorfield})$. Again, the optimization problem \cref{eq:parameter_extraction_flow} is similar to the coarse model optimization problem \cref{eq:flow_coarse_model}, where $q\subdes$ is replaced by $q\subout^i(u\subfine(\Omega\subfine))$ on $\Gamma\subout^i$.

\begin{figure}[!b]
	\centering
	\begin{subfigure}{\textwidth}
		\centering
		\includegraphics[width=0.66\textwidth, trim=0cm 0cm 0cm 1.5cm, clip]{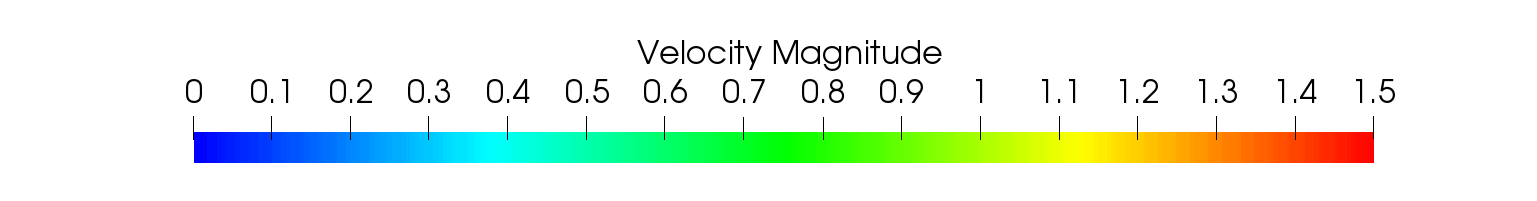}
	\end{subfigure}

	\begin{subfigure}[!t]{0.5\textwidth}
		\centering
		\begin{minipage}[!t]{0.6\textwidth}
			\includegraphics[width=\textwidth]{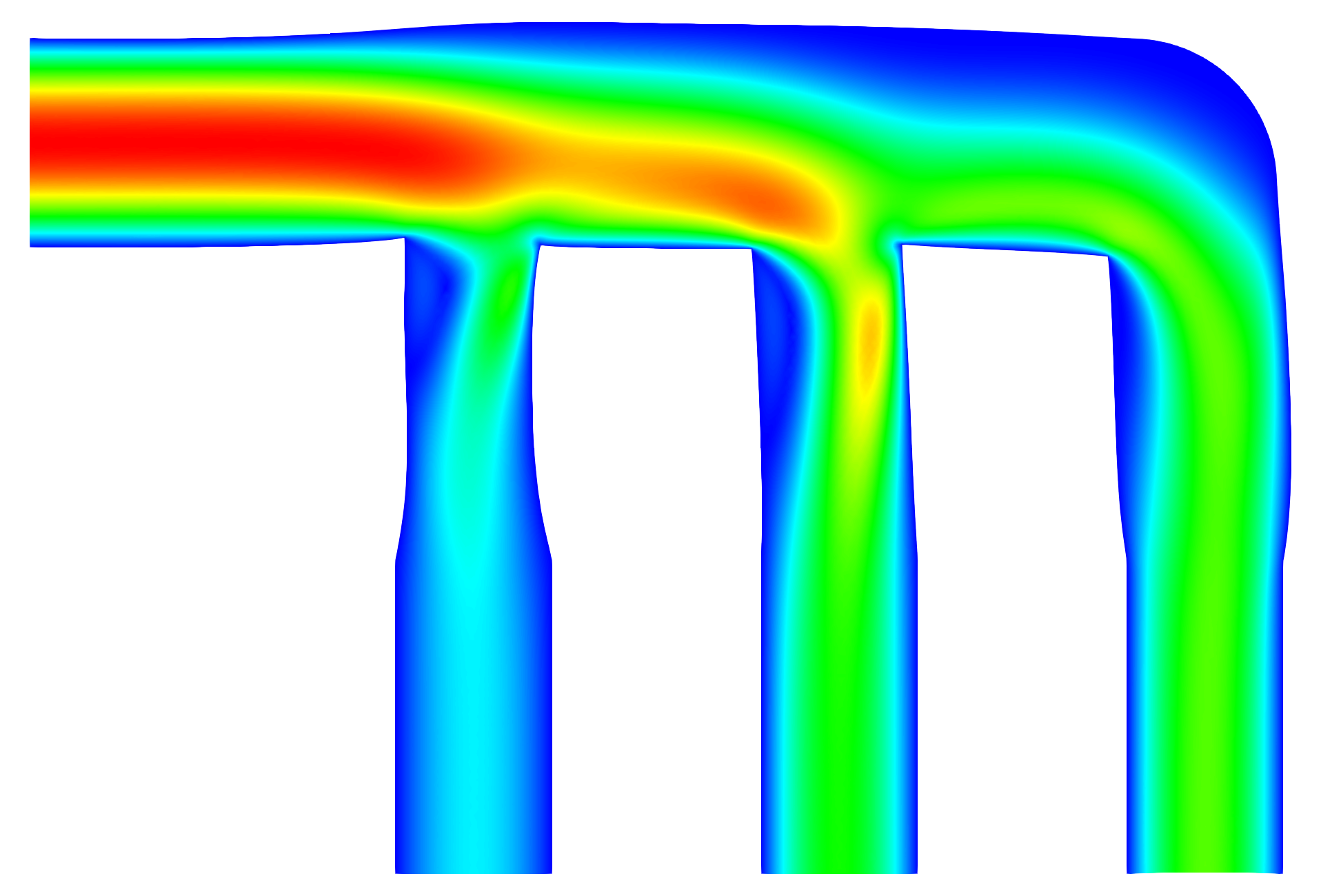}
		\end{minipage}%
		\hfil%
		\begin{minipage}[!t]{0.4\textwidth}
			\includegraphics[width=\textwidth]{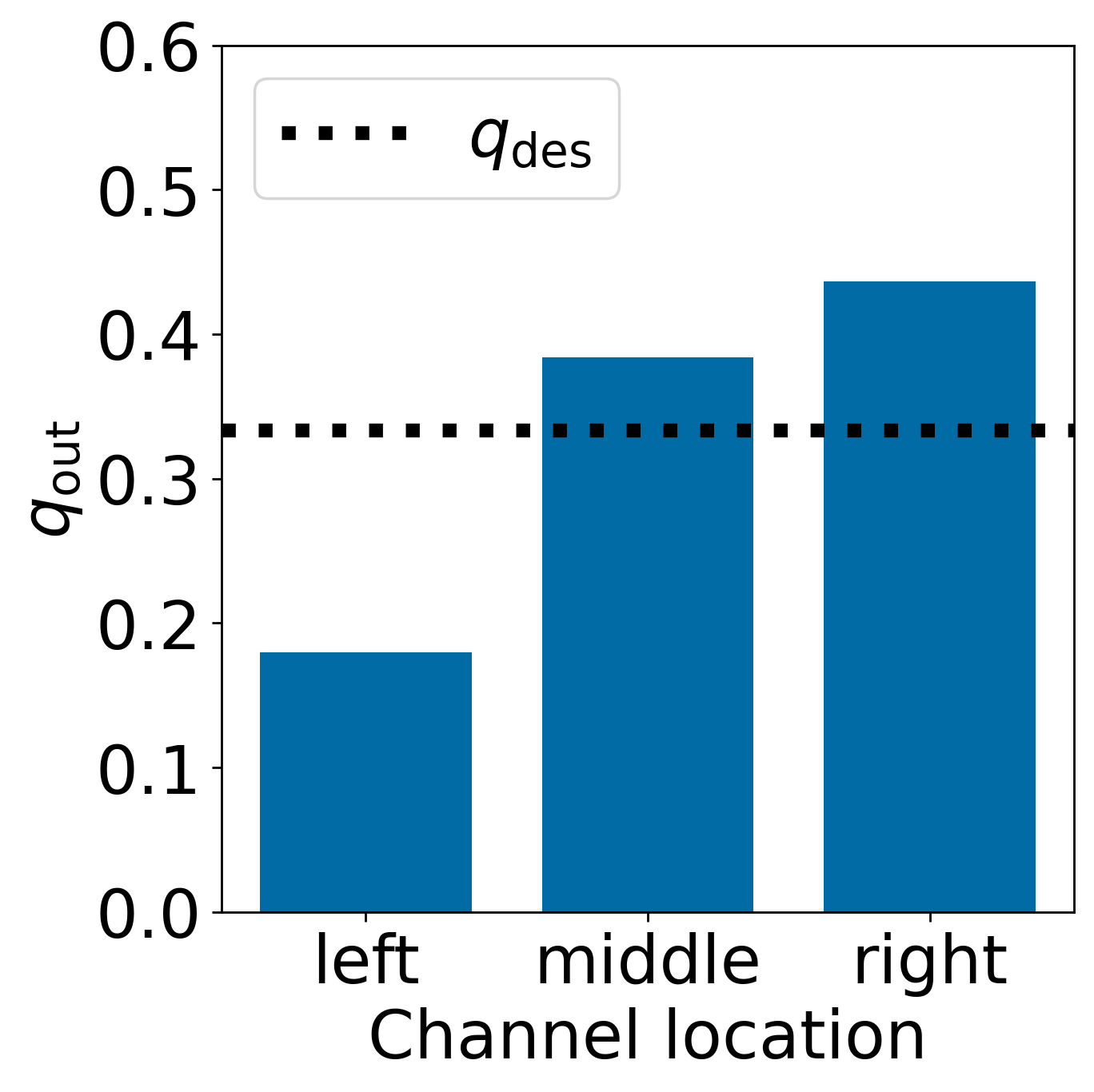}
		\end{minipage}
		\caption{Iteration 0.}
		\label{fig:flow_1}
	\end{subfigure}%
	\hfil%
	\begin{subfigure}[!t]{0.5\textwidth}
		\centering
		\begin{minipage}[!t]{0.6\textwidth}
			\includegraphics[width=\textwidth]{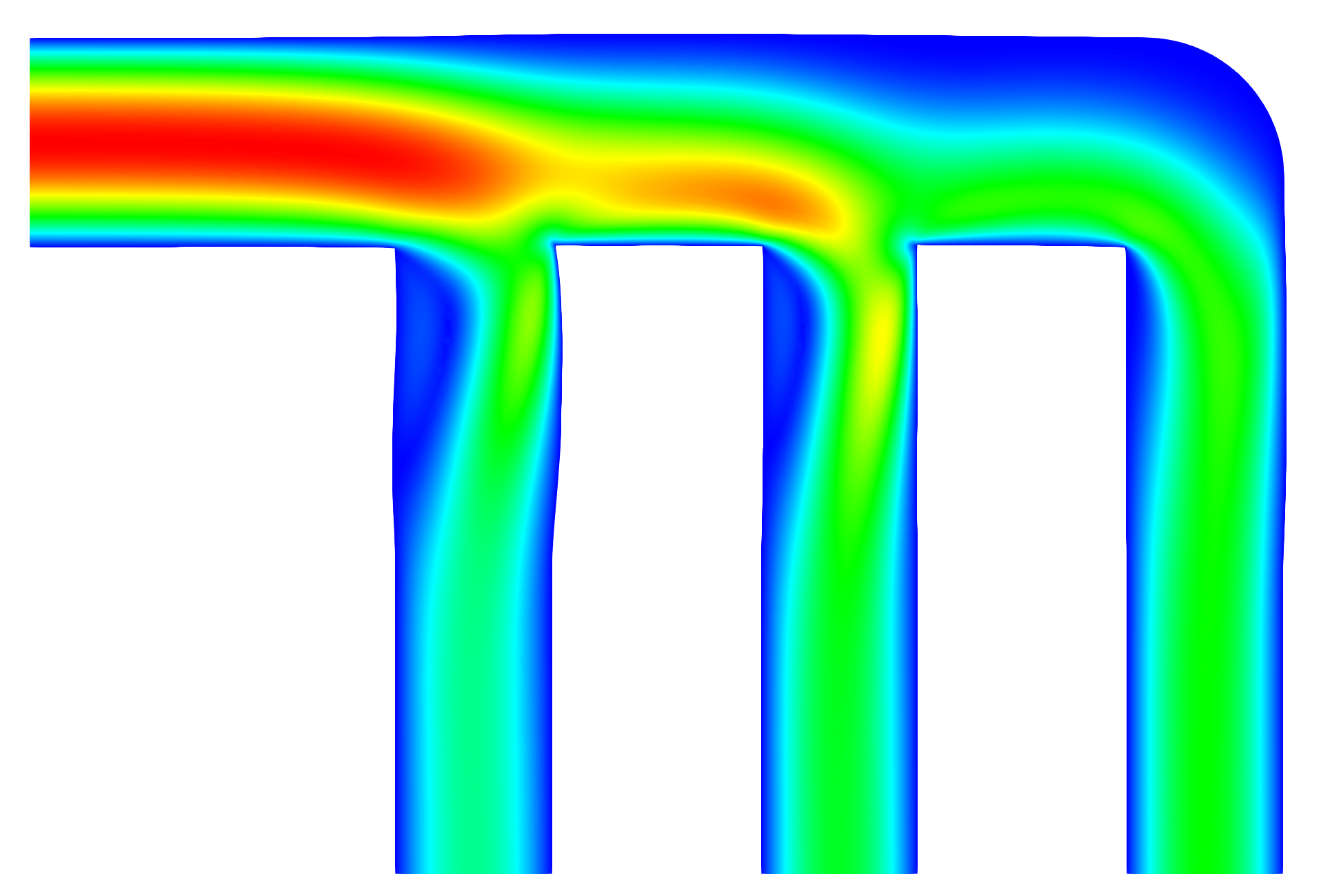}
		\end{minipage}%
		\hfil%
		\begin{minipage}[!t]{0.4\textwidth}
			\includegraphics[width=\textwidth]{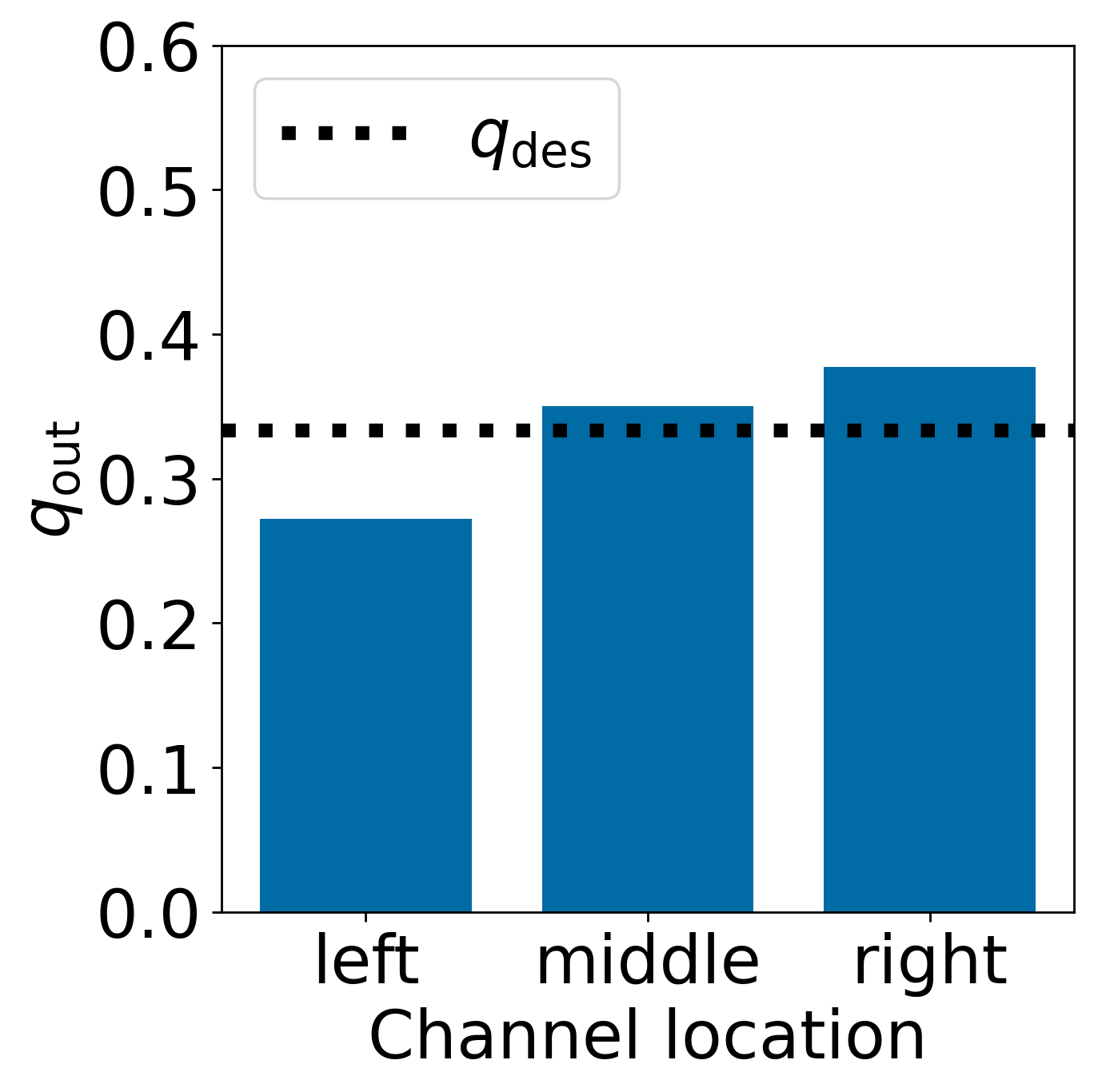}
		\end{minipage}
		\caption{Iteration 1.}
	\end{subfigure}

	\begin{subfigure}[!t]{0.49\textwidth}
		\centering
		\begin{minipage}[!t]{0.6\textwidth}
			\includegraphics[width=\textwidth]{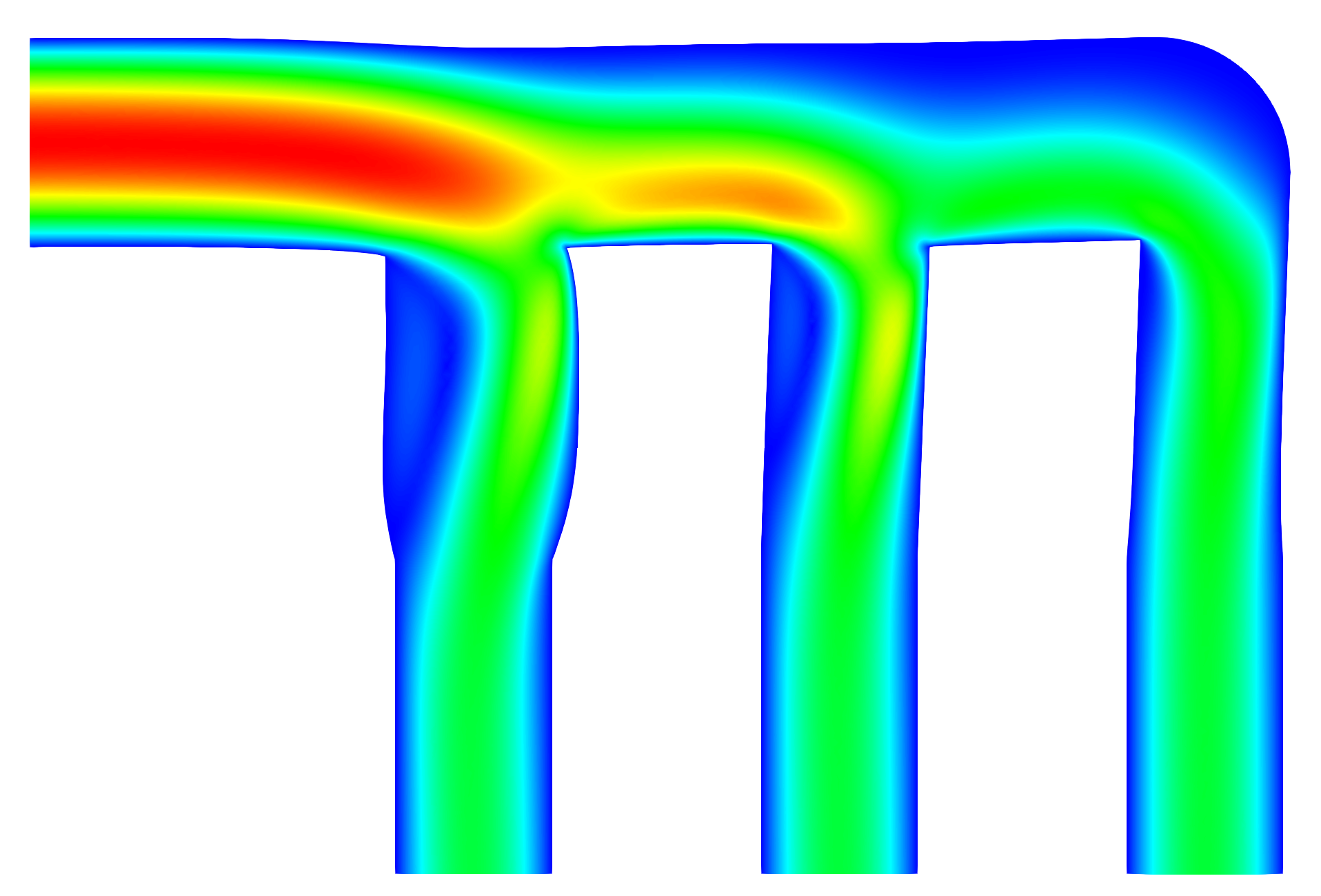}
		\end{minipage}%
		\hfil%
		\begin{minipage}[!t]{0.4\textwidth}
			\includegraphics[width=\textwidth]{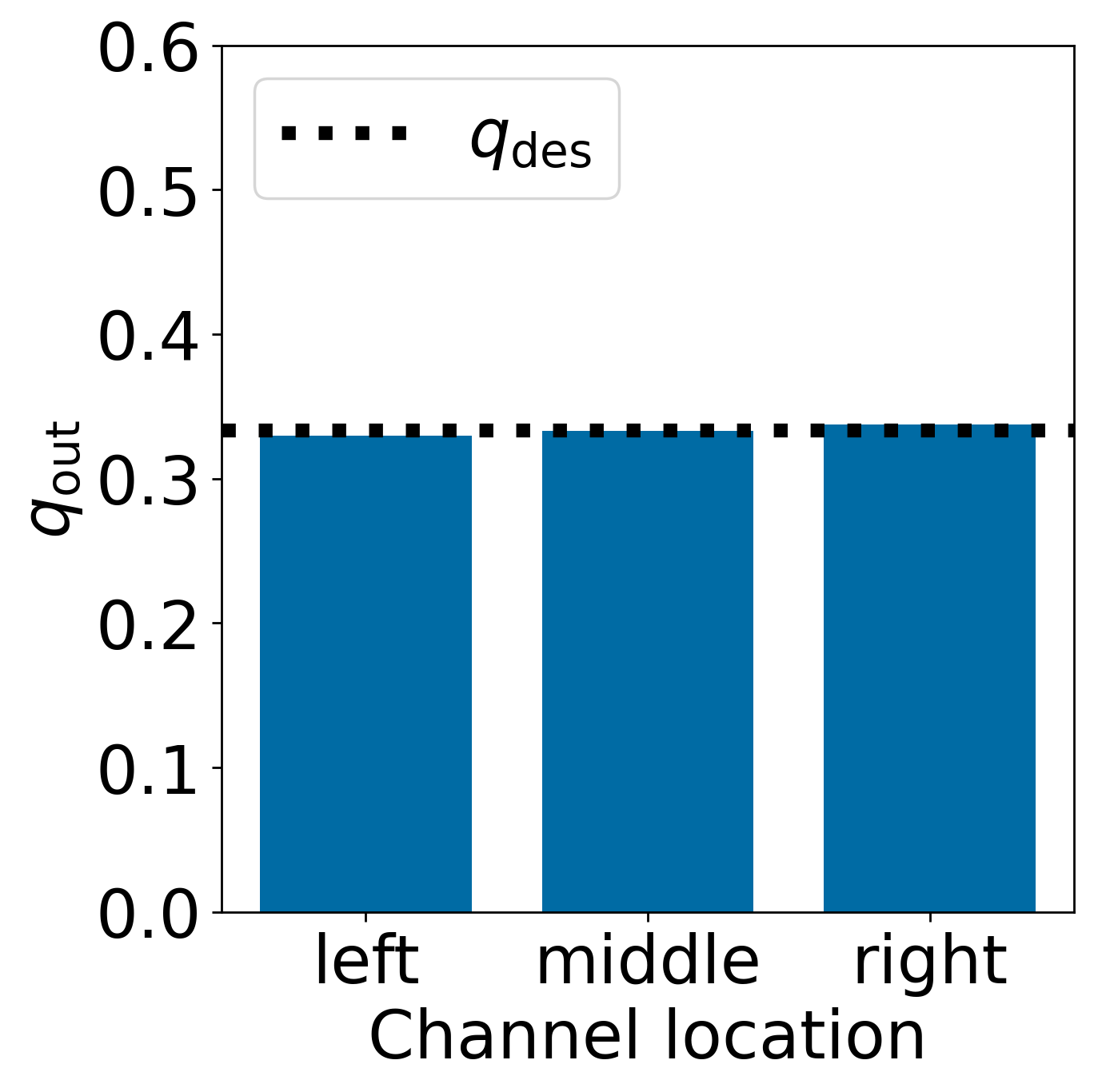}
		\end{minipage}
		\caption{Iteration 2.}
	\end{subfigure}%
	\hfil%
	\begin{subfigure}[!t]{0.49\textwidth}
		\centering
		\begin{minipage}[!t]{0.6\textwidth}
			\includegraphics[width=\textwidth]{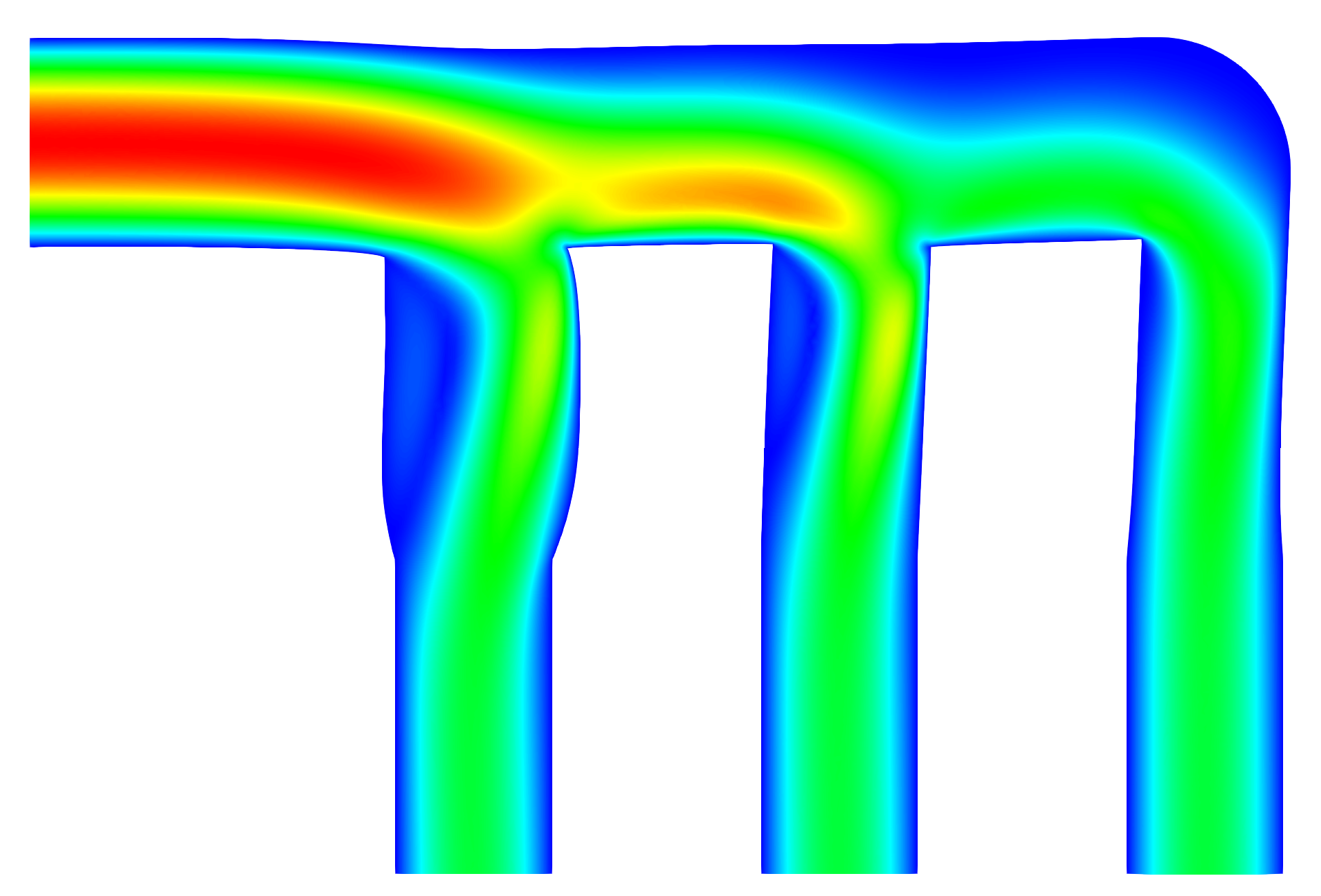}
		\end{minipage}%
		\hfil%
		\begin{minipage}[!t]{0.4\textwidth}
			\includegraphics[width=\textwidth]{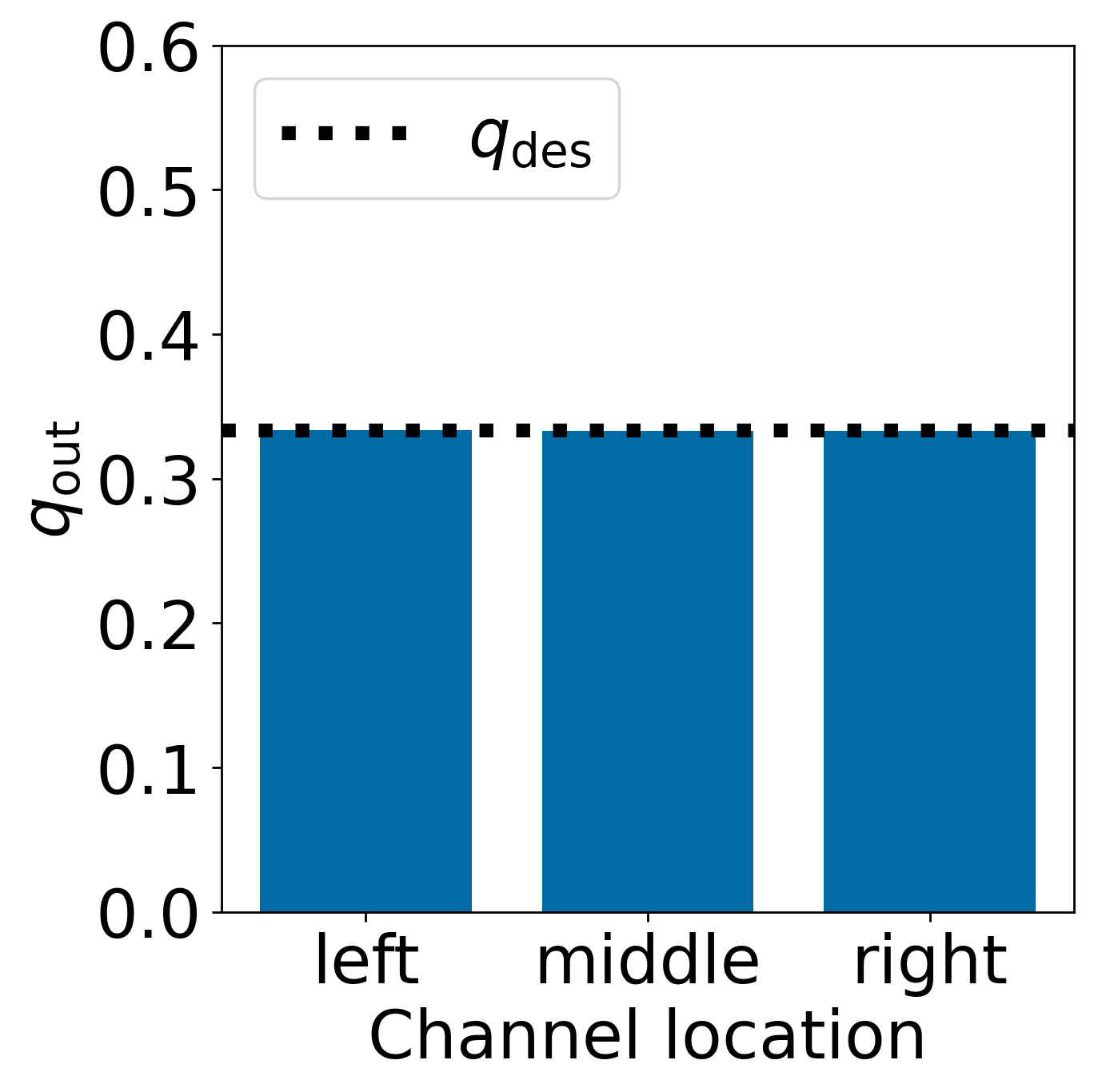}
		\end{minipage}
		\caption{Iteration 3.}
	\end{subfigure}

	\begin{subfigure}[!t]{0.49\textwidth}
		\centering
		\begin{minipage}[!t]{0.6\textwidth}
			\includegraphics[width=\textwidth]{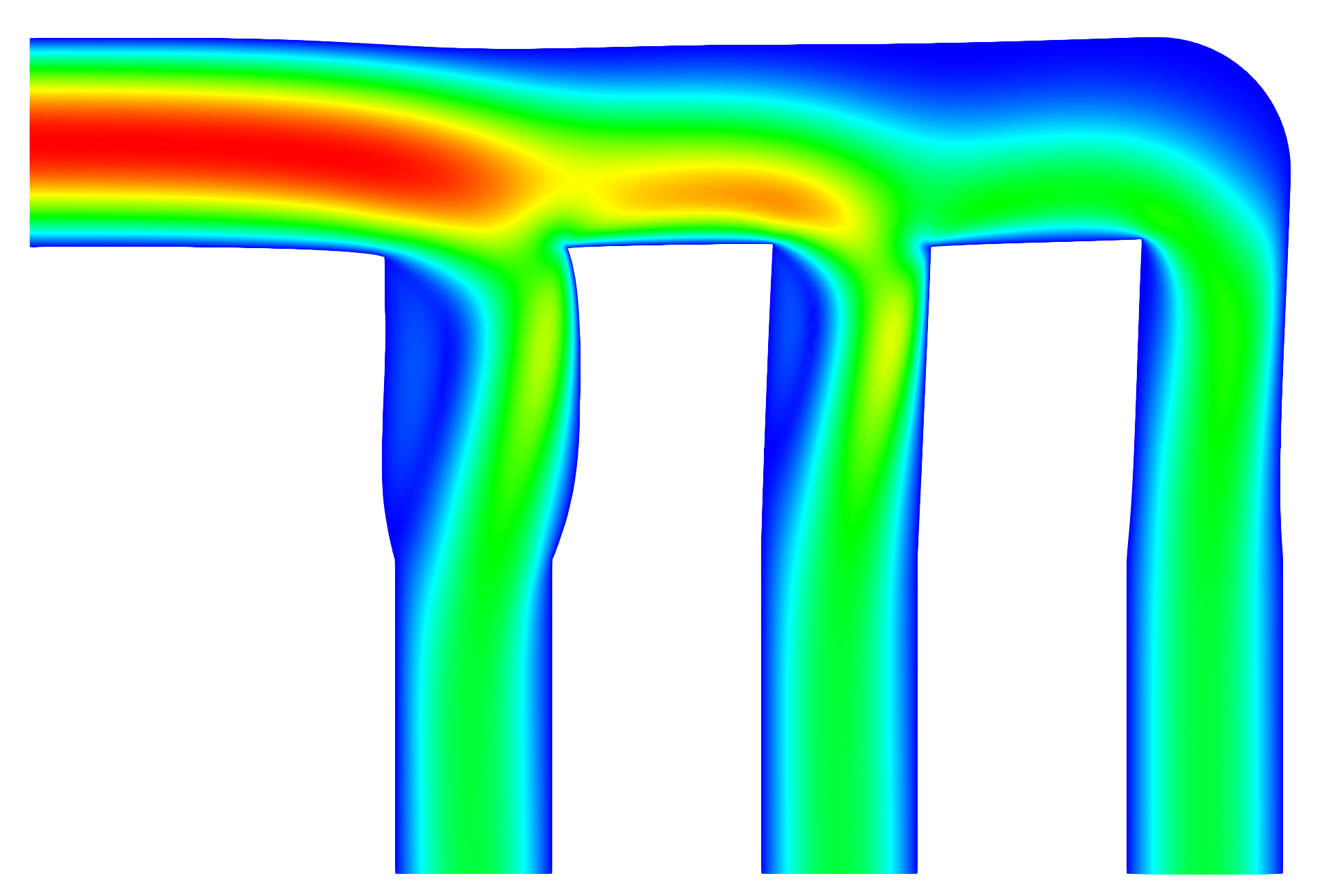}
		\end{minipage}%
		\hfil%
		\begin{minipage}[!t]{0.4\textwidth}
			\includegraphics[width=\textwidth]{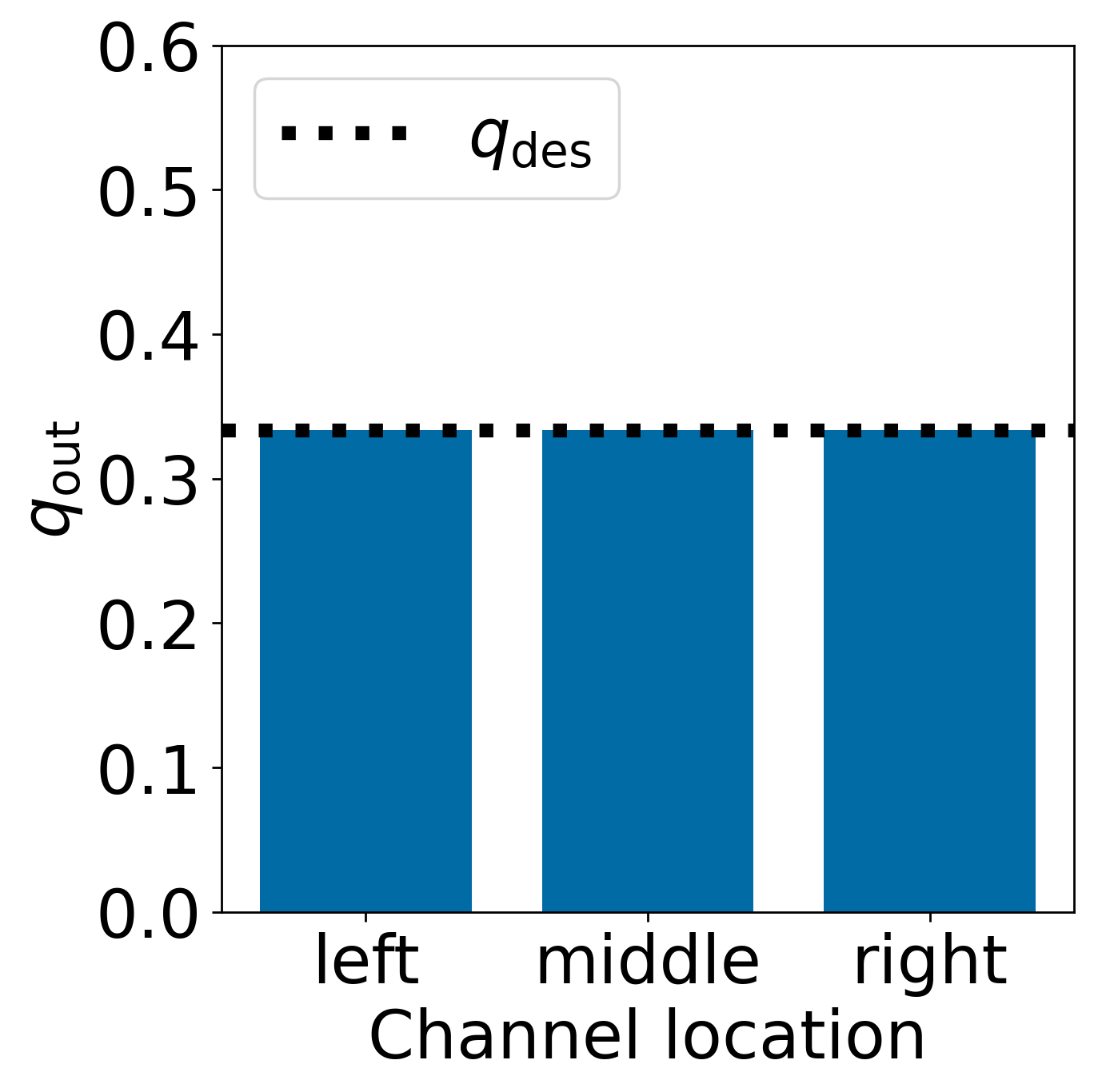}
		\end{minipage}
		\caption{Iteration 4.}
	\end{subfigure}
	\hfil
	\begin{subfigure}[!t]{0.49\textwidth}
		\centering
		\begin{minipage}[!t]{0.6\textwidth}
			\includegraphics[width=\textwidth]{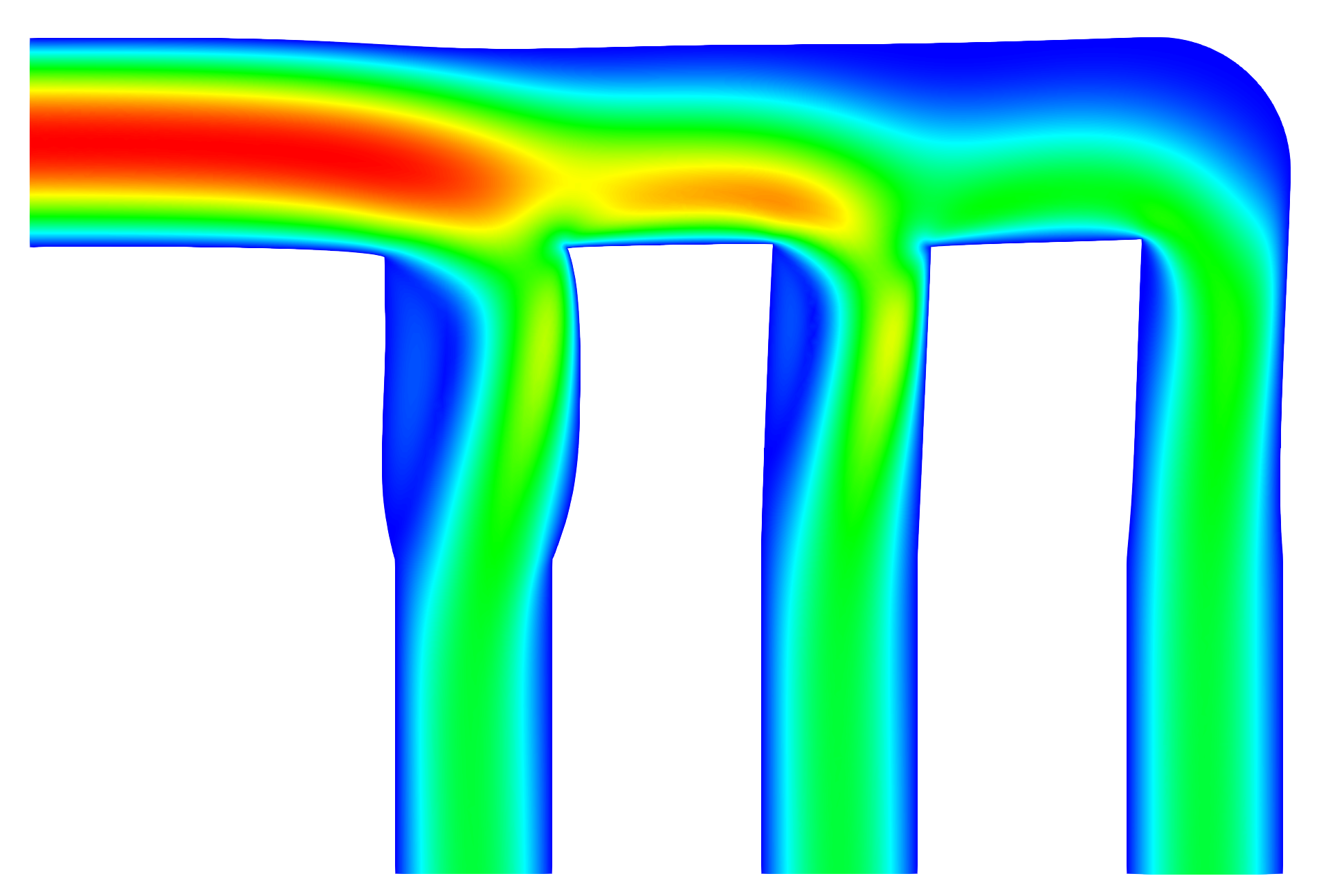}
		\end{minipage}%
		\hfil%
		\begin{minipage}[!t]{0.4\textwidth}
			\includegraphics[width=\textwidth]{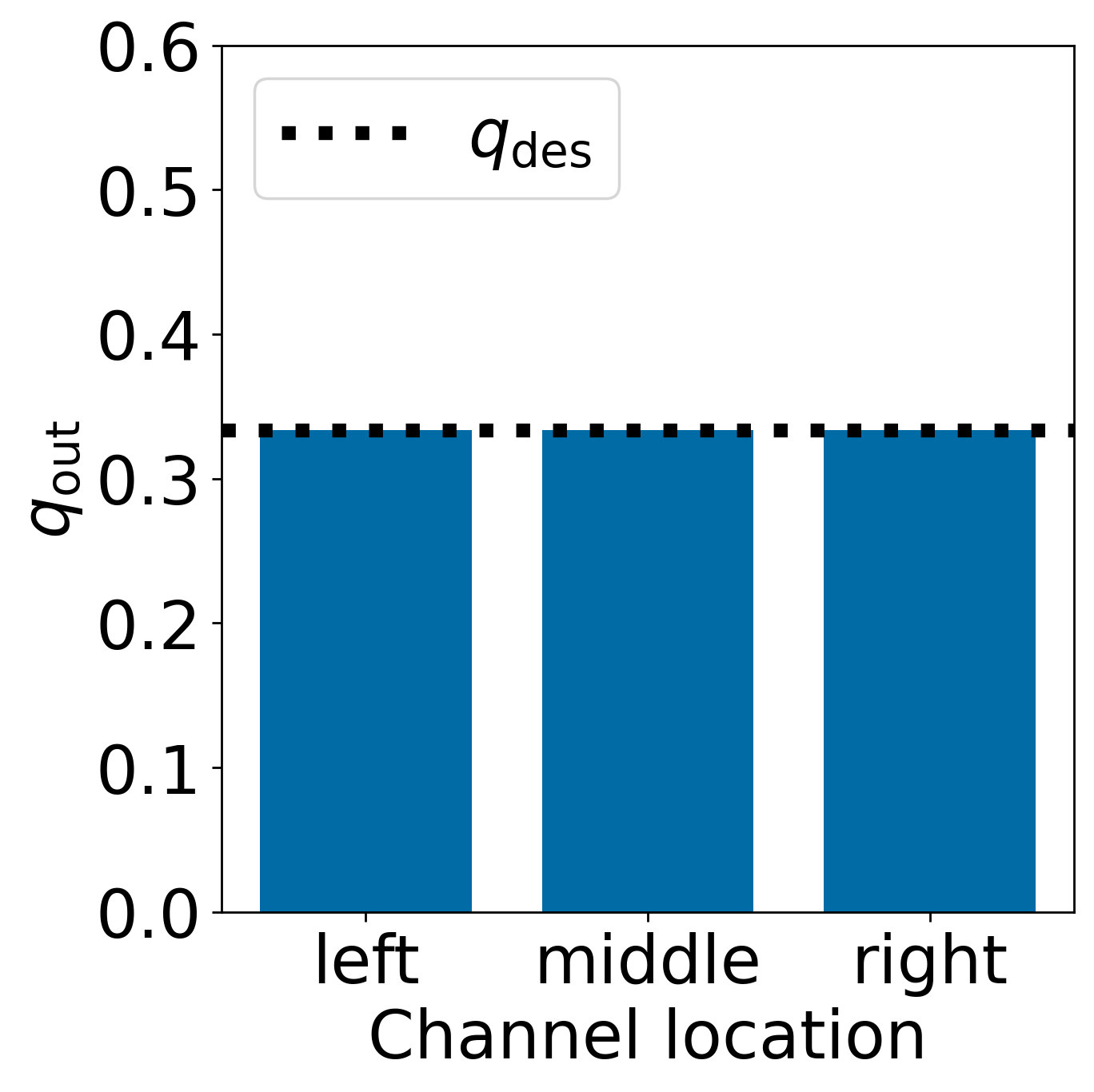}
		\end{minipage}
		\caption{Iteration 5.}
	\end{subfigure}
	
	\caption{Evolution of the geometry $\Omega$, the velocity field $u$, and the outlet flow rates $q\subout^i$ over the course of the ASM method for problem \cref{eq:flow_fine_model} with $\mathrm{Re} = \num{100}$.}
	\label{fig:flow_geometries}
\end{figure}

For our numerical solution of \cref{eq:flow_fine_model} with the ASM method, we again employ our optimization software cashocs \cite{Blauth2021cashocs}. The reference geometry for our case is shown in \cref{fig:schema_flow} and we discretize this geometry with a nonuniform mesh consisting of \num{7196} nodes and \num{12908} triangles, which has a finer resolution near the wall boundary. Again, we use the reference geometry as initial guess for all considered optimization problems. We solve the state and adjoint systems with FEniCS and discretize the velocity with piecewise quadratic Lagrange elements and the pressure with piecewise linear Lagrange elements, which are LBB-stable for \cref{eq:flow_fine,eq:flow_coarse}. For the fine model, we use a finer discretization consisting of \num{24983}~nodes and \num{46995}~triangles, which additionally increases the numerical cost for its solution. The coarse model shape optimization problems \cref{eq:flow_coarse_model,eq:parameter_extraction_flow} are solved with a limited memory BFGS method (cf.~\cite{Schulz2016Efficient}) implemented in cashocs, using a relative tolerance of \num{1e-2}. For the shape optimization, we keep the inlet $\Gamma\subin$ as well as the lower half of the outlet pipes fixed, such that only a part of $\Gamma_\mathrm{wall}$ can be deformed (cf.~\cref{fig:flow_geometries}). For the Reynolds number, we use $\mathrm{Re} = 100$ for the fine model. We choose the inlet velocity profile as $u\subin(x) = -6 x_2 (x_2 + 1)$, where the inlet is given as $\Gamma_\mathrm{in} = \Set{x \in \R^2 | x_1 = 0, x_2 \in [-1, 0]}$, which corresponds to a flow rate of $q\subin = 1$. Finally, we terminate the ASM method if the stationarity measure $\sigma$ (cf.~\cref{eq:stationarity_measure}) is smaller than the relative tolerance of $\tau = \num{1e-4}$.

The history of the ASM method is depicted in \cref{fig:asm_flow}, where the evolution of the cost functional as well as the stationarity measure is shown. Again, we observe that the ASM method is very efficient at solving the fine model optimization problem \cref{eq:flow_fine_model} as it requires only five iterations to do so. We observe that the cost functional decreases by over nine orders of magnitude and that the stationarity measure $\sigma$ decreases by about five orders of magnitude. This indicates that we have indeed found a solution to problem \cref{eq:flow_fine_model}.

In \cref{fig:flow_geometries}, we depict the evolution of the geometry over the course of the ASM method. There, the geometry $\Omega$, the resulting velocity field $u\subfine(\Omega)$, and the outlet flow rates $q\subout^i$ are shown for each iteration of the method. We observe that the ASM method is very efficient at solving \cref{eq:flow_fine_model}: We start from a quite suboptimal initial guess, given by the optimal geometry for the coarse model problem \cref{eq:flow_coarse_model}, which exhibits a bad flow distribution where the left pipe receives far too little fluid and the middle and the right one receive far too much fluid (cf.~\cref{fig:flow_1}). This is then successfully corrected by the space mapping method. For the first iterate, we can still see a mismatch in the flow distribution, but for all later ones this is not visible anymore, and we achieve a perfect flow uniformity on the numerical level in the final iteration. 
The fact that our ASM method only required five iterations for solving problem \cref{eq:flow_fine_model} shows that it is very efficient and rapidly convergent.

\subsection{Coupling with Commercial Solvers}
\label{ssec:fluent}

\begin{figure}[!b]
	\centering
	\noindent
	\begin{subfigure}{0.5\textwidth}
		\centering
		\includegraphics[width=0.8\textwidth]{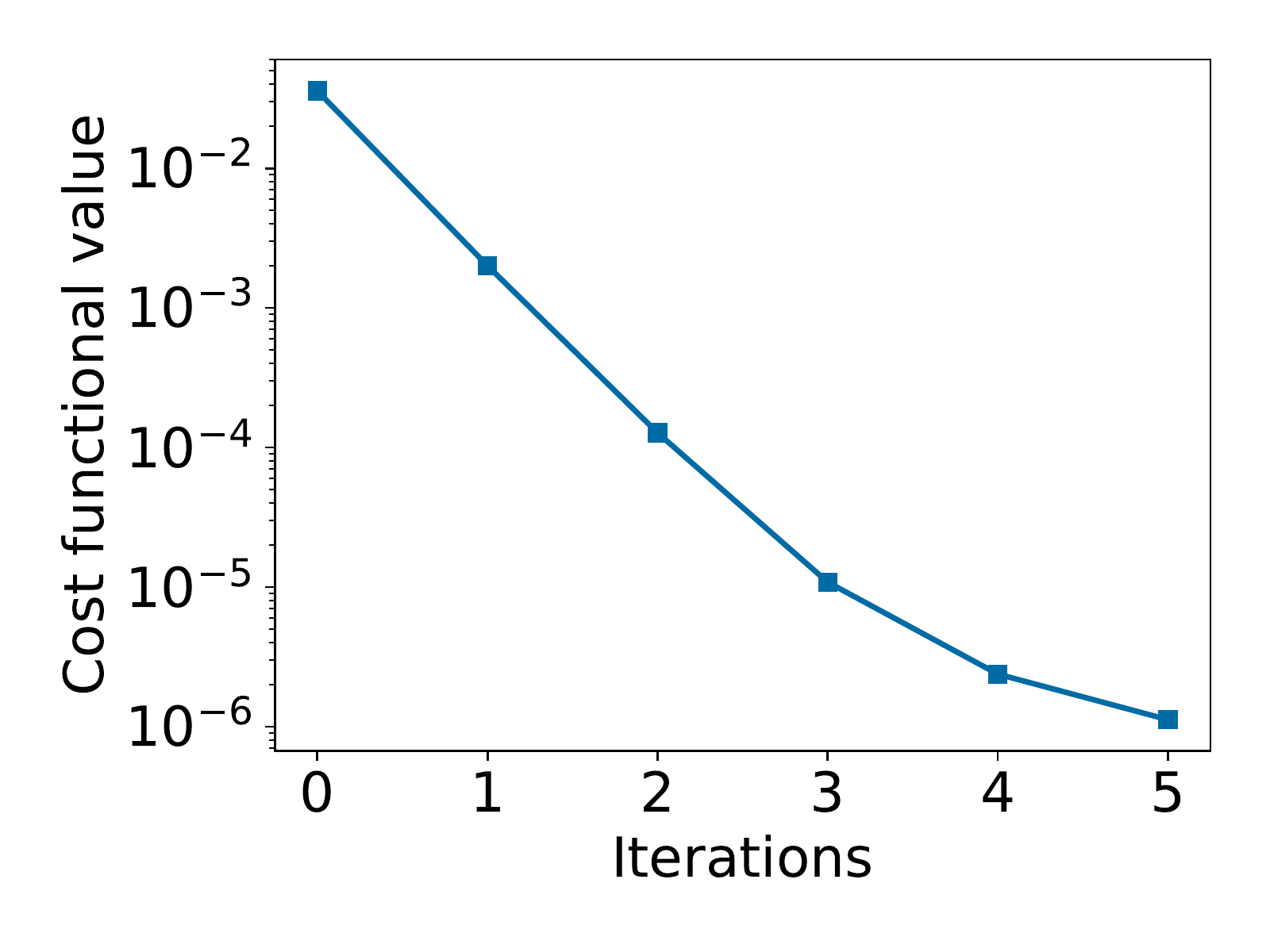}
		\caption{Evolution of the cost functional.}
	\end{subfigure}%
	\hfil%
	\begin{subfigure}{0.5\textwidth}
		\centering
		\includegraphics[width=0.8\textwidth]{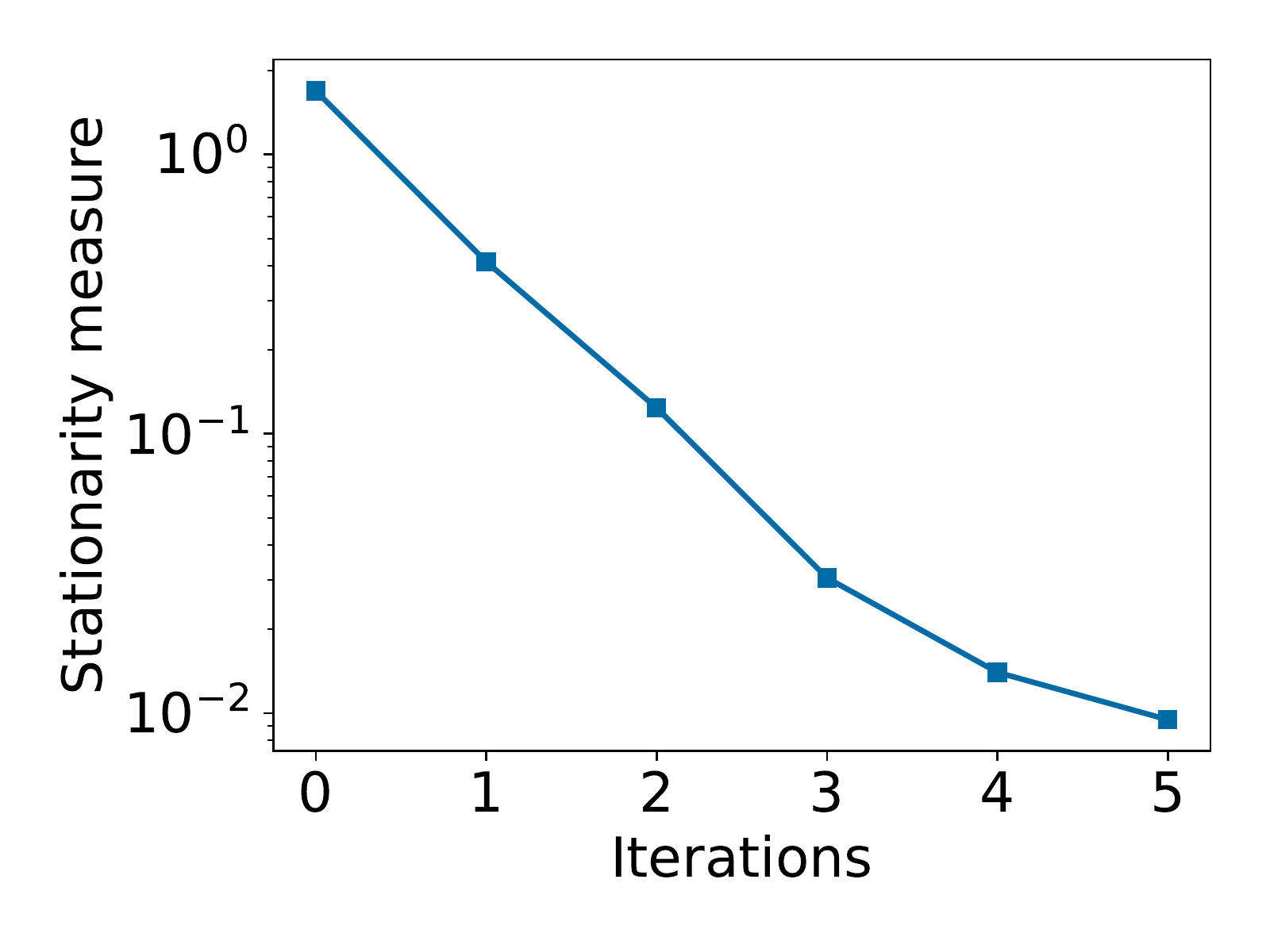}
		\caption{Evolution of the stationarity measure.}
	\end{subfigure}
	\caption{History of the ASM method for problem \cref{eq:flow_fine_model} with $\mathrm{Re} = \num{1000}$.}
	\label{fig:asm_fluent}
\end{figure}

\begin{figure}[!t]
	\centering
	\begin{subfigure}{\textwidth}
		\centering
		\begin{tikzpicture}
		\node at (0,0) {\includegraphics[width=0.55\textwidth, trim=0cm 0cm 0cm 46cm, clip]{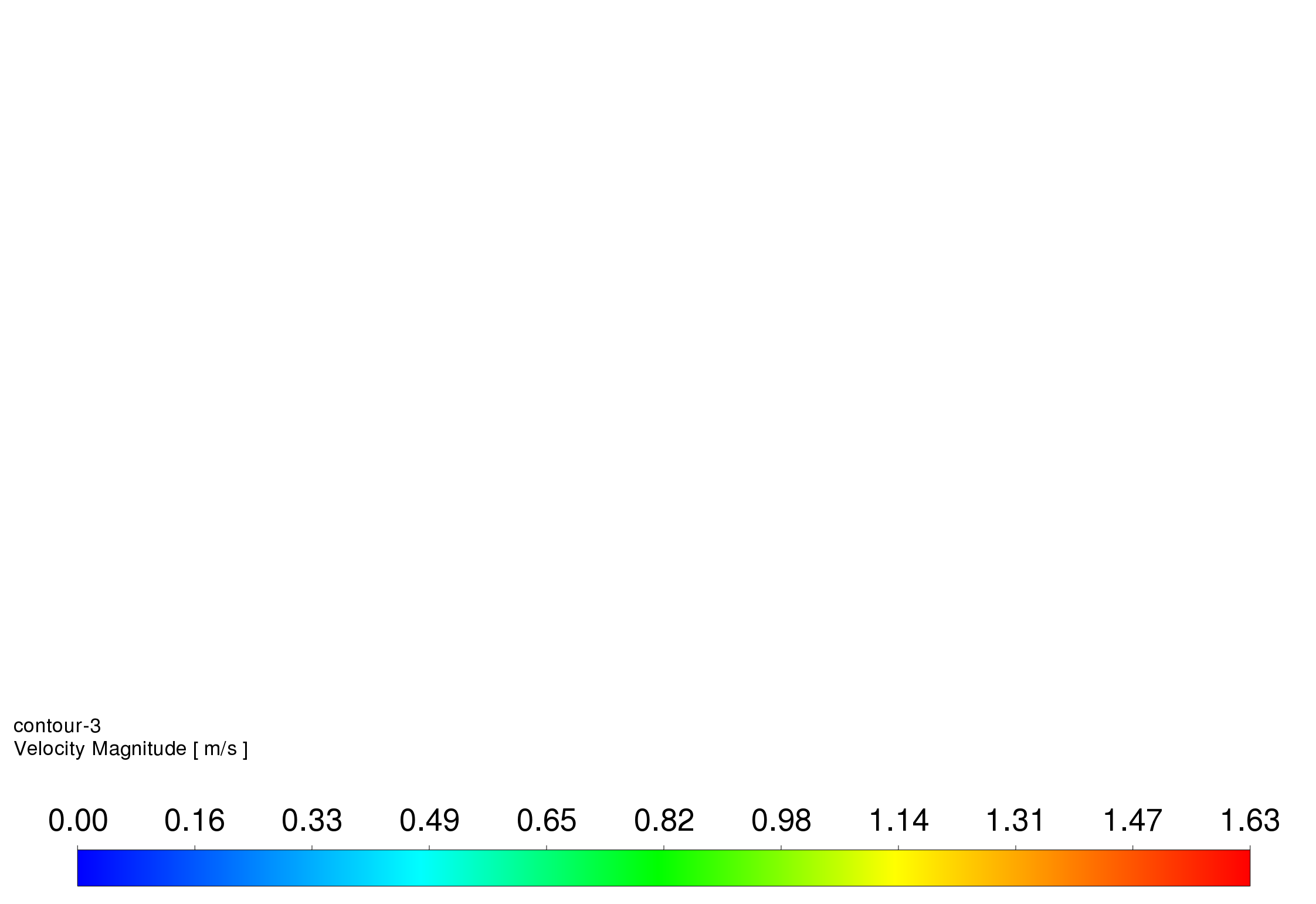}};
		
		\node at (0, 0.5) {\scriptsize Velocity Magnitude};
		\end{tikzpicture}
	\end{subfigure}

	\begin{subfigure}[!t]{0.5\textwidth}
		\centering
		\begin{minipage}[!t]{0.6\textwidth}
			\includegraphics[width=\textwidth]{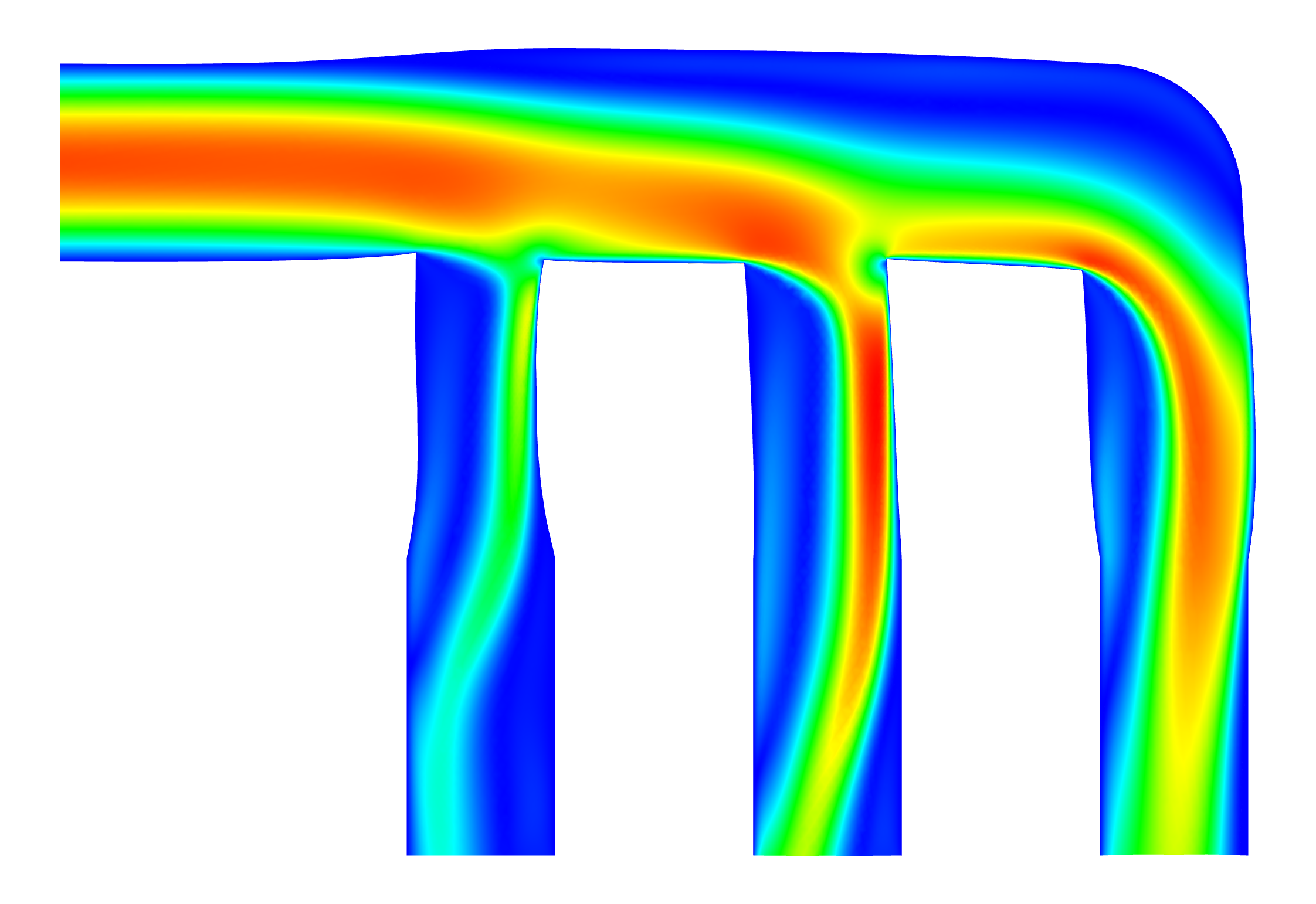}
		\end{minipage}%
		\hfil%
		\begin{minipage}[!t]{0.4\textwidth}
			\includegraphics[width=\textwidth]{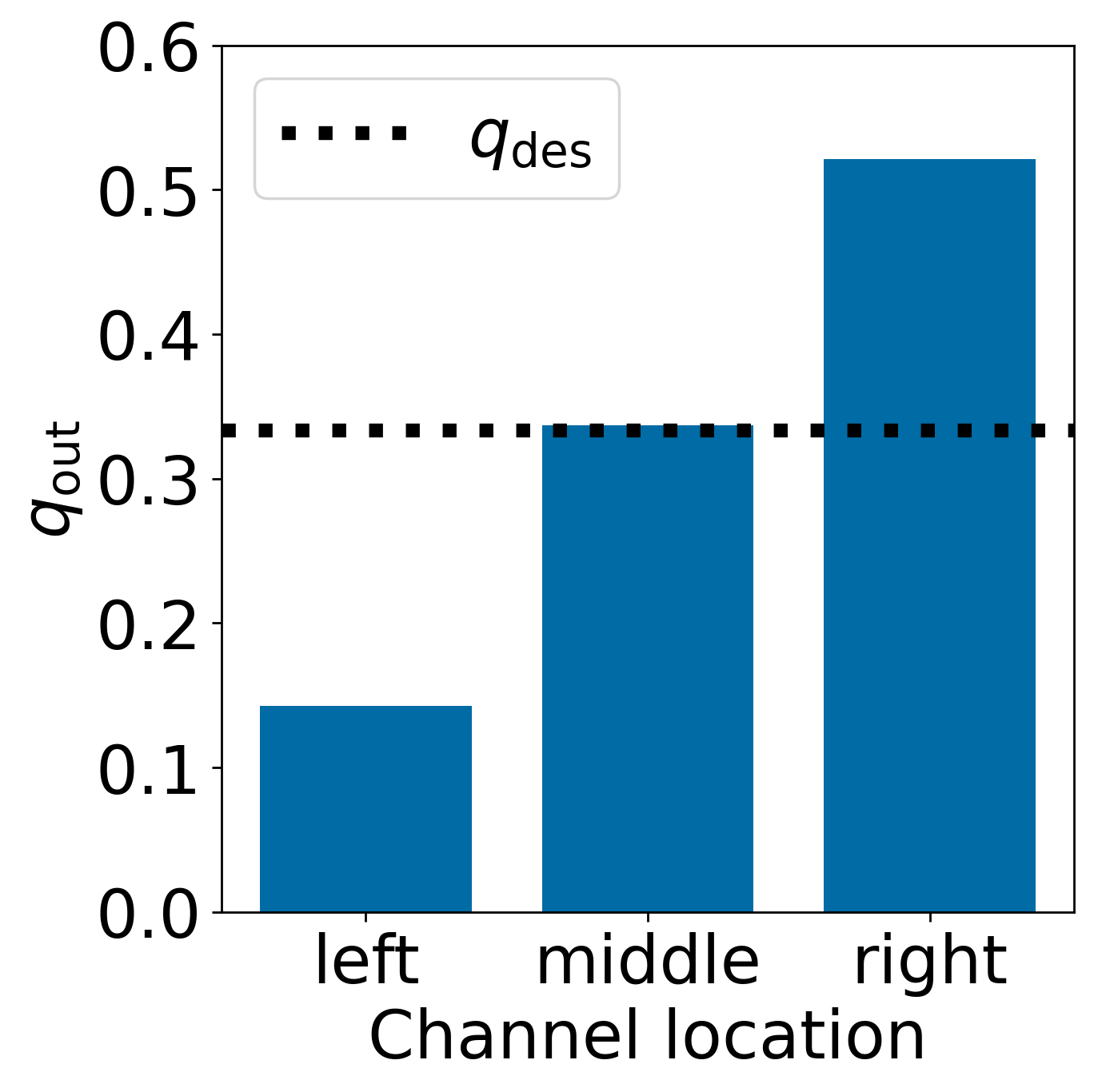}
		\end{minipage}
		\caption{Iteration 0.}
	\end{subfigure}%
	\hfil%
	\begin{subfigure}[!t]{0.5\textwidth}
		\centering
		\begin{minipage}[!t]{0.6\textwidth}
			\includegraphics[width=\textwidth]{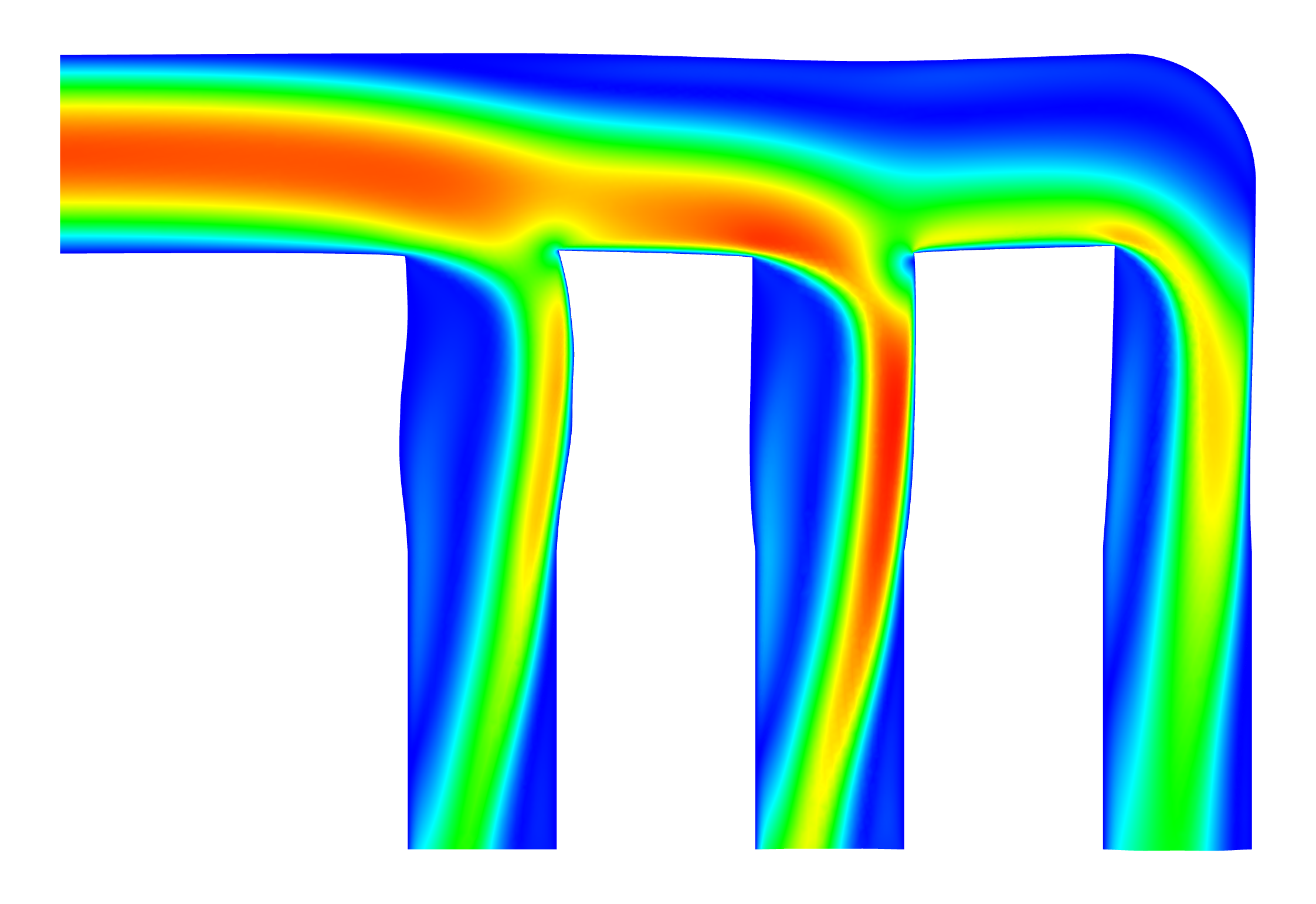}
		\end{minipage}%
		\hfil%
		\begin{minipage}[!t]{0.4\textwidth}
			\includegraphics[width=\textwidth]{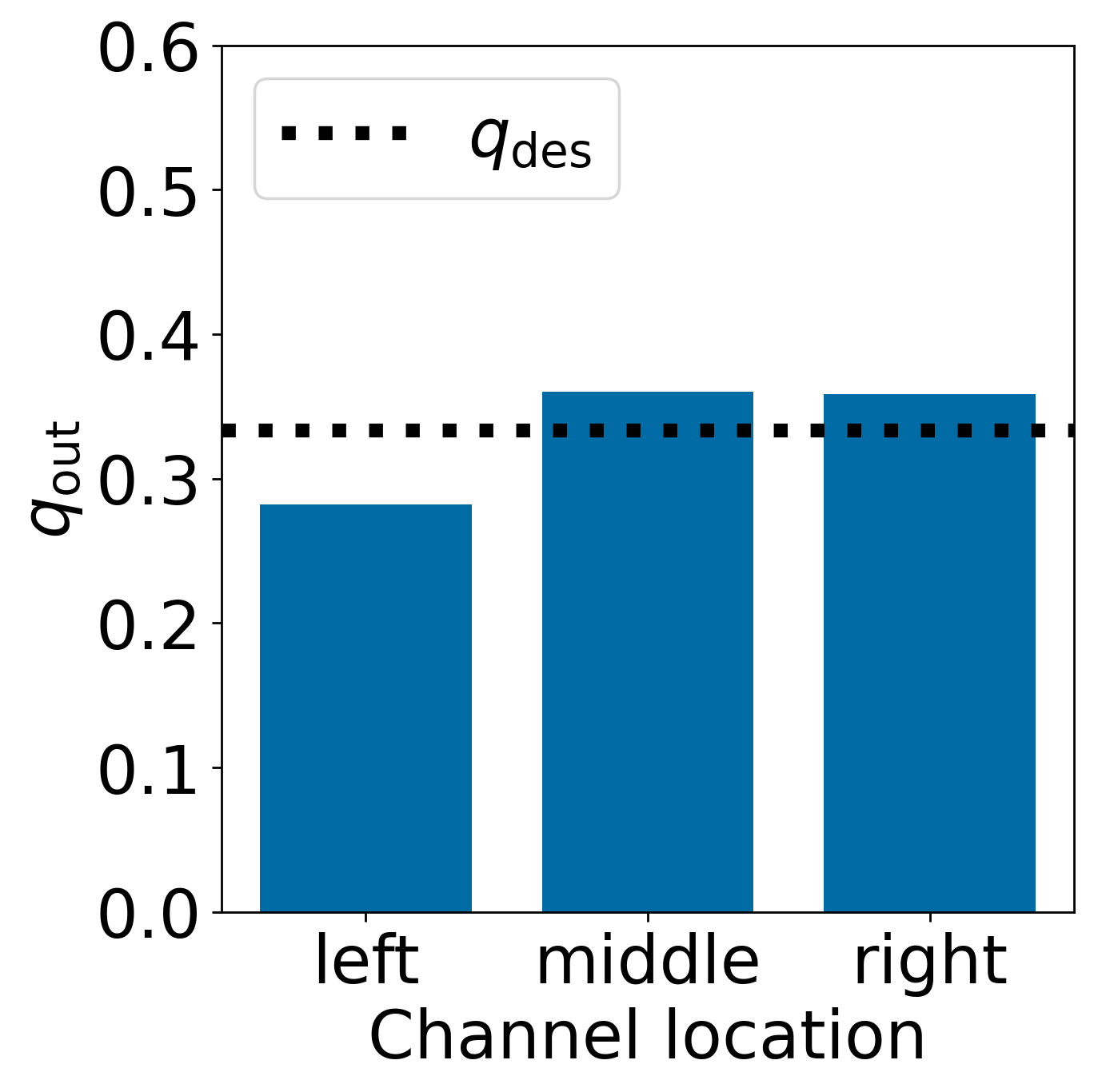}
		\end{minipage}
		\caption{Iteration 1.}
	\end{subfigure}

	\begin{subfigure}[!t]{0.49\textwidth}
		\centering
		\begin{minipage}[!t]{0.6\textwidth}
			\includegraphics[width=\textwidth]{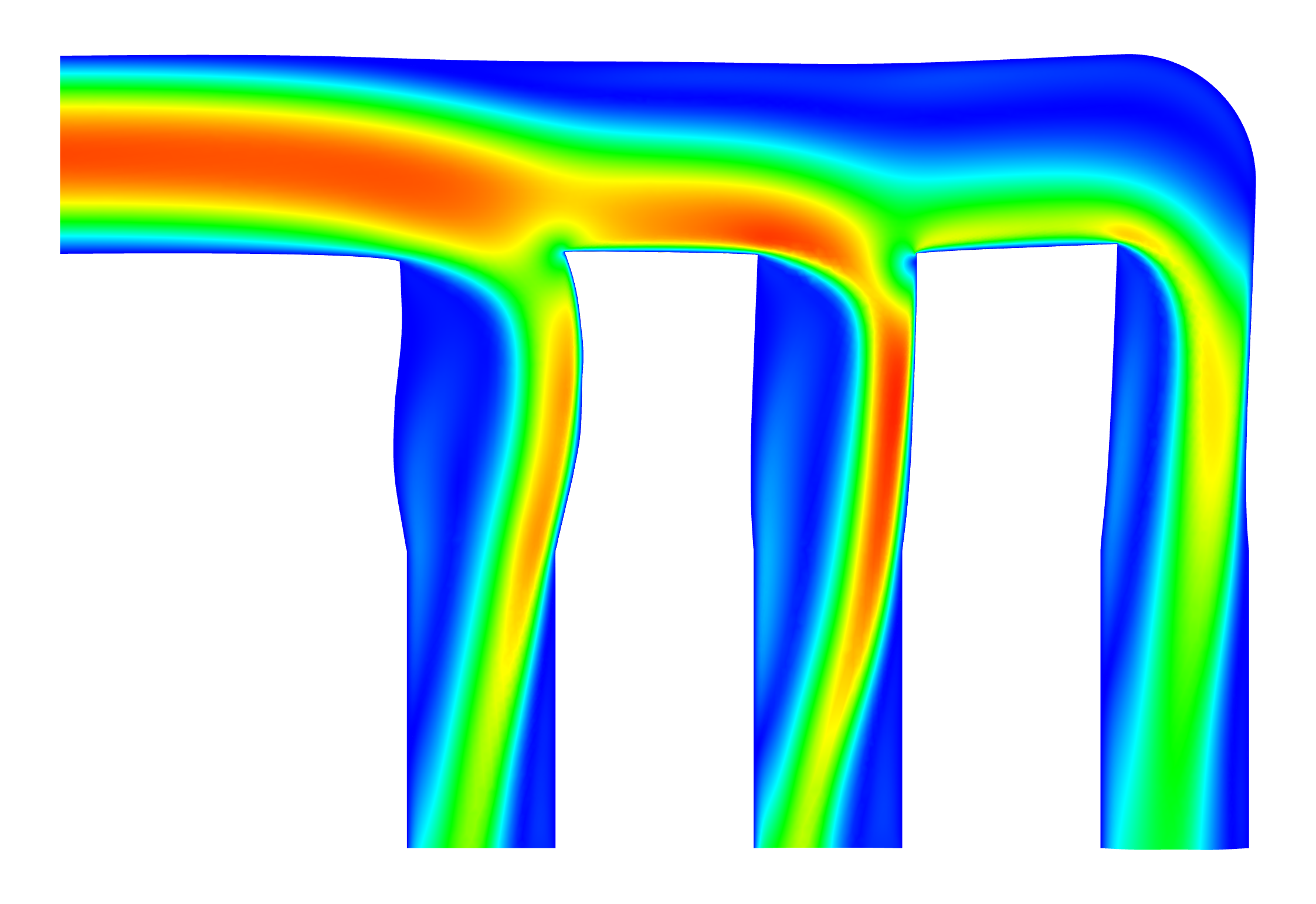}
		\end{minipage}%
		\hfil%
		\begin{minipage}[!t]{0.4\textwidth}
			\includegraphics[width=\textwidth]{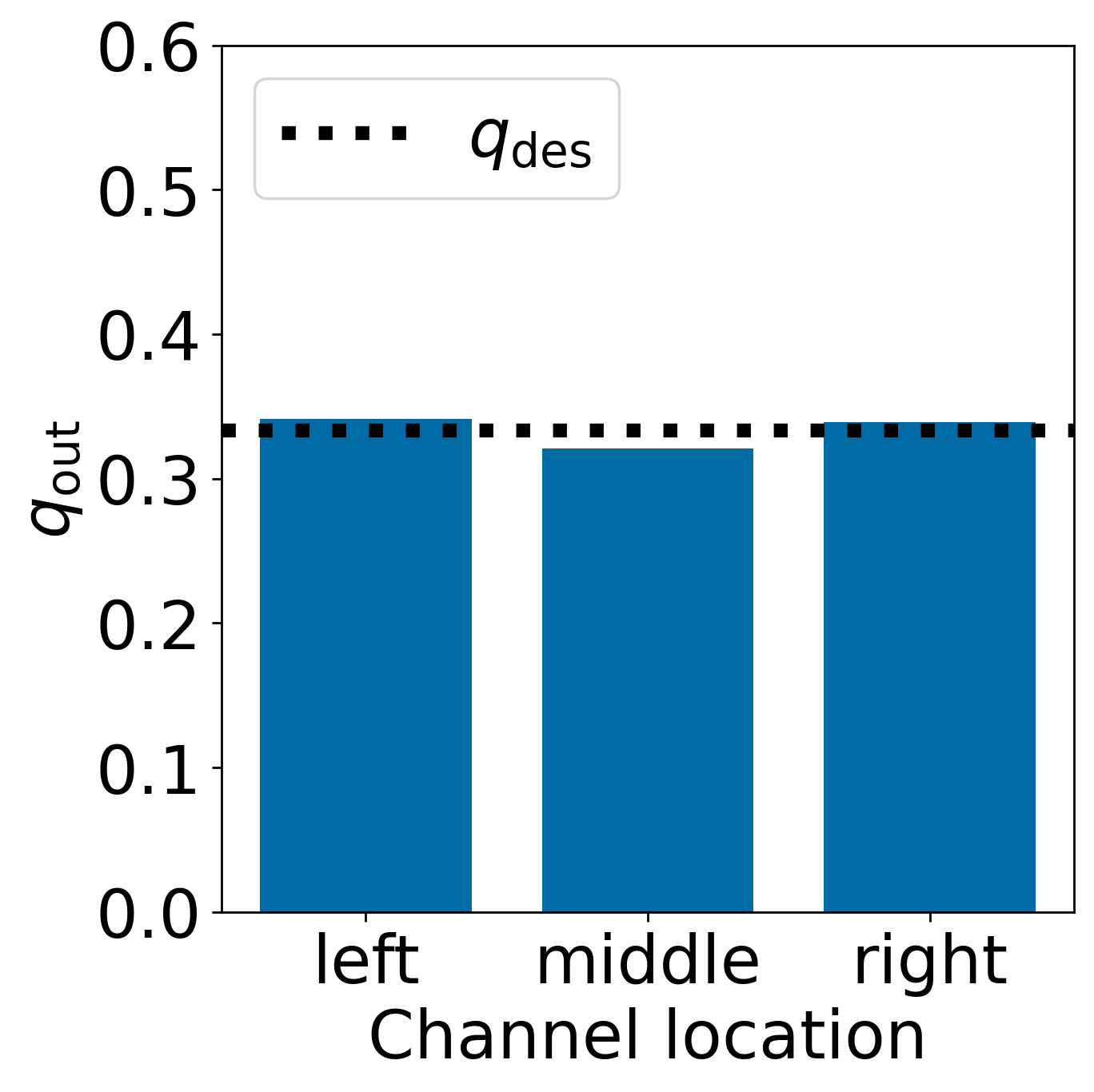}
		\end{minipage}
		\caption{Iteration 2.}
	\end{subfigure}%
	\hfil%
	\begin{subfigure}[!t]{0.49\textwidth}
		\centering
		\begin{minipage}[!t]{0.6\textwidth}
			\includegraphics[width=\textwidth]{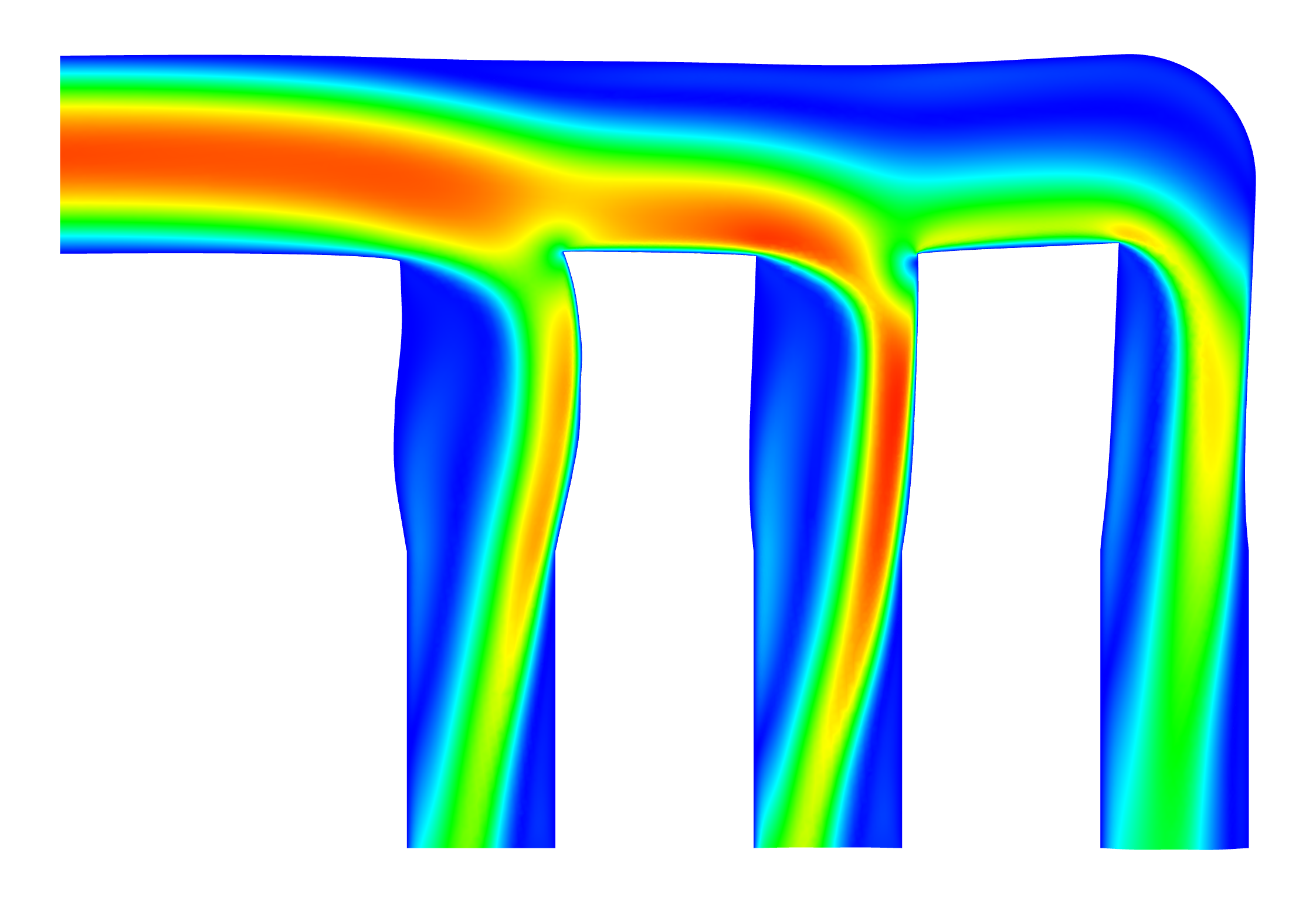}
		\end{minipage}%
		\hfil%
		\begin{minipage}[!t]{0.4\textwidth}
			\includegraphics[width=\textwidth]{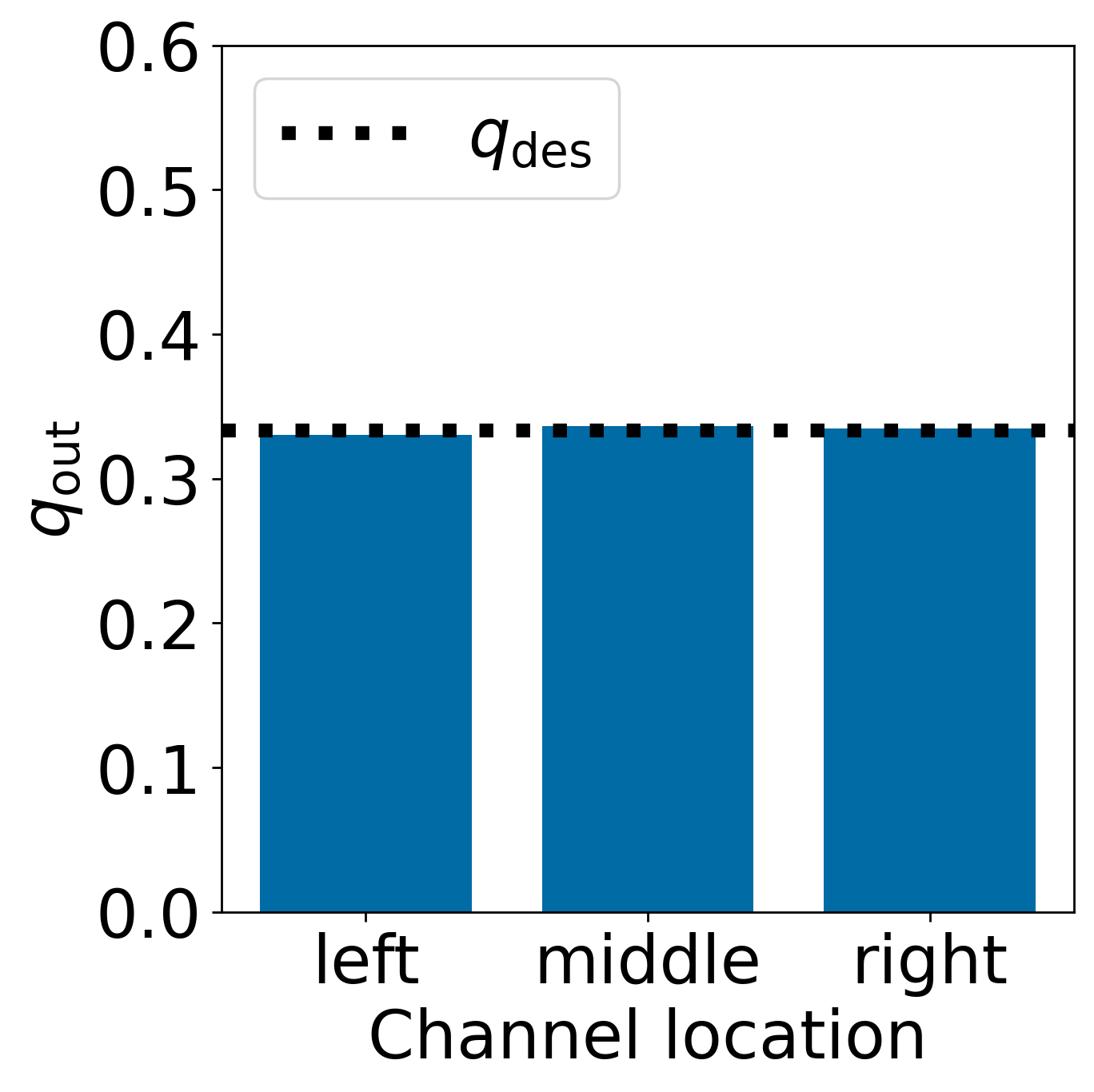}
		\end{minipage}
		\caption{Iteration 3.}
	\end{subfigure}

	\begin{subfigure}[!t]{0.49\textwidth}
		\centering
		\begin{minipage}[!t]{0.6\textwidth}
			\includegraphics[width=\textwidth]{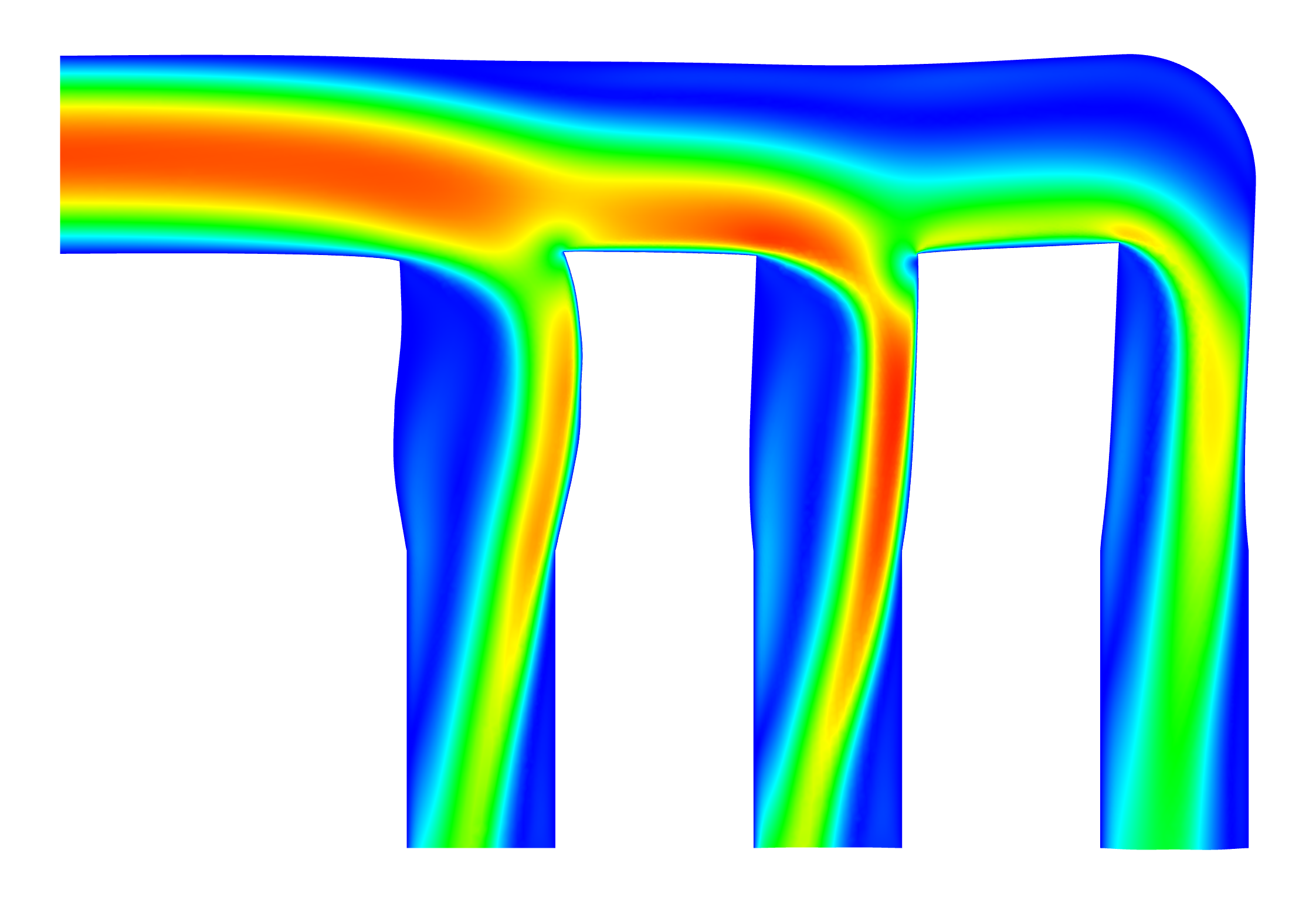}
		\end{minipage}%
		\hfil%
		\begin{minipage}[!t]{0.4\textwidth}
			\includegraphics[width=\textwidth]{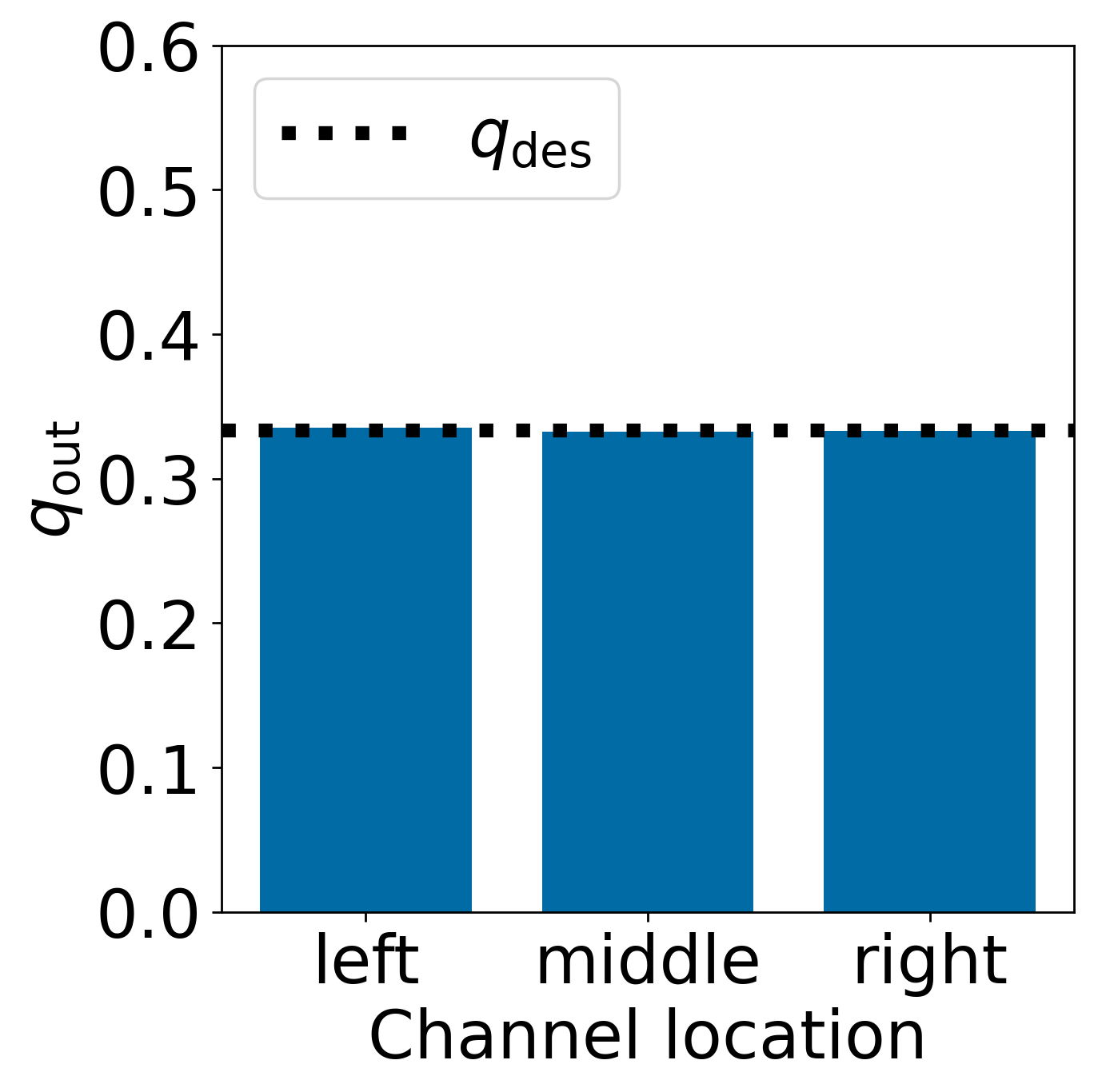}
		\end{minipage}
		\caption{Iteration 4.}
	\end{subfigure}
	\hfil
	\begin{subfigure}[!t]{0.49\textwidth}
		\centering
		\begin{minipage}[!t]{0.6\textwidth}
			\includegraphics[width=\textwidth]{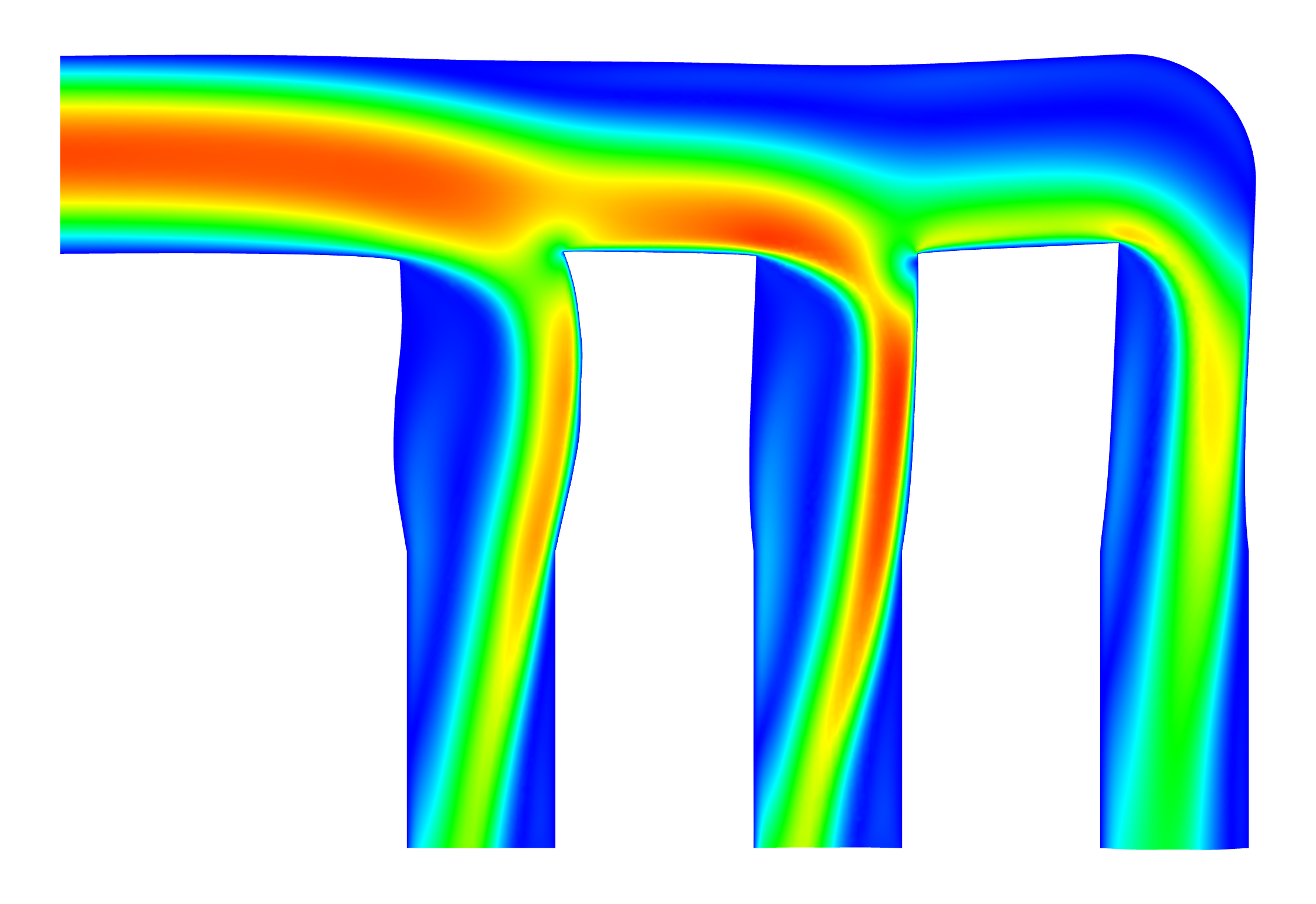}
		\end{minipage}%
		\hfil%
		\begin{minipage}[!t]{0.4\textwidth}
			\includegraphics[width=\textwidth]{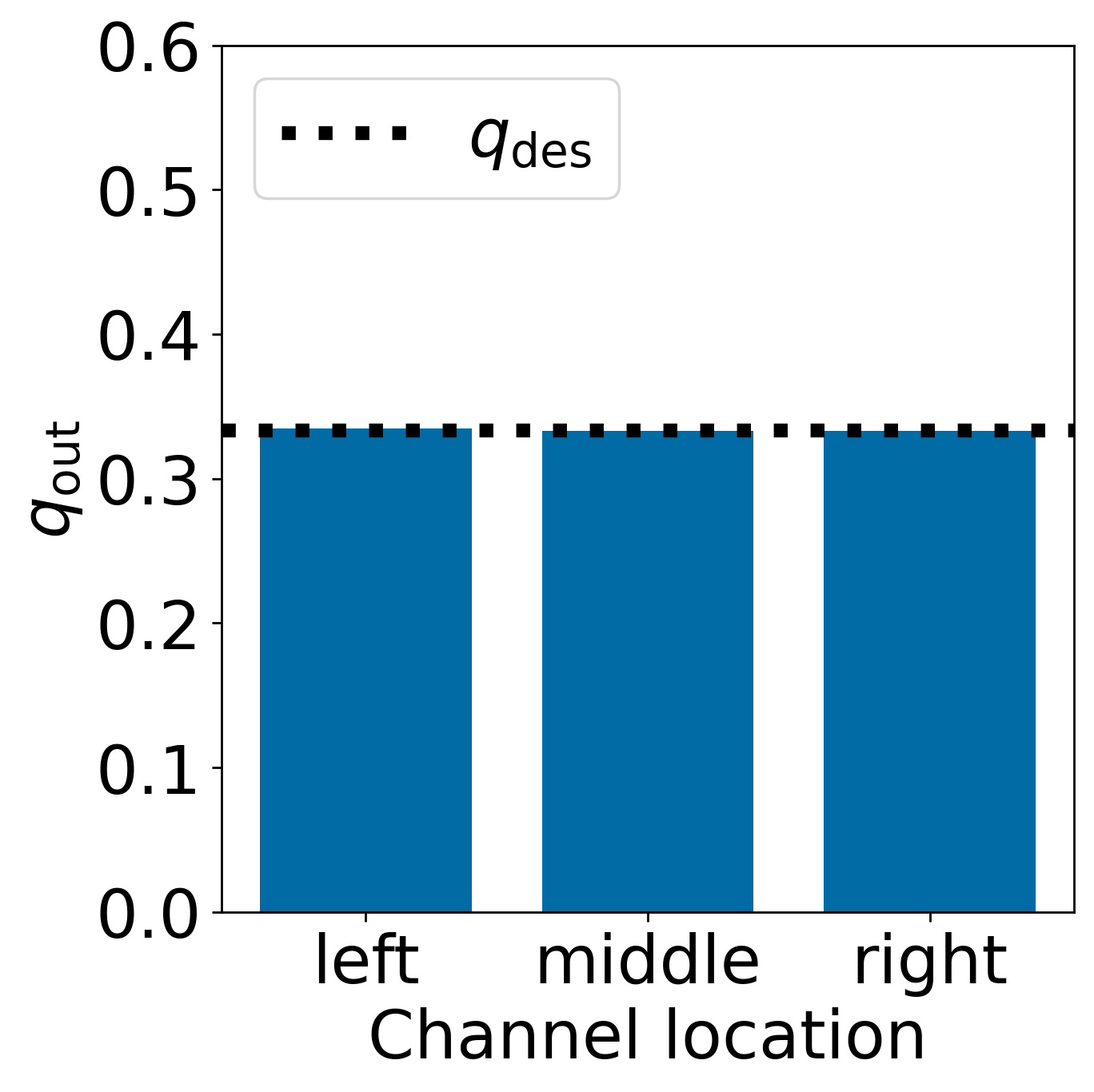}
		\end{minipage}
		\caption{Iteration 5.}
	\end{subfigure}
	
	\caption{Evolution of the geometry $\Omega$, the velocity field $u$, and the outlet flow rates $q\subout^i$ over the course of the ASM method for problem \cref{eq:flow_fine_model} with $\mathrm{Re} = \num{1000}$.}
	\label{fig:fluent_geometries}
\end{figure}

A great benefit of the space mapping technique is that it avoids the direct fine model optimization. In particular, one does not require derivative information for the fine model. This allows us to employ commercial solvers for the simulation of the fine model, which might not be able to compute derivative information or for which automatic differentiation techniques cannot be applied due to their closed-source nature. In such cases, the space mapping method is still able to solve the fine model optimization problem by employing a suitable approximation in form of the coarse model. Our implementation of the space mapping methods in our software package cashocs \cite{Blauth2021cashocs} supports this, as only the coarse model optimization problems are solved by cashocs, whereas the simulation of the fine model can be carried out with arbitrary solvers. 

Here, we demonstrate this powerful feature of the space mapping technique by investigating the problem discussed in \cref{ssec:uniform_flow_distribution}, but now we use the commercial software Ansys\textsuperscript{\textregistered} Fluent, Release 2022 R1 \cite{Ansys2022Ansys} instead of FEniCS for simulating the fine model, which also allows us to choose a higher Reynolds number of $\mathrm{Re} = 1000$ for the fine model. Additionally, we employ the SST $k-\omega$ turbulence model of Ansys\textsuperscript{\textregistered} Fluent for the fine model. Other than that, the numerical setup is identical to the one discussed in \cref{ssec:uniform_flow_distribution}. We solve the fine model optimization problem \cref{eq:flow_fine_model} with the ASM method implemented in our software cashocs, but now use a higher relative tolerance of $\tau = \num{1e-2}$ since the problem is harder to solve due to the higher Reynolds number and correspondingly stronger nonlinearity. Moreover, the differences between the fine and coarse model are significantly larger than for the problem considered in Section~\cref{ssec:uniform_flow_distribution} as the fine model now uses a Reynolds number of \num{1000} and also includes turbulence modeling, which justifies using a higher relative tolerance.

The numerical results of the ASM method can be seen in \cref{fig:asm_fluent}, where the evolution of the cost functional and the stationarity measure over the course of the method are visualized. Again, we observe a steep decrease in the cost functional of over four orders of magnitude in the five space mapping iterations. Additionally, the stationarity measure decreases by over two orders of magnitude before the relative stopping criterion is reached. These observations, again, indicate that the ASM method converged successfully and that it is very efficient even in the context of solving complex and highly nonlinear problems. 

This is further reinforced by the results shown in \cref{fig:fluent_geometries}, where the fine model geometries, velocity fields, and outlet flow rates are shown for each iteration of the ASM method. Starting from an even larger mismatch of the flow distributions in the initial space mapping iteration, the ASM method rapidly reaches a nearly uniform flow distribution already in the third iteration and continues to improve this even further until we obtain a numerically perfect flow distribution after five iterations.

\begin{remark}
	This example showcases the full potential of the space mapping method, as it was able to solve the fine model shape optimization problem by only using simulations of the fine model. If we wanted to solve the fine model optimization problem directly, we would have to implement a turbulence model, derive the corresponding adjoint system and sensitivity information for the Navier-Stokes system with the turbulence model, and implement solvers for both the state and adjoint system. 
	
	In contrast, even without utilizing our optimization software cashocs, implementing a solver for the coarse model shape optimization problem, which is constrained by just the linear Stokes system, is far less challenging, and the derivation of the adjoint system and shape derivatives can even be performed by hand (see, e.g., \cite{Blauth2020Shape}). Therefore, the space mapping approach presented in this paper significantly reduces the complexity of solving the fine model optimization problem and shows great potential for industrial applications.
\end{remark}

\section{Conclusion and Outlook}
\label{sec:conclusion}

In this paper, we have proposed and investigated novel space mapping techniques for PDE constrained shape optimization based on shape calculus. After recalling some basics of space mapping, shape calculus, and the Steklov-Poincar\'e-type metrics from \cite{Schulz2016Efficient}, we presented our framework for the space mapping methods for shape optimization. We formulated an aggressive space mapping (ASM) method and detailed its numerical discretization and implementation in our software cashocs \cite{Blauth2021cashocs}. Finally, we investigated the performance of the ASM method numerically for a shape identification problem constrained by a semi-linear transmission problem and a shape optimization problem of uniform flow distribution constrained by the incompressible Navier-Stokes equations. The corresponding numerical results show that the space mapping technique is very efficient at solving complex shape optimization problems numerically as it only requires a handful of iterations to find an (approximate) optimizer. The proposed methods are particularly well-suited for solving shape optimization problems arising from industrial problems as they can couple closed-source commercial solvers for the fine model with (open-source) software for the optimization of the coarse model, which we have successfully demonstrated in this paper. Altogether, our space mapping methods are efficient and attractive techniques for the numerical solution of shape optimization problems, particularly for large-scale industrial problems.

Future research possibilities include the mathematical analysis of the proposed space mapping methods, e.g., in analogy to the investigations of the classical space mapping method in \cite{Echeverria2005Space}. Moreover, the application of the space mapping methods to industrial problems, e.g., in the context of uniform flow distribution for meltblown processes and for the shape optimization of electrolysis cells, is planned for future work.

\section*{Acknowledgments}

The author thanks Christian Leith\"auser for his help with coupling the space mapping methods with Ansys\textsuperscript{\textregistered} Fluent.

\bibliographystyle{siamplain}
\bibliography{literature_db.bib}

\end{document}